\documentclass[review]{elsarticle}

\usepackage{lineno,hyperref}
\usepackage[fleqn]{amsmath}   %first part makes equations aligned on the left
\usepackage{amsfonts}
\usepackage{amssymb}
\usepackage{multirow}
\usepackage{multicol}
\usepackage{graphicx}
\usepackage{wrapfig}
\usepackage{bm}
\usepackage[font=footnotesize]{caption}
\usepackage{subcaption}
\usepackage[dvipsnames]{xcolor}
\usepackage{scalerel}
\usepackage{xspace}
\usepackage{dsfont}

\newcommand{\x}{\bm{x}}
\newcommand{\X}{\bm{X}}
\newcommand{\n}{\bm{n}}

\newcommand{\U}{\bm{U}}
\newcommand{\grad}{\nabla}
\newcommand{\dotp}{\boldsymbol{\cdot}}

\newcommand{\C}{\mathcal{C}}
\renewcommand{\L}{\mathcal{L}}

\newcommand{\chI}{\textcolor{black}}
\newcommand{\chII}{\textcolor{black}}
\newcommand{\chIII}{\textcolor{black}}
\newcommand{\ch}{\textcolor{black}}
\newcommand{\cha}{\textcolor{black}}
\newcommand{\chb}{\textcolor{black}}
\newcommand{\MATLAB}{\textsc{Matlab}\xspace}

\hypersetup{
  colorlinks   = true, %Colours links instead of ugly boxes
  urlcolor     = blue, %Colour for external hyperlinks
  linkcolor    = blue, %Colour of internal links
  citecolor   = green %Colour of citations
}

\modulolinenumbers[5]

\journal{Elsevier }

\begin{document}

\begin{frontmatter}

\title{\chI{Immersed Boundary Double Layer method: An introduction of methodology on the Helmholtz equation}}
\author[mymainaddress]{Brittany J. Leathers}
\author[mymainaddress]{Robert D. Guy}

\address[mymainaddress]{Department of Mathematics, University of California, Davis, Davis, CA 95616-5270, USA}

\begin{abstract}
The Immersed Boundary (IB) method of Peskin (J. Comput. Phys., 1977) is useful for problems that involve fluid-structure interactions or complex geometries. By making use of a regular Cartesian grid that is independent of the geometry, the IB framework yields a robust numerical scheme that can efficiently handle immersed deformable structures. Additionally, the IB method has been adapted to problems with prescribed motion and other PDEs with given boundary data. IB methods for these problems traditionally involve penalty forces which only approximately satisfy boundary conditions, or they are formulated as constraint problems. In the latter approach, one must find the \ch{unknown} forces by solving an equation that corresponds to a poorly conditioned first-kind integral equation. This operation can therefore require a large number of iterations of a Krylov method, and since a time-dependent problem requires this solve at each step in time, this method can be prohibitively inefficient without preconditioning. \chIII{This work introduces} a new, well-conditioned IB formulation for boundary value problems, \chIII{called} the Immersed Boundary Double Layer (IBDL) method. \chII{In order to lay the groundwork for similar formulations of Stokes and Navier-Stokes equations, this paper focuses on the Poisson and Helmholtz equations to introduce the methodology and }to demonstrate its efficiency over the original constraint method. In this double layer formulation, the equation for the \ch{unknown} boundary distribution corresponds to a well\ch{-}conditioned second-kind integral equation that can be solved efficiently with a small number of iterations of a Krylov method without preconditioning. Furthermore, the iteration count is independent of both the mesh size and spacing of the immersed boundary points.  The method converges away from the boundary, and when combined with a local interpolation, it converges in the entire PDE domain. Additionally, while the original constraint method applies only to Dirichlet problems, the IBDL formulation can also be used for Neumann boundary conditions.
\end{abstract}


\end{frontmatter}

\section{Introduction}\label{introduction}

The Immersed Boundary (IB) method \cite{Peskin77, Peskin02} is a valuable numerical tool for fluid-structure interactions. It was initially developed by Peskin for problems involving elastic, deformable structures, such as those involved in cardiac dynamics \cite{Peskin-72, Peskin-81-heart, mitral-valve, heart-model, PeskinGriffithheart}. However, the robustness and simplicity of the IB method has led to its use in many different applications (see \cite{red-blood-cells, Peskin-Fogelson, pulp}, for just a few). Additionally, since its creation, there have been many developments and variations in the IB method. For instance, the IB framework has been altered to incorporate porous boundaries \cite{porous1,porous2} and model elastic rods with a curvature or twist \cite{gIB1, gIB2}. It has also been adapted to the flow of non-Newtonian fluids \cite{nonNewtonian1,nonNewtonian2} and coupled with internal force mechanisms to model swimming organisms \cite{C-F-C-D,cilia}. In recent years, work has been done to apply the Immersed Boundary framework to problems involving prescribed boundary values \cite{Taira, GriffithDonev, GriffithBhalla, Uhlmann, SuLaiLin, Guylewis}, which is the category of problems that this paper addresses. 

The robustness of the IB method comes from its use of two coordinate systems: a Lagrangian system that moves with the structure and a fixed Eulerian system on which the fluid equations are solved. It uses convolutions with discrete delta functions to link these systems together and map forces from the structure to the grid. One can then solve the PDE on a regular Cartesian mesh. This elimination of the physical boundary makes it possible to use a PDE solver that is efficient and independent of the geometry of the structure.

As stated, this paper will focus on problems involving prescribed boundary values. In the field of fluid dynamics, this includes the motion of rigid bodies and fluid flow through a domain with stationary boundaries or boundaries with prescribed motion. Since the IB method bypasses the need for a conforming mesh, it also has obvious advantages in the broader case of solving PDEs on complex domains. \chII{In this paper, we focus on scalar elliptic boundary value problems to introduce a new IB method.}

When the IB method is applied to deformable structures, the boundary force density is found using a constitutive law from the structure characteristics \cite {Peskin02}, but in the case of rigid bodies, one needs a different way to find or interpret the boundary force. There have been several routes taken to use the IB framework in these situations. One method is to consider the Lagrangian points to be tethered to specified locations by springs. The boundary force is then interpreted as the spring restoring force, penalizing deviations from the prescribed boundary position \cite{Beyer, Goldstein, formal, TeranPeskin}. By using the IB framework, such penalty methods can be efficiently implemented, but they use parameters, such as spring constants, to approximate a rigid limit, and the required magnitude of the appropriate parameter leads to numerically stiff equations that necessitate very small timesteps. 

Another approach for applying the IB method to a rigid body problem consists of viewing the boundary force density as a Lagrange multiplier, used to enforce the no-slip boundary condition \cite{Taira}. Section \ref{IB method} gives a description of such a method as applied to a scalar PDE. In this  IB constraint method, the velocity and force are both unknowns in an algebraic system. One way to solve this sytem is to invert the Schur complement to solve first for the Lagrange multiplier force and then for the velocity. The main disadvantage of such a method is that this operator suffers from poor conditioning \cite{GriffithDonev}. Furthermore, the conditioning of the discrete problem worsens when the Cartesian grid is refined or when the Lagrangian point spacing is refined relative to the grid. To avoid solving the fully discretized constraint problem, most numerical methods using this constraint approach rely on some form of time step splitting \cite{Taira, GriffithBhalla, Uhlmann, SuLaiLin, Guylewis}. Taira and Colonius \cite{Taira}, for example, solve a simpler unconstrained system for an intermediate velocity, use this velocity to find pressure and the unknown boundary force, and then complete a projection step to remove non-divergence-free and slip components of the velocity. However, like the penalty methods, these fractional step methods only satisfy the constraint equations \textit{approximately}, which can result in fluid penetration into a rigid body. Additionally, such a method cannot be used for steady Stokes or other time-independent PDEs. 

However, if instead of a time-splitting scheme, one inverts the Schur complement with a Krylov method, the number of iterations required could be very large due to its poor conditioning. The computational cost then becomes prohibitive for a time-dependent problem, when this solve would be required at each time step. One way around this is to design a preconditioner \cite{Ceniceros, guyphilipgriffith, Stein}, such as the physics-based approximation of the Schur complement constructed by Kallemov et al. \cite{GriffithDonev}. A preconditioner can allow for a more efficient solution to the constraint problem, but developing and implementing one is a nontrivial undertaking, and such a preconditioner can involve computing the inverse of a dense matrix, which can be computationally expensive when the number of boundary points is large. Another limitation is that in order to control the conditioning, there is often a requirement that the boundary points be spaced to about twice the grid spacing, which, depending on the application and the Eulerian-Lagrangian coupling scheme, may result in decreased accuracy of the solution \cite{GriffithDonev, Griffithpointspacing, Griffithpressure}. 

In this paper, we present a reformulated Immersed Boundary method for prescribed boundary values. \chII{In order to lay the groundwork for similar formulations of Stokes and Navier-Stokes equations, in this paper, we focus on the Poisson and Helmholtz equations to introduce the methodology. }Like the IB constraint method, our method enforces the boundary conditions \textit{exactly}, but the resulting linear system is very well-conditioned. We can therefore solve it with an unpreconditioned Krylov method with very few iterations. This conditioning does not worsen as we refine the mesh, nor as we tighten the spacing of our boundary points relative to the grid. We are therefore able to avoid the need for a preconditioner altogether. 

We formulate this new IB method by utilizing a connection between the IB constraint method and a single layer boundary integral equation. Boundary integral methods rely on reformulating a boundary value problem as an integral equation with an unknown density on the boundary. These methods are particularly useful for linear, elliptic, and homogeneous PDEs \cite{Pozblue}. One advantage of these methods, which we exploit, is that there are well-conditioned integral representations available. \chII{For example, for the Helmholtz equation or Stokes equations, we can derive single and double layer integral representations.  When combined with Dirichlet boundary conditions, a double layer representation results in a second-kind integral equation. In this case, the operator we need to invert is well-conditioned. Therefore, these double layer representations have been frequently utilized within boundary integral methods \cite{Mayo1984, ShilpaBI, powermiranda, ingber}.} There are, however, some disadvantages to \chII{directly using boundary integral equations}. Firstly, \chII{for a typical boundary integral method,} once the boundary density is found, it is expensive to \ch{evaluate} the solution on an entire grid, whereas by using efficient solvers, the IB method can do this quickly. Additionally, in order to directly implement integral methods, one needs an analytical Green's function for each specific PDE and problem domain, making nontrivial exterior boundary conditions complicated to implement. The IB method, on the other hand, does not \ch{require} a Green's function and can be used on more general domains. The method we present in this paper maintains the flexibility of the IB method while capturing the better conditioning of a double layer integral equation.

The connection between the IB constraint method and a regularized single layer integral equation has been identified in a few recent works. Usabiaga et al. illustrate this connection and discuss that their rigid multiblob method can be seen as a technique for solving a regularized first-kind integral equation for Stokes flow \cite{GriffithDonev2}. Eldredge \cite{eldredge} made a more general connection between IB methods and boundary integral equations by \ch{extending} the form of a PDE to govern a variable that defines a different function for each side of an immersed boundary. The resulting PDE contains jumps in field quantities across the boundary, which correspond to the strengths of single and double layer potentials. This general connection allows for solutions on either side of the boundary. However, in this work, we look specifically at the IB method for prescribed boundary values for which the PDE domain exists on only one side of the boundary, and we create an IB formulation that corresponds to the use of a double layer potential alone in order to take advantage of the better conditioning. 

\chII{In this paper, we present an introduction to the new \textit{Immersed Boundary Double Layer (IBDL) method}, as it applies to scalar elliptic PDEs. This is a first step towards applying the method to Stokes and Navier-Stokes flows with rigid bodies, which will be presented in a forthcoming work. The method is} able to achieve the same order of accuracy as the original IB constraint method away from the boundary, while only requiring a small number of iterations of a Krylov method. \chI{In particular, as we show in the results section, we see iteration counts that range from $10$ to $2500$ times smaller than those of the IB constraint method, depending on the domain and boundary value problem. Furthermore, unlike those of the IB constraint method, these iteration counts are essentially independent of the grid point spacing.}  We derive, implement, and analyze this method for the 2-D Helmholtz and Poisson problems with Dirichlet boundary conditions\cha{, and we also implement the method for the 3-D Helmholtz problem.} With the drastically improved efficiency, the method has the potential to be utilized in moving domains without the need for stiff parameters or preconditioners. Additionally, while the original constraint method applies only to Dirichlet problems, our new method can also be used for Neumann boundary conditions. 

The paper is organized as follows. In Section \ref{background}, we give the necessary background material on the IB constraint method and boundary integral methods. In Section \ref{connection}, we explicitly demonstrate the connection between the IB constraint method and a single layer integral equation. The main contribution of this paper is contained in Section \ref{ibdl method}, where we formulate the new Immersed Boundary Double Layer method. In Section \ref{numerical implementation}, we discuss the numerical implementation. Section \ref{Results} gives the results for the Helmholtz and Poisson problems with comparisons to the IB constraint method. We also make numerical observations concerning pointwise convergence near the boundary and convergence of the potential strength. Lastly, we end with a discussion in Section \ref{discussion}.

%%%%%%%%%%%%%%%%%%%%%%%%%%%%%%%%%%%%%%%%%%%%%%%%%%%%%%%%%%%%%%%%    BACKGROUND    %%%%%%%%%%%%%%%%%%%%%%%%
%%%%%%%%%%%%%%%%%%%%%%%%%%%%%%%%%%%%%%%%%%%%%%
\section{Background methods}\label{background}

In order to present our method, we will restrict our attention to the 2-D Dirichlet Helmholtz equation, 
\begin{subequations} \label{pde}
\begin{alignat}{2}
& \Delta u - k^2 u = g \qquad && \text{in } \Omega  \label{pde1}\\
&u=U_b \qquad && \text{on } \Gamma,  \label{pde2}
\end{alignat}
\end{subequations}
where $\Omega$ is the domain of our PDE, and $\Gamma = \partial \Omega$. We will explore several interior and exterior problem domains for $\Omega$. 

Sections \ref{IB method} and \ref{BI method} will discuss two types of methods that can be used to solve this PDE, the Immersed Boundary constraint method and boundary integral equations. The two families of methods will be connected in Section \ref{connection} and then used in the formulation of a new scheme in Section \ref{ibdl method}.

%%%%%%%%%%%%%%%%%% IB METHOD %%%%%%%%%%%%%%%%%%%%%%%%
\subsection{Immersed Boundary constraint method}\label{IB method}

\begin{figure}
\centering
\begin{subfigure}{0.495\textwidth}
\centering
\includegraphics[width=\textwidth]{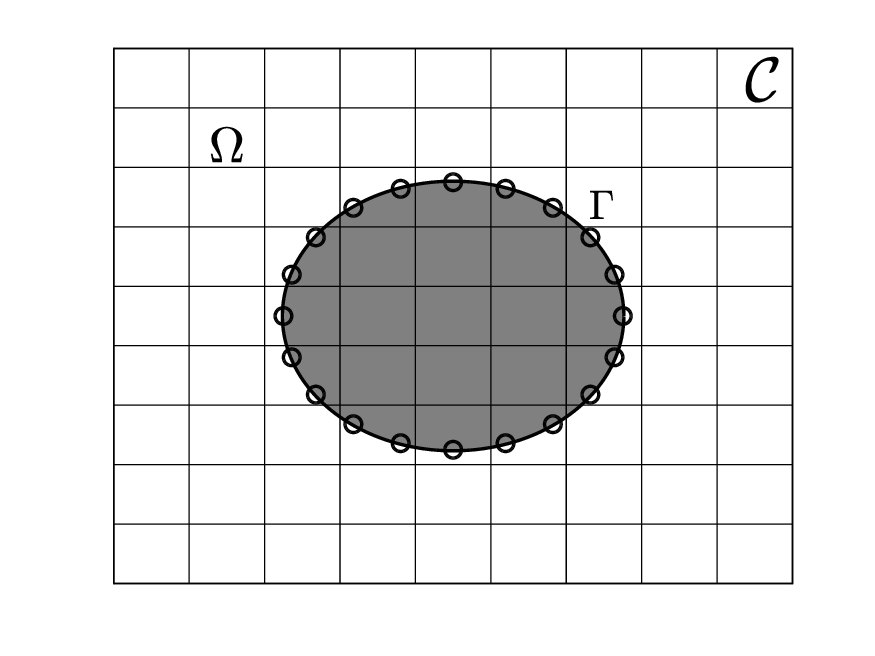}
\caption{\normalsize Exterior Domain}
\label{exterior domain}
\end{subfigure}
\begin{subfigure}{0.495\textwidth}
\centering
\includegraphics[width=\textwidth]{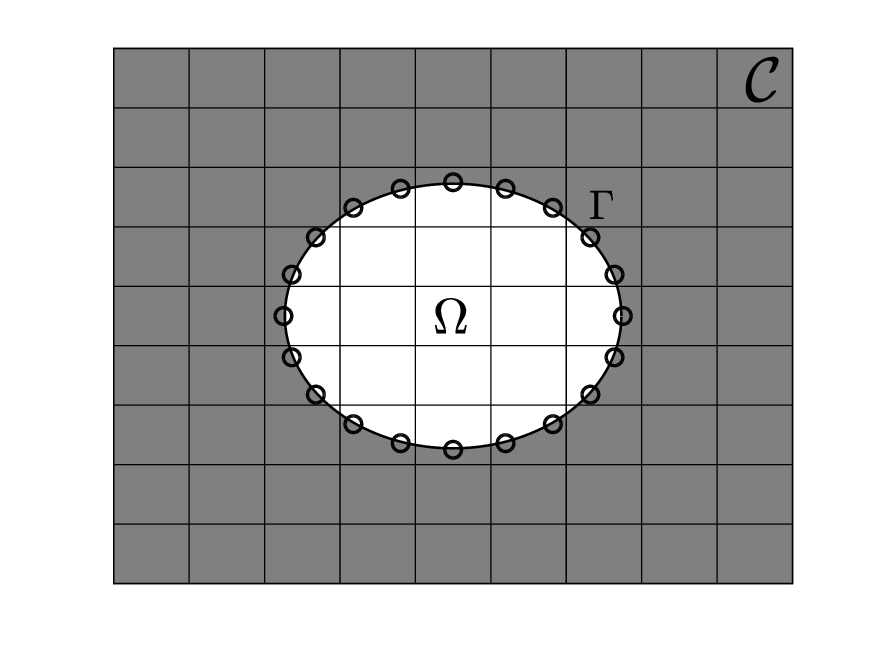}
\caption{\normalsize Interior Domain}
\label{interior domain}
\end{subfigure}
\caption{Diagram showing the two IB coordinate systems for an immersed structure $\Gamma$, embedded in a 2-D computational domain $\C$. The physical domain or PDE domain, given by $\Omega$, is shown in white, and the non-physical domain is shown in grey. \ref{exterior domain} shows an exterior PDE domain, and \ref{interior domain} shows an interior PDE domain.}\label{IB domains}
\end{figure}

One approach for solving Equation \eqref{pde} is the Immersed Boundary constraint method \chI{\cite{Taira, GriffithDonev}}. In this method, the original PDE domain $\Omega$ is embedded into a larger computational domain $\C$. The method then utilizes a Lagrangian coordinate system located on the immersed boundary, as well as an Eulerian coordinate system, in the form of a regular grid. Figure \ref{IB domains} illustrates these coordinate systems with two possible PDE domains.

In the original IB framework, the density of the force exerted by the structure on the fluid is represented by $F$, and it is supported on the boundary $\Gamma$. $F$ is mapped to the Eulerian grid through a convolution with a delta function. The resulting Eulerian force density is incorporated into the PDE, which is then solved in $\C$. With the elimination of the physical boundary, the PDE solve can be done efficiently on a regular grid. The resulting solution $u$ can then be mapped back to the boundary through another convolution with a delta function, and this gives the boundary values. 

In the case of prescribed boundary values, as in Equation \eqref{pde}, the force density, $F$, is a Lagrange multiplier, used to enforce the boundary condition. We will use $s$ and $\x$ as our Lagrangian and Eulerian coordinates, respectively. \ch{Additionally}, $\X(s)$ will be used as a parametrization of the boundary $\Gamma$. We then get the following equations for the continuous IB formulation of Equation \eqref{pde}.
\begin{subequations} \label{constraint IB}
\begin{alignat}{2}
& \Delta\chII{\tilde  u(\x)} - k^2 \chII{\tilde u(\x)} +\int_{\Gamma} F(s)\delta(\x-\X(s))ds = \chII{\tilde g(\x)} \qquad && \text{in } \mathcal{C}  \label{constraint IB1}\\
&\int_{\C} \chII{\tilde u(\x)}\delta(\x-\X(s))d\x=U_b(s) \qquad && \text{on } \Gamma  \label{constraint IB2}
\end{alignat}
\end{subequations}

\chII{In the above formulation, our unknowns are $F$ and $\tilde u$. Note that $\tilde u(\x)$ and $\tilde g (\x)$ represent extensions of the solution $u$ and the function $g$ from $\Omega$ to the larger computational domain $\C$. Once the solution $\tilde u(x)$ is found in $\C$, the solution to the original PDE is then given by $\tilde u|_{\Omega}$. As such, we will drop the tilde notation on $u(x)$. On the other hand, we can define the extension of $g$ by $\tilde g = g\chi_{\scaleto{\Omega}{4.5pt}} + g_e \chi_{\scaleto{\C \setminus \Omega}{6pt}}$, where $\chi$ is an indicator function. Since $\C\setminus\Omega$ is not in the PDE domain, there is flexibility in our choice for $g_e$. For example, \chIII{we} can simply use a smooth extension of $g$ or use $g_e=0$.} 

For mappings between the coordinate systems, the IB method uses a regularized delta function, $\delta_h$, where the regularization lengthscale, $h$, is generally chosen to be on the order of grid point spacing. We can then view the first convolution as \textit{spreading} the force density to the nearby grid points, and the second convolution as \textit{interpolating} the solution values onto the boundary points. Therefore, let us define the spread and interpolation operators, respectively, as
\begin{subequations} \label{operators}
\begin{alignat}{1}
&(SF)(\x) =\int_{\Gamma} F(s)\delta_h(\x-\X(s))ds  \label{spread operator}\\
&(S^*u)(s)=\int_{\C} u(\x)\delta_h(\x-\X(s))d\x.   \label{interpolation operator}
\end{alignat}
\end{subequations}
These operators are adjoint in the sense that 
\begin{equation}\label{adjoint}
\langle SF, u\rangle_{\C} = \langle F, S^*u\rangle_{\Gamma},
\end{equation}
where the inner products are the usual $L^2$ inner products on $\Omega$ and $\Gamma$, respectively. Discretization and numerical implementation will be discussed more in Section \ref{numerical implementation}, but note that we will also use $S$ and $S^*$ to refer to the discretized versions of these operators. Similarly, let $\L$ represent our continuous or discretized Helmholtz operator. We then get the following system of equations. 
\begin{subequations} \label{constraint IB lL}
\begin{alignat}{2}
& \L u +SF = \chII{\tilde g }\qquad && \text{in } \mathcal{C}  \label{constraint IB L1}\\
&S^* u = U_b \qquad && \text{on } \Gamma,  \label{constraint IB L2}
\end{alignat}
\end{subequations}
or
\begin{equation} \label{saddle point equation}
\begin{pmatrix}
\L & S\\
S^* & 0 
\end{pmatrix}
\begin{pmatrix}
u\\
F 
\end{pmatrix}= 
\begin{pmatrix}
\chII{\tilde g}\\U_b
\end{pmatrix}.
\end{equation}

We can solve this saddle point problem by first solving the equation 
\begin{equation}
-(S^* \L^{-1}S)F=U_b -S^* \L^{-1} \chII{\tilde g }\label{schur complement}
\end{equation}
for $F$ by inverting the Schur complement and then obtaining $u$ from Equation \eqref{constraint IB L1}. 

%The stuff about the condition number growth comes from page 27 of the pdf of citation benzisaddlepoint

\begin{figure}
\centering
\begin{subfigure}{0.495\textwidth}
\centering
\includegraphics[width=\textwidth]{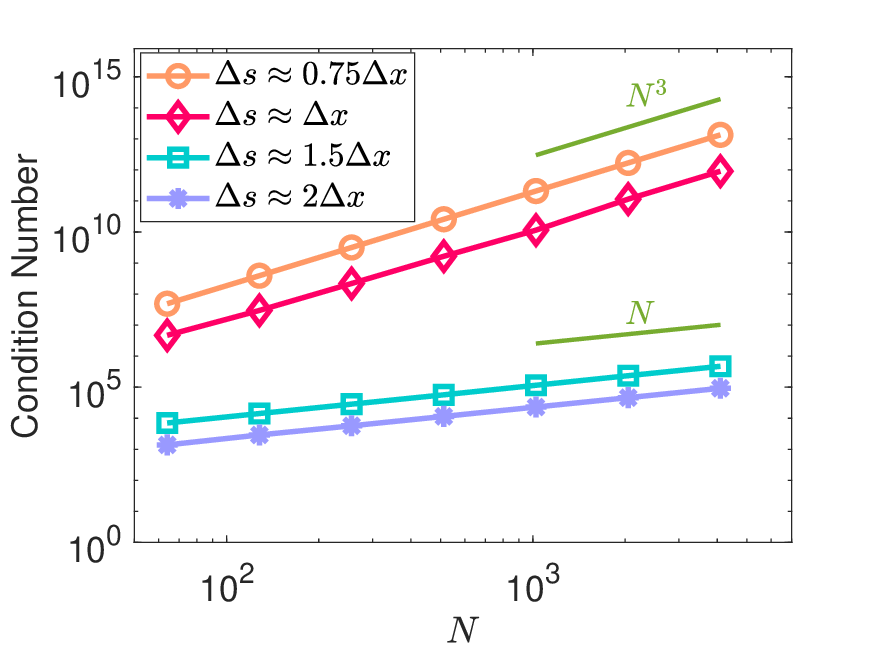}
\caption{\normalsize Condition number growth}
\label{condplot}
\end{subfigure}
\begin{subfigure}{0.495\textwidth}
\centering
\includegraphics[width=\textwidth]{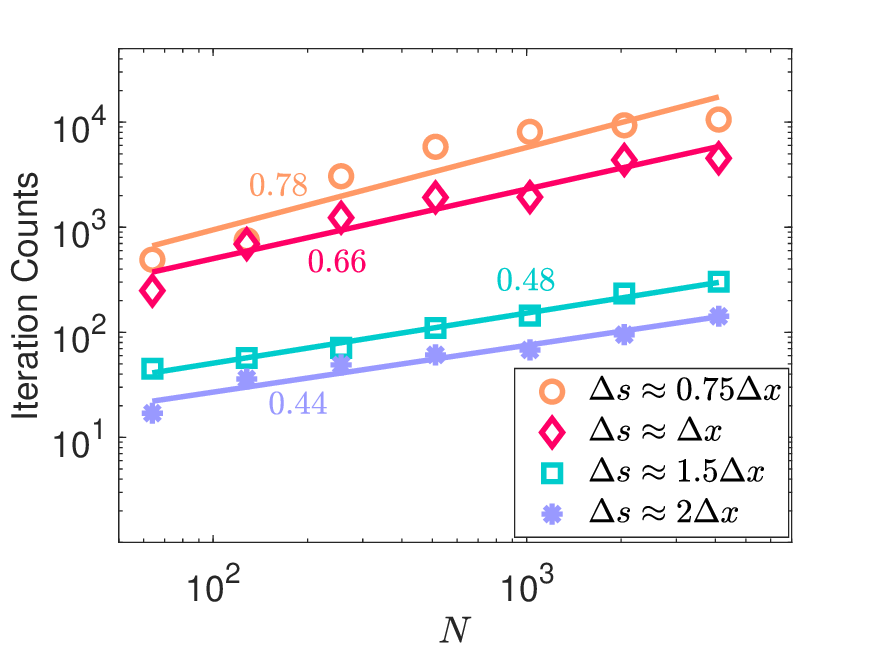}
\caption{\normalsize Iteration count growth}
\label{iteration plot}
\end{subfigure}
\caption{\chI{Condition numbers and iteration counts for the solution of Equation \eqref{schur complement} with $k=1$ and $g=0$ on the periodic computational domain $[-0.5, 0.5]^2$, where the circular boundary of radius 0.25 has prescribed boundary values given by $U_b=\sin{2\theta}$. \ref{condplot} gives the condition number of the Schur complement, given by the left-hand-side of Equation \eqref{schur complement} for various values of $N$, the number of grid points in one direction and for various immersed boundary point spacings. \ref{iteration plot} gives the subsequent number of iterations of \texttt{minres}, with tolerance $10^{-8}$, needed to solve the system, along with the best-fit lines, labeled with their slopes.  }} \label{condandit}
\end{figure}

However, \chI{generically} these saddle point problems resulting from Dirichlet boundary conditions are poorly conditioned \cite{Benzisaddlepoint}. \chI{For example, for the 3-D Stokes equation with a spherical immersed boundary, Kallemov et al.\ found that the condition number grows linearly in the number of immersed boundary points, for $\Delta s \geq  \Delta x$ \cite{GriffithDonev}. Furthermore, as $\Delta s$ is decreased, the increase in condition number is rapid, resulting in an inability to solve the system in double-precision arithmetic for $\Delta s \approx 0.5 \Delta x$ \cite{GriffithDonev}.} \chI{Figure \ref{condandit} gives the condition numbers for the Schur complement and the subsequent iteration counts of \texttt{minres}\footnote{\chI{We utilize \texttt{minres} from \MATLAB R2022a \cite{Matlab22a}.}} needed for solving Equation \eqref{schur complement} without preconditioning, where  $k=1$ and $\Gamma$ is a circle. As will be discussed in Section \ref{numerical implementation}, our computational domain is a periodic box, and we use a finite difference method for discretizing the Helmholtz operator. This PDE will be revisited in Section \ref{interior circle}. We can see that for wider boundary point spacing, the condition number grows as $\mathcal{O}(N)$, where $N$ is the number of grid points in one direction. For tighter boundary point spacings, we see the condition numbers grow even more rapidly, as $\mathcal{O}(N^3)$. The subsequent iteration counts for this problem grow between $\mathcal{O}(\sqrt{N})$ and $\mathcal{O}(N)$.} With an eye toward time-dependent PDEs in applications such as fluid dynamics, where this system would be solved in each time step, these iteration counts can be prohibitively large without proper preconditioning. Therefore, as mentioned in the introduction, much work has been devoted to forming preconditioners \cite{GriffithDonev, Ceniceros, guyphilipgriffith, Stein}. \chIII{However}, preconditioning itself can be computationally expensive, as it often involves inverting a dense matrix. In Section \ref{ibdl method}, we present an alternative IB formulation that achieves the same order of accuracy as this constraint method while avoiding these high iteration counts and therefore eliminating the need for preconditioning.\\

%\begingroup   
%\begin{center}
  %\renewcommand*\arraystretch{0.8}
 %\begin{tabular}{||c | c | c | c ||} 
% \hline
% \multicolumn{4}{|c|}{Iteration Counts - Elliptical Boundary} \\
% \hline
 % $\Delta x$ &$\Delta s = 2\Delta x $& $\Delta s = 1.5 \Delta x$ &$ \Delta s = 1 \Delta x$  \\ [0.5ex] 
% \hline
%$2^{-6}$ &   31  &   24    &        \\
%$2^{-7} $&   61   &  71     &        \\
%$2^{-8 }$&   102  &   413    &           \\
%$2^{-9}$&   145   &    579   &            \\
%$ 2^{-10}$& 202   & 1202      &            \\
%$ 2^{-11}$ &  256     &  1617     &           \\
%$ 2^{-12}$  &  315  &  1909     &           \\
 % \hline
 % \end{tabular}
%\end{center}
%\captionof{table}{Number of iterations of \texttt{minres}, with tolerance $10^{-8}$, needed to solve Equation \eqref{pde} with $k=1$ on a periodic computational domain of length $L=1$, where the boundary is given by $x(t)=\frac14 \cos{\theta} , y(t)= \frac{5}{36} \sin{\theta}$, and the prescribed boundary values are given by %$U_b=\sin{2\theta}$. Boundary points were chosen equally spaced in the parameter $\theta \in[0,2\pi]$.  }
%\label{iteration table 2}
%\endgroup
%\vspace{0.2cm}

%%%%%%%%%%%%%%%%%% INTEGRAL METHODs %%%%%%%%%%%%%%%%%%%%%%
\subsection{Integral methods}\label{BI method}

Boundary integral methods can also be used to solve Equation \eqref{pde} by reformulating the boundary value problem as an integral equation, using the appropriate Green's function \cite{Pozblue}. The Green's function for this PDE satisfies\begin{equation}
\Delta G(\x,\x_0)-k^2G(\x,\x_0) = - \delta (\x-\x_0) \hspace{0.6cm} \text{in } \Omega.  \label{greens function}
\end{equation}
For convenience, let $g=0$, and let our problem domain $\Omega$ be an interior domain, with boundary $\Gamma$. Similar steps can be taken for an exterior domain.

First, Green's second identity simplifies for this PDE to 
\begin{equation}
u(\x) \delta(\x-\x_0)= \grad \dotp \Big(G(\x, \x_0) \grad u(\x) - u(\x) \grad G(\x,\x_0)\Big). \label{greens second identity}
\end{equation} 
Let $\x_0 \in \cha{\Omega\setminus \Gamma} $. Then by integrating \chIII{Equation} \eqref{greens second identity} over $\Omega$ and using the divergence theorem, we get the integral representation of $u$, 
\begin{equation}
u(\x_0) = \int_{\Gamma} G(\x, \x_0) \grad u(\x) \dotp \n (\x) dl(\x) \\- \int_{\Gamma} u(\x) \grad G(\x,\x_0) \dotp \n(\x) dl(\x), \label{integral rep}
\end{equation} 
where $\n$ is the unit normal of $\Gamma$, pointing out of $\Omega$. We use $dl(\x)$ to denote an integral with respect to arclength. 

The first term of Equation \eqref{integral rep} is called a \emph{single layer potential}, which can be viewed as a distribution of point sources on the boundary, and the second term is called a \emph{double layer potential}, which can be viewed as a distribution of point source dipoles on the boundary. Notice that the strength of the single layer potential is given by the boundary distribution of the normal derivative of $u$, and the strength of the double layer potential is given by the boundary values of $u$. 

There are, however, other integral representations of the solution $u$ in $\Omega$. We can, for instance, represent $u$ solely with a single layer potential or a double layer potential. Let us consider two solutions: $u^{Int}$ that satisfies the PDE on $\Omega$ and $u^{Ext}$ that satisfies the PDE on $\Omega^{Ext}$, the region exterior to $\Gamma$. Let these two solutions share the same values on $\Gamma$. In order to find a single layer representation for $u^{Int}(\x_0)$ for $\x_0\in\Omega$ we need two equations. Firstly, we use Equation \eqref{integral rep} for $u^{Int}$. Secondly, by the same process that gave us Equation \eqref{integral rep}, we can get a similar equation for $u^{Ext}$, where the corresponding normal vector points outward from $\Omega^{Ext}$, and the left-hand side vanishes because $\x_0\notin \Omega^{Ext}$. By negating the expression in order to use the same unit normal as that in Equation \eqref{integral rep}, we get the following equation for $u^{Ext}$:   
\begin{equation}
0= -\int_{\Gamma} G(\x, \x_0) \grad u^{Ext}(\x) \dotp \n(\x) dl(\x) \\+ \int_{\Gamma} u^{Ext}(\x) \grad G(\x,\x_0) \dotp \n(\x) dl(\x). \label{integ rep 2}
\end{equation} 
Then, by adding Equation \eqref{integral rep}, for $u^{Int}$, and Equation \eqref{integ rep 2},  we get the \emph{generalized single layer integral representation} for $u^{Int}$ given by 
\begin{equation}
u^{Int} (\x_0) = \int_{\Gamma}  \sigma (\x)  G(\x, \x_0)dl(\x)  ,        \label{single layer rep}
\end{equation}
where $\sigma = \grad(u^{Int}-u^{Ext}) \dotp \n \Big|_{\Gamma}$ is the strength of the single layer potential. 

If we instead consider two such solutions which match normal derivatives on the boundary, we get the \emph{generalized double layer integral representation} for $u^{Int}$ given by 
\begin{equation}
u^{Int} (\x_0) = \int_{\Gamma} \gamma(\x) \grad G(\x, \x_0)\dotp \n(\x) dl(\x)  ,        \label{double layer rep}
\end{equation}
where $\gamma = (u^{Ext}-u^{Int})\Big|_{\Gamma}$ is the strength of the double layer potential. 

In order to use our boundary data, we then take $\x_0 \longrightarrow \Gamma$. To have consistent notation across methods, we will denote a boundary point with a capital letter, $\X_0$. Then, to generalize for exterior domains, we emphasize again that $\n$ is the unit normal pointing out of $\Omega$, regardless of whether $\Omega$ is an exterior or interior domain. Using our boundary condition, Equations \eqref{single layer rep} and \eqref{double layer rep} give us 
\begin{equation}
U_b(\X_0) = \int_{\Gamma}  \sigma (\x)  G(\x, \X_0)dl(\x)   ,     \label{single layer eqn}
\end{equation}
\vspace{-0.3cm}
\begin{equation}
U_b (\X_0) = \int_{\Gamma}^{PV} \gamma(\x) \grad G(\x, \X_0)\dotp \n(\x) dl(\x) - \frac12 \gamma(\X_0),   \label{double layer eqn}
\end{equation}
where $PV$ denotes the principal-value integral, for which $\X_0$ is placed exactly on $\Gamma$. It is computed by integrating over $\Gamma \setminus B_{\epsilon}(\X_0)$, in which a disk around $\X_0$ has been removed from the boundary. The $PV$ integral then comes from taking the limit as $\epsilon$ goes to 0. The $1/2$ in Equation \eqref{double layer eqn} is obtained using the properties of the normal derivative of the Green's function. For a reference, see Section 2.4 in \cite{Pozblue}.

The \ch{unknown} quantities in Equations \eqref{single layer eqn} and \eqref{double layer eqn} are the potential strengths, $\sigma$ and $\gamma$, respectively. Equation \eqref{single layer eqn} has the form of a Fredholm integral of the first kind, and Equation \eqref{double layer eqn} has the form of a Fredholm integral of the second kind \cite{atkinson}. The operators given by  
\begin{subequations} \label{integral operators}
\begin{alignat}{1}
&K_1 \sigma = \int_{\Gamma}  \sigma (\x)  G(\x, \x_0)dl(\x)    \label{sl integral operator}\\
& K_2\gamma = \int_{\Gamma} \gamma(\x) \grad G(\x, \x_0)\dotp \n(\x) dl(\x)    \label{dl integral operator}
\end{alignat}
\end{subequations}
have eigenvalues in the interval $(-0.5, 0.5)$, and the only limit point of the eigenvalues is $0$. \cite{coltonkress, spectralproperties}. The condition number, defined as the ratio of the largest and smallest eigenvalues, is therefore \ch{infinite} in the continuous case and large in the discretized case. On the other hand, Equation \eqref{double layer eqn} can be rewritten as
\begin{equation}
\Big(K_2 - \frac12 I\Big)\gamma = U_b, 
\end{equation}
where $I$ is the identity operator. Shifting the operator in this way shifts the eigenvalues, and the only limit point for the eigenvalues becomes -1/2. Therefore, the condition number of the discretized operator is finite and does not grow with refinement of the discretization. The better conditioning of this operator is the characteristic of the double layer representation that we exploit to form our new Immersed Boundary method.

%%%%%%%%%%%%%%%%%%%%%%%%%%%%%%%%%%%%%%%%%%%%%%%%%%%%%%%%%%%%%%%%    CONNECTION    %%%%%%%%%%%%%%%%%%%%%%%%
%%%%%%%%%%%%%%%%%%%%%%%%%%%%%%%%%%%%%%%%%%%%%%
\section{Connection between IB method and integral methods}\label{connection}
In this section, we relate the Immersed Boundary constraint method as described in Section \ref{IB method} to the single layer integral equation in Equation \eqref{single layer eqn}. As discussed in Section \ref{introduction}, this connection has been established in several recent works \cite{GriffithDonev2, eldredge}, but here, we present this connection explicitly from the IB constraint system of Equation \eqref{constraint IB lL}. After establishing this connection between the two types of methods, we propose that a solution to the conditioning problem of the Schur complement seen in Equation \eqref{schur complement} is to formulate an Immersed Boundary version of a double layer integral equation. Such a formulation will be presented in Section \ref{ibdl method}. 

For Sections \ref{connection} and \ref{ibdl method}, we will simplify the presentation by focusing on the homogeneous case, or $g=0$. However, this connection holds for the inhomogeneous case as well, and we use the new method on several \ch{inhomogeneous} problems in Section \ref{Results}.

Let us define a regularized Green's function by 
\begin{equation}
G_h(\x,\x_0) \equiv G * \delta_h = \int_{\C} G(\mathbf{y}, \x_0)\delta_h(\x-\mathbf{y})d\mathbf{y}. \label{Gh}
\end{equation}
The linearity of our operator then gives us that 
\begin{equation}
\L G_h(\x,\x_0) \cha{ = \L(G * \delta_h) = G * \L\delta_h = \L G * \delta_h = -\delta * \delta_h }= -\delta_h(\x-\x_0).    \label{Gh2}
\end{equation}

The homogeneous IB constraint equations are given by 
\begin{subequations} \label{constraint IB lL again }
\begin{alignat}{2}
& \L u +SF = \chI{0} \qquad && \text{in } \mathcal{C}  \label{constraint IB L1 again}\\
&S^* u = U_b \qquad && \text{on } \Gamma . \label{constraint IB L2 again}
\end{alignat}
\end{subequations}
Starting with the first equation of the IB constraint method, Equation \eqref{constraint IB L1 again}, we have 
\begin{equation}
\L u = - SF = - \int_{\Gamma} F(s)\delta_h(\x-\X(s)) ds. \label{connection1}
\end{equation}
Inverting the operator and using Equation \eqref{Gh2}, we get
\begin{equation}
u(\x) =  \int_{\Gamma} F(s)G_h(\x,\X(s)) ds \label{connection2}
\end{equation}
for $\x\in\Omega\setminus \Gamma$. The second equation of the IB constraint method, Equation \eqref{constraint IB L2 again}, gives us 
\begin{equation}
U_b(s') =\int_{\C} u(\x)\delta_h(\x-\X(s')) dx,  \label{connection3}
\end{equation}
where we are using $s'$ to distinguish from our previous variable of integration or alternately, viewing $\X(s')$ as selecting a particular boundary point. Combining this with Equation \eqref{connection2} and changing the order of integration, we get
\begin{equation}
U_b(s') =\int_{\Gamma} F(s) \int_{\C} G_h(\x,\X(s)) \delta_h(\x-\X(s')) d\x ds. \label{connection4}
\end{equation}
Notice that the second integral is equivalent to $G_h * \delta_h$, evaluated at $\X(s')$. Therefore, by denoting a twice-regularized Green's function by $G_{hh}$, this becomes
\begin{equation}
U_b(s') =\int_{\Gamma} F(s)  G_{hh}(\X(s'),\X(s)) ds.\label{connection5}
\end{equation}
The symmetry of the Green's function, which is preserved through convolutions with the regularized delta function, gives us that $G_{hh}(\X(s'),\X(s)) = G_{hh}(\X(s),\X(s'))$. Using this to switch the arguments of $G_{hh}$ and appropriately redefining $F$ and $U_b$ as functions of $\x$, we get 
\begin{equation}
U_b(\X(s')) =\int_{\Gamma} F(\X(s)) G_{hh}(\X(s),\X(s')) ds. \label{connection6}
\end{equation}

Recall that the single layer integral equation for $\X_0$ on the boundary is given by 
\begin{equation}
 U_b(\X_0) =\int_{\Gamma} \sigma(\x)  G(\x, \X_0) dl(\x).  \label{connection7}
\end{equation}
 By parametrizing $\Gamma$ by $\X(s)$, this becomes 
 \begin{equation}
 U_b(\X(s')) =\int_{\Gamma} \sigma(\X(s))  G(\X(s),\X(s')) \bigg|\frac{\partial \X(s)}{\partial s}\bigg| ds.  \label{connection8}
\end{equation}
By comparing Equations \eqref{connection6} and \eqref{connection8}, we can associate the IB constraint force $F$ with the potential strength modified by the parametrization term, $\big|\partial \X(s)/\partial s\big| $. In other words, in the limit that the regularization width, $h$, approaches $0$, we have
\begin{equation}
F(\X)=\sigma(\X)\bigg|\frac{\partial \X}{\partial s}\bigg|. \label{F eqn}
\end{equation}
Specifically, in the case of an arclength parametrization, where $ \big|\partial \X(s)/\partial s\big| =1$, we see from Section \ref{BI method} that $F$ gives the jump in the normal derivative of the solution across the boundary. We have now established that the IB constraint method is equivalent to a \emph{regularized} single layer integral equation. As such, in order to distinguish it from the new method we propose in Section \ref{ibdl method}, we will henceforth refer to this as the Immersed Boundary Single Layer (IBSL) method.

%%%%%%%%%%%%%%%%%%%%%%%%%%%%%%%%%%%%%%%%%%%%%%%%%%%%%%%%%%%%%%%%    FORMULATION OF IBDL    %%%%%%%%%%%%%%%%%%%%%%%%
%%%%%%%%%%%%%%%%%%%%%%%%%%%%%%%%%%%%%%%%%%%%%%
\section{The Immersed Boundary Double Layer (IBDL) method}\label{ibdl method}

\subsection{Formulation of the IBDL method} \label{formulation}

The discussion at the end of Section \ref{BI method} suggests that to improve the conditioning of the IBSL method, we could use an Immersed Boundary framework to instead solve a regularized \emph{double} layer integral equation, which has better conditioning. We will call this the Immersed Boundary Double Layer (IBDL) method, and it is given by 
\begin{subequations} \label{ibdl}
\begin{alignat}{2}
& \L u +\widetilde S Q = \chII{\tilde g} \qquad && \text{in } \mathcal{C}  \label{ibdl 1}\\
&S^* u + \frac12 Q = U_b \qquad && \text{on } \Gamma,   \label{ibdl 2}
\end{alignat}
\end{subequations}
where we assume an arclength parametrization of $\Gamma$, and where $\widetilde S Q \equiv \grad\dotp (SQ\n)$ for $\n$, the unit normal \chIII{on $\Gamma$}, pointing out of $\Omega$. The form of $\widetilde S$ implies that $Q$ gives the strength of a dipole force distribution on the boundary, and this is our Lagrange multiplier used to enforce the boundary condition. We will next demonstrate that this formulation corresponds to a double layer integral equation. 

Starting with the Equation \eqref{ibdl 1}, \chIII{and taking $g=0$ for simplicity,} we have 
\begin{equation}
\L u = - \widetilde SQ = - \grad \dotp \int_{\Gamma} Q(s)\n(s) \delta_h(\x-\X(s)) ds. \label{ibdl connect1}
\end{equation}
Inverting the operator and using Equation \eqref{Gh2}, we get
\begin{equation}
u(\x) = \grad \dotp  \int_{\Gamma} Q(s)\n(s)G_h(\x,\X(s)) ds \label{ibdl connect2}
\end{equation}
for $\x\in\Omega\setminus \Gamma$. Bringing in the divergence and manipulating the expression, we get
\begin{equation}
u(\x) =   \int_{\Gamma} Q(s) \grad G_h(\x,\X(s)) \dotp\n(s) ds .\label{ibdl connect3}
\end{equation}

Then, the second equation of the IBDL method, Equation \eqref{ibdl 2}, gives us 
\begin{equation}
U_b(s') =\int_{\C} u(\x)\delta_h(\x-\X(s')) dx +\frac12Q(s'). \label{ibdl connect4}
\end{equation}
Combining Equation \eqref{ibdl connect3} with Equation \eqref{ibdl connect4}, changing the order of integration, and again recognizing the presence of $G_{hh} = G_h*\delta_h$, we get 
\begin{equation}
U_b(s') =\int_{\Gamma} Q(s) \grad G_{hh}(\X(s'),\X(s)) \dotp\n(s) ds +\frac12Q(s'). \label{ibdl connect5}
\end{equation}

The odd symmetry of the gradient of the Green's function \cite{Pozblue}, which is again preserved through the convolutions with the regularized delta function, gives us that $\grad G_{hh}(\X(s'),\X(s)) = -\grad G_{hh}(\X(s),\X(s'))$.  Using this to switch the arguments of $\grad G_{hh}$ and appropriately redefining $Q$ and $U_b$ as functions of $\x$, we get
\begin{equation}
U_b(\X(s')) =-\int_{\Gamma} Q(\X(s)) \grad G_{hh}(\X(s),\X(s'))\dotp \n(\X(s)) ds + \frac12Q(\X(s')). \label{ibdl connect5}
\end{equation}

Recall that the double layer integral equation for $\X_0$ on the boundary is given by 
\begin{equation}
 U_b(\X_0) =\int_{\Gamma}^{PV} \gamma(\x) \grad G(\x, \X_0) \dotp \n(\x) dl(\x) - \frac12 \gamma(\X_0).  \label{ibdlconnect6}
\end{equation}
Using our arclength parametrization, we can rewrite this as
\begin{equation}
 U_b(\X(s')) =\int_{\Gamma}^{PV} \gamma(\X(s)) \grad G(\X(s), \X(s')) \dotp \n(\X(s)) ds - \frac12 \gamma(\X(s)).  \label{ibdlconnect7}
\end{equation}
By comparing Equations \eqref{ibdl connect5} and \eqref{ibdlconnect7}, we see that in the limit that the regularization width, $h$, approaches $0$, we have 
\begin{equation}
Q(\X)= - \gamma(\X).
\end{equation}
Therefore, we can see that the IBDL method is equivalent to a \emph{regularized} double layer integral equation. Note from Section \ref{BI method} that $Q$ gives the jump in the solution across the boundary. 

We reiterate that the formulation, specifically Equation \eqref{ibdl 2}, assumes that $X(s)$ gives an arclength parametrization of $\Gamma$. In the IBSL method, if \chIII{we use} a different paramterization of the boundary, the rescaling of $F$ seen in Equation \eqref{F eqn} accounts for the transformation, and no alteration in the method is needed to find the solution. However, in the IBDL method, since the $1/2$ in Equation \eqref{ibdl 2} is derived using an arclength parametrization, one would need to alter this term in the case that  $  \big|\partial \X(s)/\partial s\big| \neq 1$. \chIII{This will be discussed further in Section \ref{arbparam}, but otherwise, we will assume an arclength parametrization throughout this paper.}

%%%%%%%%%%%%%%%%%% DISCONTINUITY %%%%%%%%%%%%%%%%%%%%%%
\subsection{Discontinuity of the solution}\label{discontinuity}
\begin{figure}
\centering
\begin{subfigure}{0.495\textwidth}
\centering
\includegraphics[width=\textwidth]{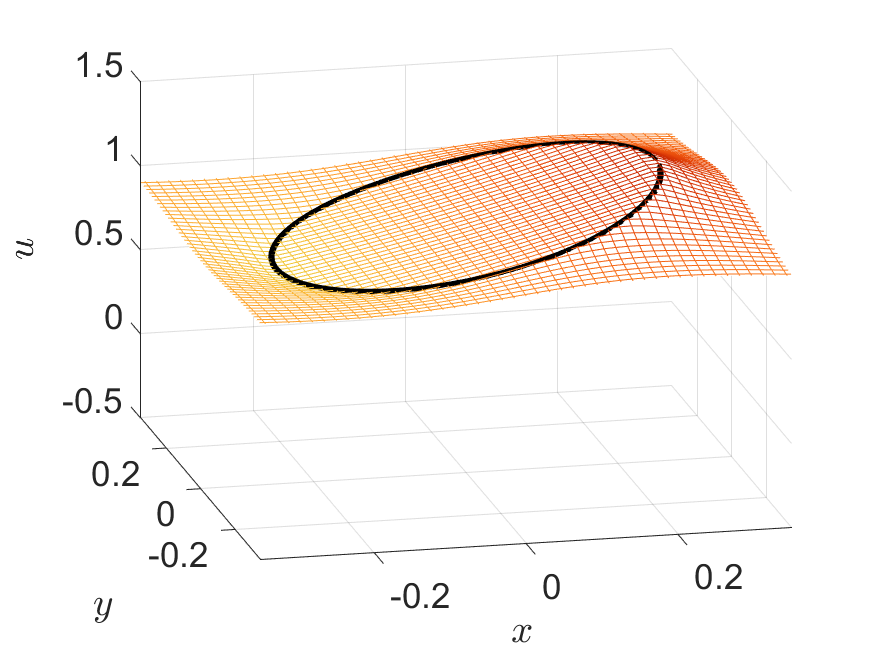}
\caption{\normalsize IBSL Solution}
\label{ibsl solution}
\end{subfigure}
\begin{subfigure}{0.495\textwidth}
\centering
\includegraphics[width=\textwidth]{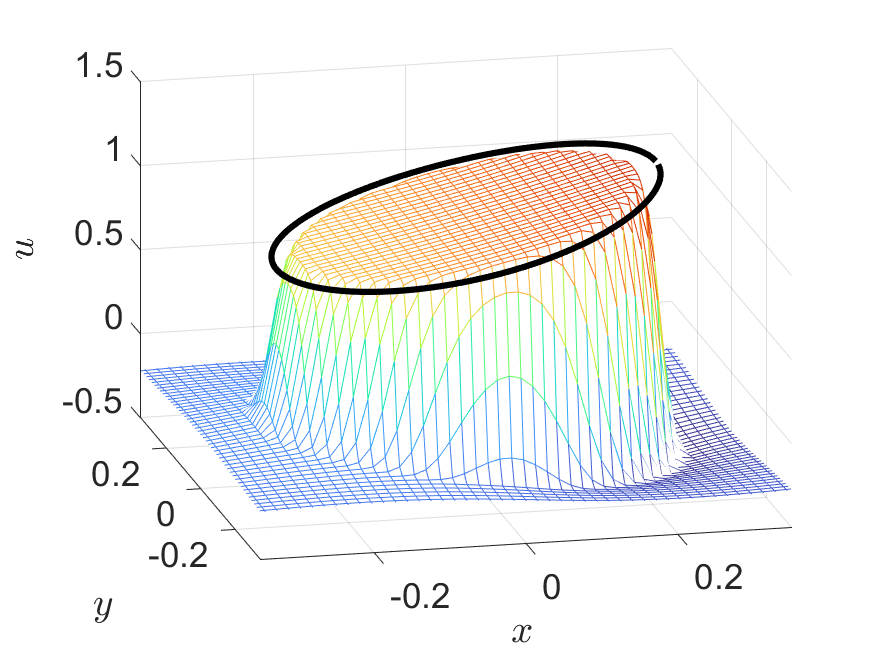}
\caption{\normalsize IBDL Solution}
\label{ibdl solution}
\end{subfigure}
\begin{subfigure}{0.65\textwidth}
\centering
\includegraphics[width=\textwidth]{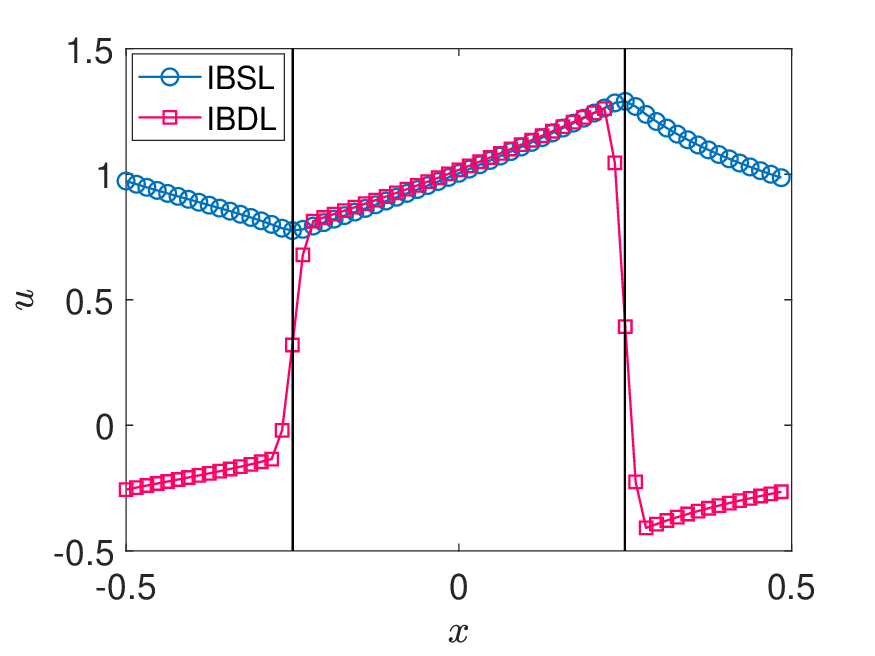}
\caption{\normalsize Solution Slices}
\label{slices}
\end{subfigure}
\caption{Solution plots found using the IBSL and IBDL methods to solve Equation \eqref{pde} with $k=1$ and $g=0$, and where the circular boundary of radius 0.25 has prescribed boundary values given by $U_b=e^x$. The periodic computational domain is $[-0.5, 0.5]^2$, and the grid and boundary point spacing is $\Delta s \approx \Delta x = 2^{-6}$. Figure \ref{ibsl solution} gives the solution from using the IBSL method, and Figure \ref{ibdl solution}, the IBDL method. $U_b$ is shown with a black curve. \chI{Figure \ref{slices} shows solution values for $y=-0.015625$ for both methods, and the black vertical lines show the location of the immersed boundary.}}\label{solution plots} 
\end{figure}

Since our boundary density $Q$ corresponds to the jump in the solution values across the boundary, clearly our solution $u$ will be \ch{discontinuous}. In the original IB constraint method, the solution was continuous, but the normal derivative was not. Figure \ref{solution plots} illustrates the solutions to Equation \eqref{pde} produced by the IBSL and IBDL methods on the entire computational domain, $\C=[-0.5,0.5]^2$. The \chIII{PDE} domain $\Omega$ is the interior of a circle, and the boundary \ch{condition} is $U_b=e^x$. We use the numerical implementation discussed in Section \ref{numerical implementation}.  We can see that the two solutions match on the portion of the problem domain, $\Omega$, that is about a couple meshwidths away from the boundary. However, the IBDL method gives a solution that is discontinuous across the boundary.

The lack of smoothness in the derivative across the boundary causes the IBSL method to achieve only first-order accuracy. Since the IBDL method will instead yield a discontinuous function, we will not see pointwise convergence near the boundary if we use the solution values given directly from the method in Equation \eqref{ibdl}. \chII{However, we do achieve first-order convergence in all norms in the region away from the boundary by a distance that is $\mathcal{O}(\Delta x)$. If we then wish to recover pointwise convergence up to the boundary, we can} replace the values of $u$ for grid points \textit{near} the boundary. 

Integral equation methods are generally able to achieve higher accuracy than the Immersed Boundary method. However, this lack of smoothness related to the singularity in the Green's function derivative results in the need to employ analytical techniques in order to achieve the same level of accuracy for points near $\Gamma$. For instance, Beale and Lai use a regularized Green's function and then analytically derive the correction terms \cite{BealeLai}. Kl\"{o}ckner et al. use analytic expansions centered at points several meshwidths from the boundary to \chIII{evaluate} the solution \chIII{at} points nearer to the boundary \cite{QBX}, and Carvalho et al. use asymptotic analysis to approximate the solution at near-boundary points with known boundary data and a nonlocal correction \cite{ShilpaBI}. 

\chI{For problems involving prescribed boundary values on sharp interfaces, the lack of smoothness of a solution from an IB method generically results in only first-order accuracy. While higher-order accuracy has been possible on problems with smooth solutions across the boundary \cite{Breugem, moripeskin, formally, Stein}, here we are restricting our attention to achieving the pointwise first-order convergence achieved by the IBSL method applied to a sharp interface.} \chII{Therefore, while there are many options for replacing the solution values near the boundary, we have chosen to use a simple linear interpolation because with this choice, we are able to obtain the first-order accuracy desired. This interpolation uses} known boundary values and approximate solution values several meshwidths into the PDE domain. \ref{interpolation appendix} provides details of the numerical implementation. 

\begin{figure}
\centering
\begin{subfigure}{0.495\textwidth}
\centering
\includegraphics[width=\textwidth]{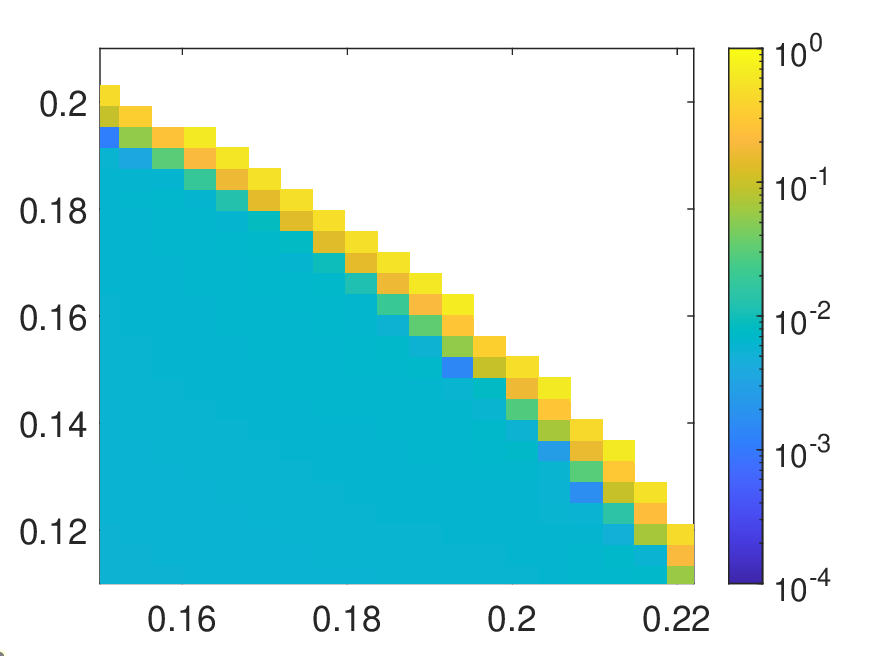}
\caption{\normalsize $\Delta x = 2^{-8}$}
\label{course}
\end{subfigure}
\begin{subfigure}{0.495\textwidth}
\centering
\includegraphics[width=\textwidth]{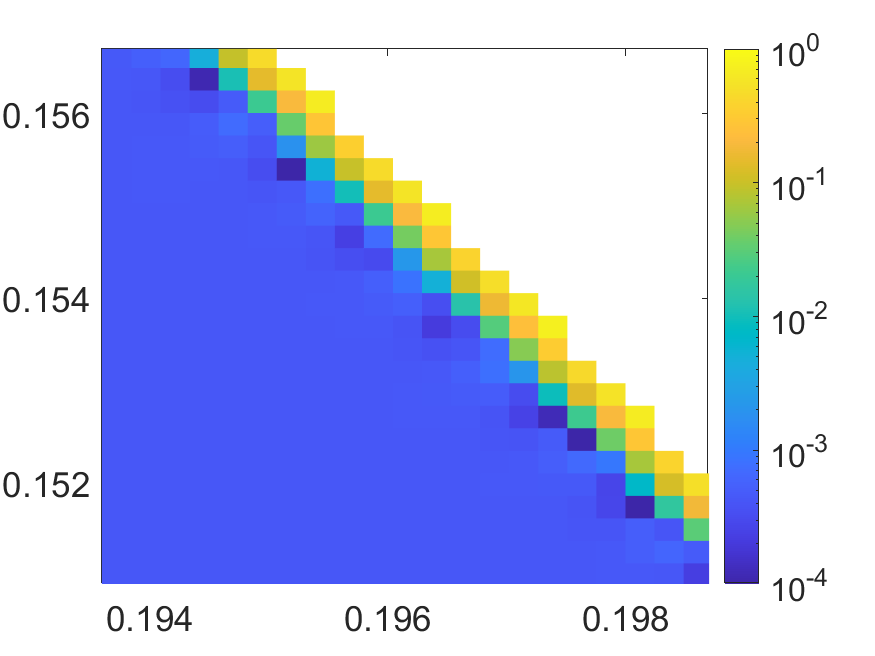}
\caption{\normalsize  $\Delta x = 2^{-12}$}
\label{fine}
\end{subfigure}
\caption{Plots of errors from using IBDL method to solve Equation \eqref{pde} with $k=1$ and $g=0$ on the periodic computational domain $[-0.5, 0.5]^2$, where the circular boundary of radius 0.25 has prescribed boundary values given by $U_b=e^x$. The errors are shown for two grid point spacings and with $\Delta s \approx \Delta x$. The figure is zoomed-in to view the large errors that lie in a small region near $\Gamma$.}\label{near boundary zoomed in}
\end{figure}

There are two regions of grid points for which we do not see pointwise convergence. The first region contains the largest errors, and it is illustrated in Figure \ref{near boundary zoomed in} for an interior circular domain with boundary data $U_b=e^x$. We can see that this region remains localized to only a small number of grid points, so that as we refine the grid, the width of this region quickly goes to 0. These large errors are the direct result of the smoothing of the dipole forces, and the number of meshwidths is therefore determined by the support of the regularized delta function used. \chI{As we discuss in Section \ref{numerical implementation}, we use the traditional Peskin four-point delta function in Equations \eqref{delta eqn 1}-\eqref{delta eqn 2} \cite{Peskin02}}, and this results in about $2-3$ meshdwidths of large errors on one side of the boundary, which can be seen in Figure \ref{near boundary zoomed in}. 

Secondly, the numerical method used to discretize the PDE spreads the error from the discontinuity into a larger region near the boundary. Therefore, if we use the interpolation to only correct for the large errors within the first few meshwidths, we will still see a lessening of first-order convergence in the max norm once the max norm is small enough. The width of the region on which the solution fails to converge pointwise still approaches 0, but the number of meshwidths affected can increase. \chIII{We} can recover pointwise convergence with further interpolation, and in practice, \chIII{we} can generally obtain a comparable maximum error to the IBSL method by using interpolation for only a relatively small number of meshwidths. In Section \ref{max norm stuff}, we further discuss factors that contribute to the number of meshwidths needed.

The IBDL method can now be described in the following steps:
\begin{enumerate} 
\item Use a Krylov method to solve the following system for $Q$. 
\begin{equation} 
-(S^* \L^{-1}\widetilde S)Q + \frac12 Q=U_b - S^*\L^{-1} \chII{\tilde g}\label{ibdl schur}
\end{equation}
\item Use $Q$ to solve the following equation for $u$ in $\C$. 
 \begin{equation} 
\L u + \widetilde S Q=\chII{\tilde g} \label{ibdl 1 again}
\end{equation}
\item Use an interpolation to replace the values for $u$ at $\x \in \Omega$ that are near $\Gamma$. Factors that determine the number of meshwidths to include in this step are discussed in Section \ref{max norm stuff}. \chII{Note: This step is \textit{optional}, and only required if \chIII{we wish} to recover pointwise convergence \textit{up to the boundary}.}
\end{enumerate}

%%%%%%%%%%%%%%%%%%%%% Arb param
\subsection{\chIII{Arbitrary boundary parametrization}}\label{arbparam}

\chIII{In this paper, we utilize an arclength parametrization of the immersed boundary $\Gamma$. Numerically, this means that in the discretization of the spread operator, $\Delta s_i$ is the size of the arclength connecting $\X(s_i)$ and $\X(s_{i+1})$. In practice, this choice is not a limiting factor because even if an analytical arclength parametrization of the immersed boundary curve is not available, \chIII{we can} still approximate this value with the distance $\Delta s_i\approx ||\X(s_{i+1})-\X(s_i)||_2$. However, to be thorough, in this section, we discuss the use of an arbitrary parametrization. }

\chIII{Let $\X(s)$ be an arbitrary smooth parametrization of $\Gamma$. Using this parametrization, the double layer integral equation for the Helmholtz equation would have the form 
\begin{equation}
U_b(\X(s'))= \int_{\Gamma}^{PV} \gamma(\X(s)) \grad G(\X(s),\X(s'))\dotp\n(\X(s)) \Big|\frac{\partial\X}{ds}\Big|d s-\frac12 \gamma(\X(s')).\label{param1}
		\end{equation}}

\chIII{In this case, we can form the IBDL method as
\begin{subequations}\label{ibdl param}
\begin{alignat}{2}
& \L u +\widetilde S Q = \tilde g \qquad && \text{in } \mathcal{C}  \label{ibdl 1 param}\\
&S^* u + \frac{1}{2|\partial\X/\partial s|} Q = U_b \qquad && \text{on } \Gamma.   \label{ibdl 2 param}
\end{alignat}
\end{subequations}
Then, using the same steps as in Section \ref{formulation}, we would obtain the equation 
\begin{multline}
U_b(\X(s'))= -\int_{\Gamma}^{PV} Q(\X(s)) \grad G_{hh}(\X(s),\X(s'))\dotp\n(\X(s)) ds\\+ \frac{1}{2|\partial\X/\partial s|}Q(\X(s')).\label{param2}
\end{multline}
By comparing Equations \eqref{param1} and \eqref{param2}, we see that in the limit $h\longrightarrow 0$, we get 
\begin{equation}
Q(\X)=-\gamma(\X) \Big|\frac{\partial\X}{\partial s}\Big|.
\end{equation}}

\chIII{As discussed in Section \ref{BI method}, the constant $1/2$ shifts all of the eigenvalues of the integral operator so that the only limit point is away from $0$. This ensures a small, constant condition number as the mesh is refined and therefore small iteration counts for the Krylov method. An arbitrary paramatrization may not shift the eigenvalues in this uniform manner. As such, the iteration counts may increase. In Section \ref{interior circle}, we revisit the Helmholtz equation from Section \ref{IB method} given by 
\begin{subequations} 
\begin{alignat}{2}
& \Delta u -  u = 0 \qquad && \text{in } \Omega  \\
&u=\sin{2\theta} \qquad && \text{on } \Gamma, 
\end{alignat}
\end{subequations}
where $\Omega$ is the interior of a circle of radius 0.25, centered at the origin. We will see that the iteration count for solving this equation with the IBDL method using arclength parametrization is about 4-5 iterations. If we instead use the parametrization given by 
\begin{equation}
\begin{pmatrix} x \\y\end{pmatrix} = \begin{pmatrix} 0.25 \cos{(s^2+s)} \\ 0.25\sin{(s^2+s)}\end{pmatrix}, 
\end{equation}
for $0\leq s\leq \frac12(-1+\sqrt{1+8\pi})$, we have 
\begin{equation}
\Big|\frac{\partial \X}{\partial s}\Big|=0.25(2s+1).
\end{equation}
Using immersed boundary points that are equally spaced in the parameter $s$, the iteration count is about 20. The iteration count is larger, but it does remain constant as the mesh is refined. This is likely to be the case with an arbitrary parametrization as the eigenvalues will still be shifted in some manner away from $0$. It is therefore acceptable to implement the IBDL method in this manner, but we then need the values or approximations of $|\partial \X/\partial s|$. To maintain consistency in the formulation of the IBDL method, we will use an arclength parametrization in all other sections of this paper. }

%%%%%%%%%%%%%%%%%% NEUMANN SET UP %%%%%%%%%%%%%%%%%%%%%%
\subsection{Neumann boundary conditions}\label{Neumann formulation}

We now use the IB framework and our connection to boundary integral equations to solve a PDE with Neumann boundary conditions. The original IBSL method was unable to handle such a problem since the derivatives are not convergent at the boundary. 

Let us look at the PDE given by
\begin{subequations} \label{pde neumann}
\begin{alignat}{2}
& \Delta u - k^2 u = g \qquad && \text{in } \Omega  \label{pde1 neumann}\\
&\frac{\partial u}{\partial n} =V_b \qquad && \text{on } \Gamma,  \label{pde2 neumann}
\end{alignat}
\end{subequations}
where we will first take $g=0$ for simplicity. In Section \ref{BI method}, we introduced the full integral representation of $u$. For $\x_0 \in \Omega\setminus \Gamma $, we have
\begin{equation}
u(\x_0) = \int_{\Gamma}  \grad u(\x) \dotp \n (\x) G(\x, \x_0)dl(\x) \\- \int_{\Gamma} u(\x) \grad G(\x,\x_0) \dotp \n(\x) dl(\x). \label{integral rep again}
\end{equation}  
The first term in Equation \eqref{integral rep again} is a single layer potential with the strength given by the known boundary derivatives, $V_b$. The second term is a double layer potential whose strength is given by the \ch{\textit{unknown}} boundary values. If we let $U_b\equiv u\big|_{\Gamma}$, we can then write the PDE in the IB framework as
\begin{equation}
 \L u +\widetilde S U_b+SV_b = 0 \label{homog neumm 1}
\end{equation}

To determine the interpolation equation, we can take $\x_0\longrightarrow \Gamma$ in Equation \eqref{integral rep again} and denote it as $\X_0$. Using the same property of the Green's function derivative that was used in Section \ref{BI method}, we get
\begin{multline}
u(\X_0) = \int_{\Gamma}  \grad u(\x) \dotp \n (\x)G(\x, \X_0) dl(\x) - \int_{\Gamma}^{PV} u(\x) \grad G(\x,\X_0) \dotp \n(\x) dl(\x) \\+\frac{1}{2} u(\X_0). \label{integral eqn again}
\end{multline} 
By combining the $u(\X_0)$ terms, we can identify that the interpolated solution values will correspond to $1/2u(\X_0)$. Allowing for a non-zero $g$, we can therefore write the PDE in the Immersed Boundary framework as 
\begin{subequations} \label{ib neumann}
\begin{alignat}{2}
& \L u +\widetilde S U_b+SV_b = \chII{\tilde g} \qquad && \text{in } \mathcal{C}  \label{ib neumann 1}\\
&S^* u = \frac12 U_b \qquad && \text{on } \Gamma.   \label{ib neumann 2}
\end{alignat}
\end{subequations}
For Neumann boundary conditions, $U_b$ is the \textit{unknown} potential strength on the boundary, corresponding to the \ch{unknown} boundary values. A detailed demonstration of the connection between Equation \eqref{integral eqn again} and Equation \eqref{ib neumann} is presented in \ref{neumann appendix}, and it uses steps similar to those used in Sections \ref{connection} and \ref{formulation}.  

We can solve this system by first solving the equation 
\begin{equation}
(-S^* \L^{-1}\widetilde S)U_b-\frac12 U_b=S^*\L^{-1}SV_b-S^*\L^{-1}\chII{\tilde g} \label{Neumann saddle point}
\end{equation} 
for the boundary values, $U_b$, and then obtaining $u$ from Equation \eqref{ib neumann 1}. We can see by comparing to Equation \eqref{ibdl schur} that the operator is similar to that of the IBDL method for Dirichlet boundary conditions, where the only difference is the sign on the $1/2$. Therefore, we again get a well-conditioned problem that can be solved with a small number of iterations of a Krylov method. An example of a Neumann problem will be presented in Section \ref{neumann}.

%%%%%%%%%%%%%%%%%%%%%%%%%%%%%%%%%%%%%%%%%%%%%%%%%%%%%%%%%%%%%%%%    NUMERICAL IMPLEMENTATION    %%%%%%%%%%%%%%%%%%%
%%%%%%%%%%%%%%%%%%%%%%%%%%%%%%%%%%%%%%%%%%%%%%
\section{Numerical implementation}\label{numerical implementation}

\cha{This section provides a detailed description of numerical implementation for the method in 2-D. The extension to 3-D is straightforward, and we give an example in Section 6.7.} 

\subsection{\chIII{Discretization of space and differential operators}} \label{computational domain}

In this paper, we take the computational domain to be a periodic box of length $L$ in order to make use of efficient finite difference and Fourier spectral methods to solve Equation \eqref{ibdl 1 again}.  As is typical for IB methods, there is flexibility in the choice of PDE solver. However, as is discussed in Section \ref{max norm stuff}, the discontinuity of the solution affects the accuracy near the boundary differently for different discretizations of the PDE. Therefore, we use a finite difference method for all applications other than those presented to compare methods in Section \ref{max norm stuff}. We use the standard five-point, second-order accurate approximation for the Laplacian. 

The computational domain is discretized with a regular Cartesian mesh with $N$ points in each direction, giving us $\Delta x = \Delta y = L/N$. To discretize the boundary, we use a set of $N_{IB}$ boundary points, given by $\{X(s_i)\}_{i=1}^{N_{IB}}$. 

\chb{For our refinement studies, we rescale the discrete $L^1(\Omega)$ and $L^2(\Omega)$ norms by the area or volume of $\Omega$, denoted $|\Omega|$. Our refinement studies therefore use the discrete function norms given by 
\begin{subequations}
\begin{alignat}{2}
&||u||_1 &&= \frac{\Delta x \Delta y}{|\Omega|} \sum_{i,j} u(\x_{i,j})\chi_{\scaleto{\Omega}{4.5pt}}\\
&||u||_2 &&= \Bigg(\frac{\Delta x \Delta y}{|\Omega|} \sum_{i,j} (u(\x_{i,j}))^2\chi_{\scaleto{\Omega}{4.5pt}}\Bigg)^{1/2}.
\end{alignat}
\end{subequations}}

\subsection{Discretization of spread and interpolation operators}\label{discret}%%%%%%%%%%%%%%%%%%%%%%%%%%%%%%%%%%%%%%%%%%%%%%%%%%%%%%%%%%%%%

We will first discuss the implementation of our operator defined by 
\begin{equation}
\widetilde S Q = \grad \dotp (S Q\n). \label{Stilde}
\end{equation}
This operation consists of first multiplying our discretized boundary density $Q$ by the outward unit normal vectors  at the immersed boundary points, whose locations will be discussed more in depth shortly. We then spread these vectors to the grid and take the divergence. We discuss the implementation of each of these steps below. 

The first component of the operator $\widetilde S$ is the divergence. To compute this, we use a finite difference  method to match the method used to discretize the PDE. We use the standard second-order accurate centered difference for the derivatives.

Next, we need the unit normal vectors to the immersed boundary. These can be treated as an input to the method, or we can approximate them if they are unavailable. Approximating the curve with linear elements and using the resulting normal vectors is sufficient for maintaining first-order accuracy. For instance, for an interior domain and a counterclockwise curve parametrization, we can use 
\begin{equation}
\n(s_i)= \begin{bmatrix}Y(s_{i+1})-Y(s_{i-1})\\   -\big(X(s_{i+1})-X(s_{i-1})\big)\end{bmatrix} \Bigg/ \big|\big|\X(s_{i+1})-\X(s_{i-1})\big|\big|_2.\label{normal vector approximations}
\end{equation}
For an exterior domain, we would negate this.

Next, we look at the spread operator, $S: \Gamma \rightarrow \C$, which is defined as
\begin{equation}
(SF)(\x) =\int_{\Gamma} F(s)\delta_h(\x-\X(s))ds .    \label{spread operator again}
\end{equation}
Here, we use the traditional Peskin four-point delta function for $\delta_h$ \cite{Peskin02}, given by 
\begin{equation}
\delta_h=\frac{1}{h^2} \phi\bigg(\frac{x}{h}\bigg)\phi\bigg(\frac{y}{h}\bigg), \label{delta eqn 1}
\end{equation}
for $h=\Delta x = \Delta y$ and 
\begin{equation}
\phi (r)=\begin{cases}
 \frac{1}{8} \big( 3-2|r| + \sqrt{1+4|r|-4|r|^2}\big) & 0\leq |r|\leq 1 \\ 
      \frac{1}{8} \big( 5-2|r| - \sqrt{-7+12|r|-4|r|^2}\big) & 1\leq |r|\leq 2 \\
	0 & |r|\geq 2   \end{cases}. \label{delta eqn 2}
\end{equation}

Then, as mentioned in Section \ref{formulation}, we discretize the integral in Equation \eqref{spread operator again} with respect to arclength, so henceforth, \chI{we} will treat $s$ as an arclength parameter. With this new assumption on the coordinate $s$, $\Delta s_i$ gives the length of the curve connecting the immersed boundary points, $\X(s_i)$ and $\X(s_{i+1})$. Then, the discrete spread operator is given by
\begin{equation}
(SF)(\x) =\sum_{i=1}^{N_{IB}} F(s_i)\delta_h(\x-\X(s_i))\Delta s_i .    \label{spread operator again2}
\end{equation}

Depending on the parametrization, the boundary points may not be equally spaced, but since we will be looking at a periodic distribution on $\Gamma$, choosing equally spaced points will give higher accuracy in our integral approximation. In this case, we will space the boundary points so that $\Delta s \approx c \Delta x$ for various $c$ values, and then $\Delta s = L_{IB}/N_{IB}$, where $L_{IB}$ gives the length of the immersed boundary. If the exact values of $\Delta s_i$ are unavailable, approximations can be made. Again, using $\Delta s_i \approx ||\X(s_{i+1})-\X(s_i)||_2$ is sufficient to maintain first-order accuracy.

Lastly, the interpolation operator, $S^*: \C \rightarrow \Gamma$, is similarly discretized as follows.
\begin{equation}
(S^*u)(s)=\int_{\C} u(\x) \delta_h(\x-\X(s)) d\x \approx \sum_{i,j} u(\x_{i,j})\delta_h(\x_{i,j}-\X(s)) (\Delta x)(\Delta y),  \label{interpolation again}
\end{equation}
where, for a fixed $s$ and a four-point delta function, the sum over the grid points has at most 16 non-zero terms.

%%%%%%%%%%%%%%%%%%%%%%%%%%%%%%%%%%%%%%%%%%%%%%%%%%%%%%%%%%%%%%%%%%%%%%%%%%%%%%
\subsection{\chIII{Solution to discrete system for invertible $\L$ } }\label{5.3 invertible L}

\chIII{We now describe the method for solving the discretized system 
\begin{subequations} \label{ibdl again}
\begin{alignat}{2}
& \L u +\widetilde S Q = \tilde g \qquad && \text{in } \mathcal{C}  \label{ibdl again 1}\\
&S^* u + \frac12 Q = U_b \qquad && \text{on } \Gamma.   \label{ibdl again 2}
\end{alignat}
\end{subequations}}
\chIII{In the case that the differential operator is invertible, such as for $\L=\Delta -k^2$, for $k\neq0$, we can invert the operator to obtain
\begin{equation}
u=-\L^{-1}\widetilde S Q+\L^{-1}\tilde g  .
\end{equation}
Then, by applying the interpolation operator and using Equation \eqref{ibdl again 2}, we obtain
\begin{equation}
-(S^* \L^{-1}\widetilde S)Q + \frac12 Q=U_b - S^*\L^{-1} \tilde g . \label{ibdl schur again}
\end{equation}
We can therefore solve Equation \eqref{ibdl again} by first solving Equation \eqref{ibdl schur again} for $Q$ and then obtaining $u$ from Equation \eqref{ibdl again 1}. Here, the operator that must be inverted, $-(S^* \L^{-1}\widetilde S) + \frac12 \mathds{I}$, is not symmetric, and we therefore solve for $Q$ with \texttt{gmres}\footnote{\chI{We utilize \texttt{gmres} from \MATLAB R2022a \cite{Matlab22a}.}}, to a tolerance of $10^{-8}$, unless otherwise specified. }

%%%%%%%%%%%%%%%%%%%%%%%%%%%%%%%%%%%%%%%%%%%%%%%%%%%%%%%%%%%%%%%%%%%%%%%%%%%%%%%%%%%%
\subsection{\chIII{Solution to discrete system for Poisson equation}}\label{poisson explanation}

\chIII{In the case that the differential operator is the periodic Laplacian, $\L=\Delta$, the solution method outlined in the previous section must be adjusted, both for the nullspace of the discrete constraint system and the nullspace of $\L$ . }

\chIII{The boundary value problem given by 
\begin{subequations} \label{laplaces pde}
\begin{alignat}{2}
& \Delta u = g \qquad && \text{in } \Omega \label{l pde1}\\
&u=U_b \qquad && \text{on } \Gamma \label{l pde2}
\end{alignat}
\end{subequations}
has a unique solution in $\Omega$, but, by using the IB framework, we embed the PDE into a periodic computational domain on which $\Delta$ is not invertible, resulting in a nullspace for the differential operator $\L$. Furthermore, the discretized constraint matrix for the IBDL method also has a nullspace that becomes relevant for large exterior domains. Both issues are discussed in this section. }

\subsubsection{\chIII{Nullspace of the constraint system }}\label{nullspace1}

 \chIII{In the case of the IBSL method applied to the Poisson equation, the matrix
\begin{equation}
\begin{pmatrix} \Delta &S\\ S^* &0\end{pmatrix},
\end{equation}
is invertible. Therefore, we can use the method presented by Stein et al. \cite{Stein} to handle the nullspace of the periodic Laplacian and recover the unique solution to the constraint system. }

\chIII{In the case of the IBDL method, however, the matrix
\begin{equation}
\begin{pmatrix} \Delta &\widetilde S\\ S^* &\frac12\mathbb{I}\end{pmatrix},
\end{equation}
has a one-dimensional nullspace. The nullspace of the spatially continuous problem is spanned by 
\begin{subequations}\label{nullspace shit}
\begin{alignat}{2}
&u=\chi_{\scaleto{\C \setminus \Omega}{6pt}}\\
&Q=-1,
\end{alignat}
\end{subequations}
where $\chi_{\scaleto{\C \setminus \Omega}{6pt}}$ gives the indicator function that is $0$ on $\Omega$ and $1$ otherwise. According to Equation \eqref{nullspace shit}, since an element of the nullspace would be $0$ on the physical domain, any solution we get would be the same on $\Omega$, and this nullspace would therefore not pose a problem for the IBDL method. However, the discretized problem results in a different \textit{numerical} nullspace. In Section \ref{poisson1}, we show that for large exterior domains, the IBDL method may not converge to the correct solution. This may be due to the fact that as we increase the domain size, we approach the infinite exterior domain case, for which a double layer potential is unable to represent solutions to the Poisson equation with arbitrary boundary values \cite{Periodic}. As we approach this case, the IBDL method becomes more sensitive to the numerical nullspace. }

\chIII{To solve this issue for large exterior domains, we may use a \textit{completed} double layer representation. There are multiple ways to formulate such a representation for integral equations \cite{powermiranda, power}, but we have chosen the method presented by Hsiao and Kress  \cite{hsiaokress} and Hebeker \cite{hebeker} since it easily translates to an IB framework. We therefore formulate the completed representation by the addition of a single layer potential whose strength is a constant multiple of the strength of the double layer potential. The completed system is then given by 
\begin{subequations} \label{ibdl completed}
\begin{alignat}{2}
& \Delta u +\eta S Q+ \widetilde S Q = \tilde g \qquad && \text{in } \mathcal{C}  \label{ibdl 1 completed}\\
&S^* u + \frac12 Q = U_b \qquad && \text{on } \Gamma,  \label{ibdl 2 completed}
\end{alignat}
\end{subequations}
where $\eta$ is a positive constant. We will refer to this form of the IBDL method as the \textit{completed IBDL method}. In this case, the matrix 
\begin{equation}
\begin{pmatrix} \Delta &\widetilde S+\eta S\\ S^* &\frac12\mathbb{I}\end{pmatrix},
\end{equation}
is invertible, and then we need only deal with the nullspace of the periodic Laplacian.}

\subsubsection{\chIII{Nullspace of $\L$}} \label{nullspace2}

\chIII{As previously mentioned, the original boundary value problem on $\Omega$ has a unique solution, but, by using the IB framework, we embed the PDE into a periodic computational domain on which $\Delta$ is not invertible. }

\chIII{We can recover the unique solution to the PDE by a method similar to that presented by Stein et al.\ \cite{Stein}. In this process, we first decompose the solution $u$ into 
\begin{equation}
u=u_0+\bar u,
\end{equation}
where $u_0$ has mean $0$ on $\C$ and $\bar u$ gives the mean value of $u$ on $\C$. Then, we use the PDE resulting from either the IBDL or the completed IBDL formulation to derive a solvability condition, resulting in either a condition on $g_e$ or in the addition of another equation to the discrete IBDL system. Details for this process are given in \ref{nullspace appendix}.}

%%%%%%%%%%%%%%%%%% NEAR BOUNDARY %%%%%%%%%%%%%%%%%%%%%%
\subsection{Near-boundary points}\label{near boundary}

The replacement of solution values for grid points within $m_1$ meshwidths of  $\Gamma$ is done with a linear interpolation. For the Dirichlet problem, the interpolation for an individual grid point uses one approximate solution value located about $m_2$ meshwidths away from $\Gamma$ and one \chIII{value approximated from the known solution values} located on $\Gamma$. For the Neumann problem, the solution values located on $\Gamma$ are not known, and therefore they are also approximate. Details for this interpolation are given in \ref{interpolation appendix}.

In order to complete this step and to isolate the solution on the problem domain, it is necessary to identify grid points as interior or exterior to $\Gamma$. Depending on the boundary shape, this can be done with the exact parametrization of the curve, but \chIII{we can} also use a portion of the IBDL method itself to identify grid points interior to the \textit{discretized} curve. Details for this are given in \ref{in out appendix}.

%%%%%%%%%%%%%%%%%%%%%%%%%%%%%%%%%%%%%%%%%%%%%%%%%%%%%%%%%%%%%%%%    RESULTS  1  %%%%%%%%%%%%%%%%%%%
%%%%%%%%%%%%%%%%%%%%%%%%%%%%%%%%%%%%%%%%%%%%%%
\section{Numerical results}\label{Results}

In this section, we apply the Immersed Boundary Double Layer method and discuss numerical results and observations. In Section \ref{interior circle}, we apply the IBDL method to a \cha{2-D} \ch{Dirichlet} Helmholtz problem and demonstrate the efficiency of the method and the convergence of the solution. In Section \ref{max norm stuff}, we discuss interpolation widths and other factors affecting the max norm convergence. \chIII{In Section \ref{poisson1}, we explore the completed IBDL method for the Dirichlet Poisson problem on a periodic domain. In Section \ref{poisson2}, we demonstrate the use of the IBDL method on a non-convex exterior domain.} In Section \ref{Q convergence}, we explore the convergence of the potential strength $Q$, and in Section \ref{neumann}, we solve a Neumann Helmholtz problem. \cha{Finally, in section \ref{3D}, we solve a 3-D Helmholtz problem with Dirichlet boundary conditions.}

%%%%%%%%%%%%%%%%%%%%%%%%%%%%%%%%%%%%%%%%%%%%%%%%%%
%%%%%%%%%%%%%%%%%% INTERIOR   %%%%%%%%%%%%%%%%%%%%%%
%%%%%%%%%%%%%%%%%%%%%%%%%%%%%%%%%%%%%%%%%%%%%%%%%%

\subsection{Dirichlet Helmholtz on interior circular domain}\label{interior circle}
First, we look at the PDE
\begin{subequations} \label{pde again}
\begin{alignat}{2}
& \Delta u -  u = 0 \qquad && \text{in } \Omega  \label{pde1 again}\\
&u=\sin{2\theta} \qquad && \text{on } \Gamma,  \label{pde2 again}
\end{alignat}
\end{subequations}
where $\Omega$ is the interior of a circle of radius 0.25, centered at the origin. The analytic solution is given by 
\begin{equation}
u=\frac{I_2(r)\sin{2\theta}}{I_2(0.25)}, 
\end{equation}
where $I_2$ is the first-kind modified Bessel function of order 2. Our computational domain here is the periodic box $[-0.5, 0.5]^2$. We use equally spaced boundary points with $\Delta s \approx c \Delta x$ for various values of $c$, and the solutions are computed for grid sizes ranging from $N=2^6$ to $2^{12}$. 

\begin{figure}
\begin{center}
 \begin{tabular}{||c | c | c | c | c | c | c | c | c||} 
 \hline
 \multicolumn{9}{|c|}{Iteration Counts - Circular Boundary} \\
 \hline
  $\Delta x$ &\multicolumn{2}{|c|}{$\Delta s \approx  2\Delta x $}& \multicolumn{2}{|c|}{$\Delta s \approx  1.5 \Delta x$} &\multicolumn{2}{|c|}{$ \Delta s \approx  1 \Delta x$}& \multicolumn{2}{|c|}{$\Delta s \approx 0.75 \Delta x$}  \\ [0.5ex] 
 \hline
 &\textcolor{blue}{ IBSL} & \textcolor{cyan}{IBDL} &\textcolor{blue}{  IBSL}& \textcolor{cyan}{IBDL}&\textcolor{blue}{ IBSL} &  \textcolor{cyan}{IBDL}&\textcolor{blue}{  IBSL} &  \textcolor{cyan}{IBDL} \\
 \hline
$2^{-6}$ &  \textcolor{blue}{  17} & \textcolor{cyan}{4} &    \textcolor{blue}{45}&\textcolor{cyan}{5}   &  \textcolor{blue}{249}&\textcolor{cyan}{5}  &   \textcolor{blue}{491}&\textcolor{cyan}{4}   \\
$2^{-7} $&   \textcolor{blue}{36}& \textcolor{cyan}{5}   &   \textcolor{blue}{57}&\textcolor{cyan}{4}    &  \textcolor{blue}{691}& \textcolor{cyan}{5} &  \textcolor{blue}{752}&\textcolor{cyan}{4}   \\
$2^{-8 }$&   \textcolor{blue}{49}& \textcolor{cyan}{5}  &    \textcolor{blue}{71}&\textcolor{cyan}{4}   & \textcolor{blue}{1233}& \textcolor{cyan}{4}   &  \textcolor{blue}{3057 }&\textcolor{cyan}{4}      \\
$2^{-9}$& \textcolor{blue}{ 61}&\textcolor{cyan}{4}    &     \textcolor{blue}{110}&\textcolor{cyan}{4}     &  \textcolor{blue}{1922}& \textcolor{cyan}{4}   &   \textcolor{blue}{5829}&\textcolor{cyan}{4}        \\
$ 2^{-10}$& \textcolor{blue}{68}&\textcolor{cyan}{4}   &     \textcolor{blue}{144}&\textcolor{cyan}{4}     &  \textcolor{blue}{1936}& \textcolor{cyan}{4}   &   \textcolor{blue}{ 8084 }&\textcolor{cyan}{4}    \\
$ 2^{-11}$ & \textcolor{blue}{95}&\textcolor{cyan}{4}  &   \textcolor{blue}{234}&\textcolor{cyan}{4}     &  \textcolor{blue}{ 4364}& \textcolor{cyan}{4}  &   \textcolor{blue}{9335 }&\textcolor{cyan}{4}    \\
$ 2^{-12}$  & \textcolor{blue}{142}&\textcolor{cyan}{4} & \textcolor{blue}{303}&\textcolor{cyan}{4}       &  \textcolor{blue}{ 4535}& \textcolor{cyan}{4}  &   \textcolor{blue}{10589  }&\textcolor{cyan}{4}   \\
  \hline
  \end{tabular}
\captionof{table}{Number of iterations of \texttt{minres} and \texttt{gmres}, with tolerance $10^{-8}$, for the IBSL and IBDL methods, respectively, to solve Equation \eqref{pde again} on the periodic computational domain $[-0.5, 0.5]^2$, where the circular boundary of radius 0.25 has prescribed boundary values given by $U_b=\sin{2\theta}$.}
\label{iteration table 2}
\end{center}
\end{figure}

\begin{figure}
\centering
\begin{subfigure}{0.495\textwidth}
\centering
\includegraphics[width=\textwidth]{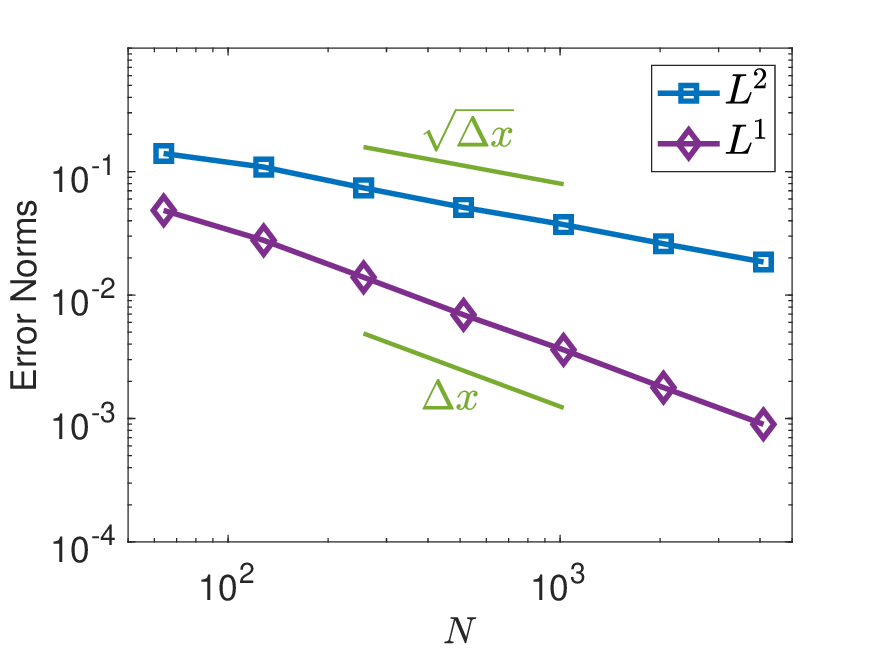}
\caption{\normalsize IBDL Refinement, no interpolation}
\label{nointerpfigure}
\end{subfigure}
\begin{subfigure}{0.495\textwidth}
\centering
\includegraphics[width=\textwidth]{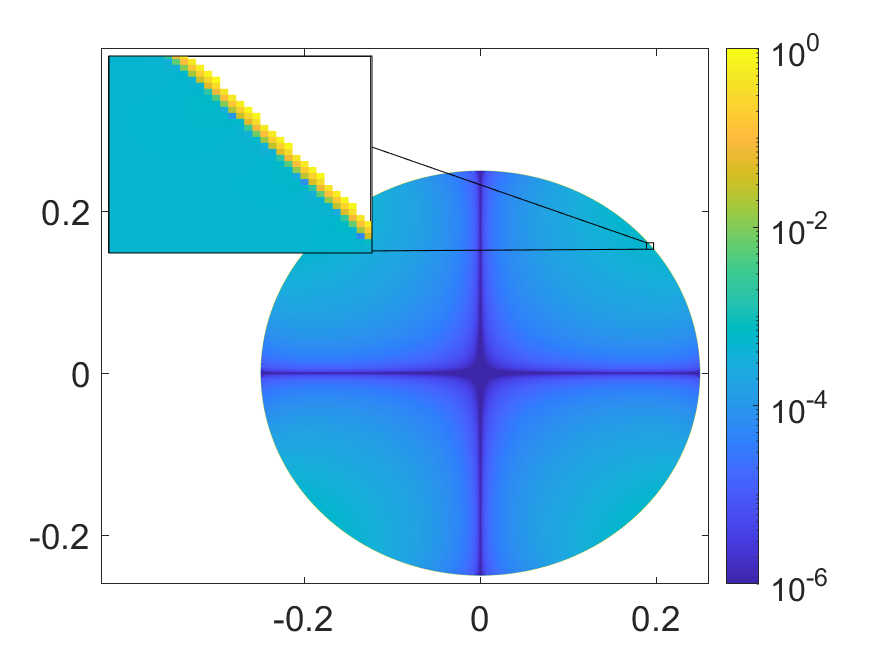}
\caption{\normalsize IBDL Refinement, no interpolation}
\label{nointerpfigure2}
\end{subfigure}
\begin{subfigure}{0.495\textwidth}
\centering
\includegraphics[width=\textwidth]{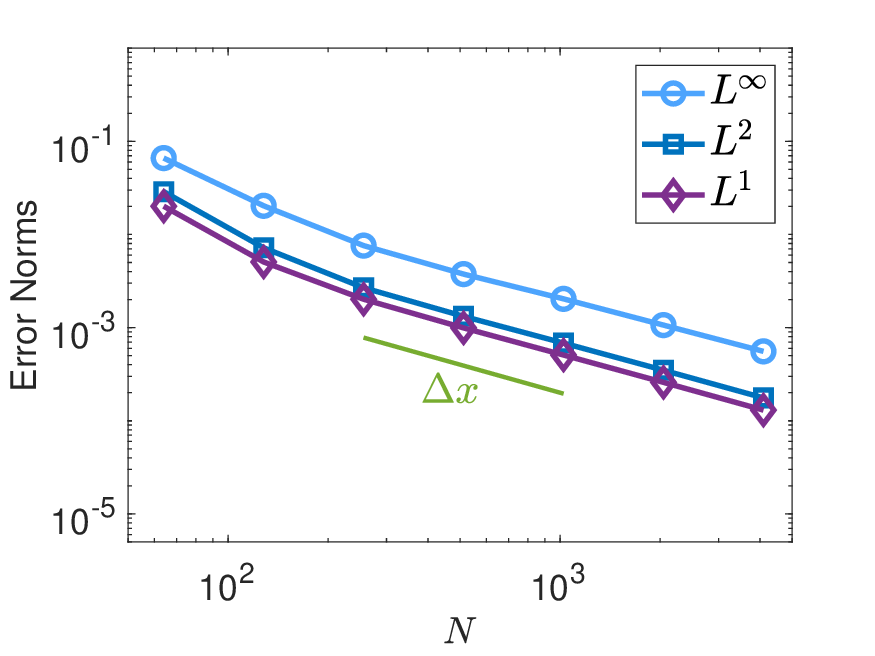}
\caption{\normalsize  IBDL Refinement}
\label{ibdl interior refinement}
\end{subfigure}
\begin{subfigure}{0.495\textwidth}
\centering
\includegraphics[width=\textwidth]{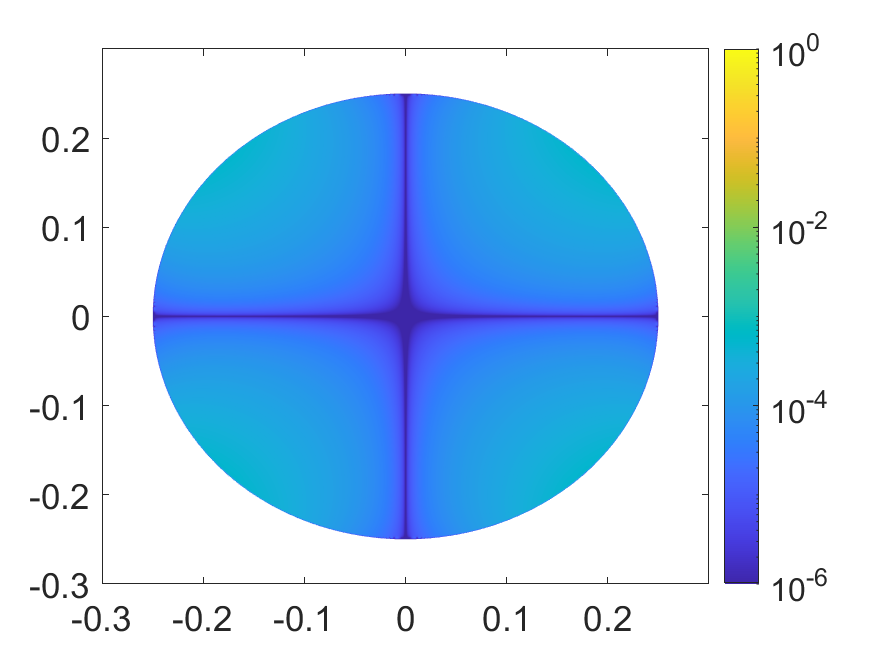}
\caption{\normalsize  IBDL Error}
\label{ibdl interior error}
\end{subfigure}
\begin{subfigure}{0.495\textwidth}
\centering
\includegraphics[width=\textwidth]{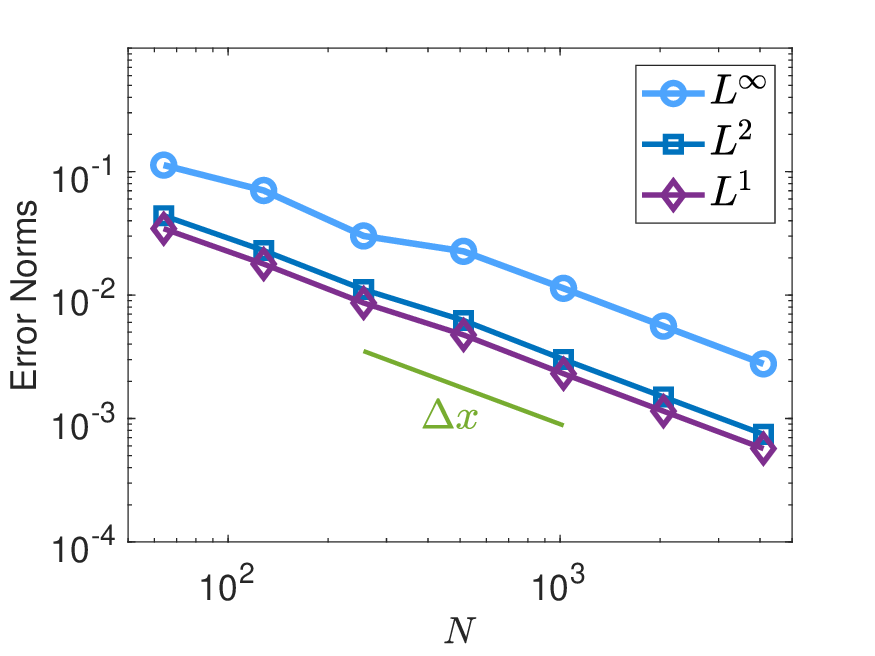}
\caption{\normalsize IBSL Refinement}
\label{ibsl interior refinement}
\end{subfigure}
\begin{subfigure}{0.495\textwidth}
\centering
\includegraphics[width=\textwidth]{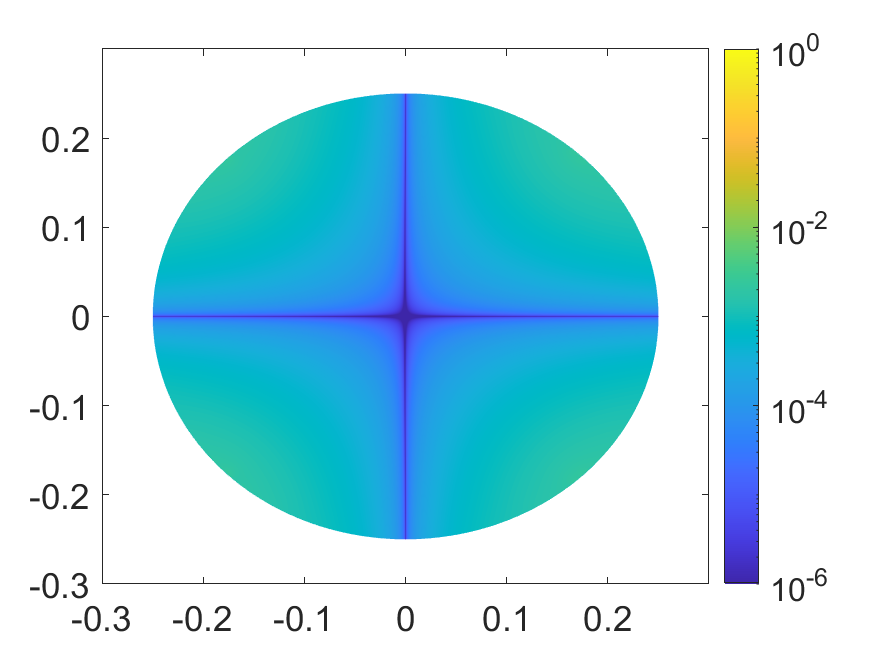}
\caption{\normalsize IBSL Error}
\label{ibsl interior error}
\end{subfigure}
\caption{\chIII{Refinement studies and error plot for solutions to Equation \eqref{pde again} found using the IBSL and IBDL methods. The computational domain is the periodic box $[-0.5, 0.5]^2$, $\Omega$ is the interior of a circle of radius 0.25, and the prescribed boundary values are given by $U_b=\sin{2\theta}$. Figures \ref{nointerpfigure}-\ref{nointerpfigure2} give the plots for the IBDL method, with no interpolation step. Figures \ref{ibdl interior refinement}-\ref{ibdl interior error} are for the IBDL method with an interpolation step, replacing solution values within $m_1= 6$ meshwidths from the boundary using an interior interpolation point $m_2=8$ meshwidths away from the boundary. Figures \ref{ibsl interior refinement}-\ref{ibsl interior error} are for the IBSL method. The error plots are shown for the grid spacing $\Delta x = 2^{-12}$. All methods use boundary point spacing of $ \Delta s  \approx 0.75 \Delta x$. }}
\end{figure}

\subsubsection{Iteration counts and point spacing}

Table \ref{iteration table 2} gives the iteration counts for the Krylov methods used to solve Equation \eqref{pde again} for the IBSL and IBDL methods. This table illustrates the drastic improvement in efficiency that is a key advantage of the IBDL method. We see that only a small number of iterations are needed to solve for the strength of the IBDL  potential. Moreover, the number of iterations remains essentially constant as we refine the mesh or tighten the boundary points relative to the grid. The IBSL method, on the other hand, requires many more iterations, and this iteration count is greatly affected by the mesh size and the boundary point spacing\chI{, as previously illustrated in Figure \ref{iteration plot}}. Tighter boundary points can result in lower quadrature error, while more widely spaced points allow for smoother force distributions and better conditioning in the case of the IBSL method \cite{GriffithDonev}. Making the IBSL method practical therefore often entails special handling, such as using more widely spaced points for discretization of the boundary, while using a denser set of points for the direct Eulerian-Lagrangian interaction \cite{Griffithpointspacing}. Moreover, to make the IBSL method efficient, one must use a preconditioner  \cite{Ceniceros, GriffithDonev, guyphilipgriffith, Stein}. On the other hand, the fact that the IBDL method presented here can be efficiently utilized for both coarsely \chIII{and} finely spaced boundary points means no such step is required for this method. 

In addition, the iteration counts for the IBSL method are affected by the tolerance for the Krylov method. For example, for $\Delta s \approx 2\Delta x$ and $\Delta x = 2^{-12}$, decreasing the tolerance from $10^{-8}$ to $10^{-10}$ increases the iteration count by more than fifty percent, from 142 to 233. For $\Delta s \approx 0.75 \Delta x$, the Krylov method cannot reach a tolerance of $10^{-10}$ for the IBSL method. For the IBDL method, on the other hand, the Krylov method reaches this smaller tolerance in only 5 iterations for either of these cases. 

\subsubsection{\chII{Accuracy}}

\chII{To illustrate that the IBDL gives first-order convergence away from the boundary by a distance that is $\mathcal{O}(\Delta x)$, Figure \ref{nointerpfigure} shows the $L^1$ and $L^2$ refinement studies for the IBDL solution to Equation \eqref{pde again}, obtained without the interpolation step. We observe $\mathcal{O}(\Delta x)$ convergence for $L^1$ and $\mathcal{O}(\sqrt{\Delta x})$ convergence for $L^2$. This indicates that we get first-order pointwise convergence in the solution away from the boundary. This is also illustrated in Figure \ref{nointerpfigure2}, where the large errors are seen only within a few meshwidths of the boundary. Therefore, as previously discussed, the interpolation step is not necessary unless we wish to recover pointwise convergence in this small region near the boundary. which is $\mathcal{O}(\Delta x)$ in  width. }

\chIII{We then implement the IBDL method, including the interpolation step, in which we replace solution values within $m_1=6$ meshwidths of the boundary and use interpolation points located $m_2=8$ meshwidths from the boundary.}  Figures \ref{ibsl interior refinement}-\ref{ibdl interior refinement} compare the IBSL and IBDL refinement studies, and \chIII{we see that with the interpolation step, we maintain first-order convergence of the solution on the entire PDE domain.} Figures \ref{ibsl interior error}-\ref{ibdl interior error} illustrate the solution errors for the two methods, and the errors from the IBDL method are actually lower for this problem. In general, the two methods typically give errors of comparable size.

%%%%%%%%%%%%%%%%%%%%%%%%%%%%%%%%%%%%%%%%%%%%%%%%%%
%%%%%%%%%%%%%%%%%% 6.2 Max norm stuff   %%%%%%%%%%%%%%%%%%%%%%
%%%%%%%%%%%%%%%%%%%%%%%%%%%%%%%%%%%%%%%%%%%%%%%%%%

\subsection{Interpolation width and pointwise convergence}\label{max norm stuff}

As was discussed in Section \ref{discontinuity}, the discontinuity in the solution across the boundary prevents the method from achieving first-order convergence in the max norm near the boundary. As seen in Figure \ref{near boundary zoomed in}, the largest errors are confined to a region of about 2-3 meshwidths due to the support of the discrete delta function. However, there is also a wider region to which the numerical method spreads the errors. The \chIII{physical} length of the region on which the pointwise error fails to converge approaches $0$ as the grid is refined, so if \chIII{we only need} the solution away from the boundary, \chIII{we} can proceed with the method as illustrated, or even omit the interpolation step altogether. However, in this section, we will examine this issue more closely and investigate choices that can be made to recover pointwise convergence for the entire PDE domain.

%%%%%%%%%%%%%%%%%%%%%%%%%%%%%%%%%%%%%%%%%%%%%%%%%%
\subsubsection{Finite difference vs. Fourier spectral discretization}
%%%%%%%%%%%%%%%%%%%%%%%%%%%%%%%%%%%%%%%%%%%%%%%%%%

\begin{figure}
\centering
\begin{subfigure}{0.495\textwidth}
\centering
\includegraphics[width=\textwidth]{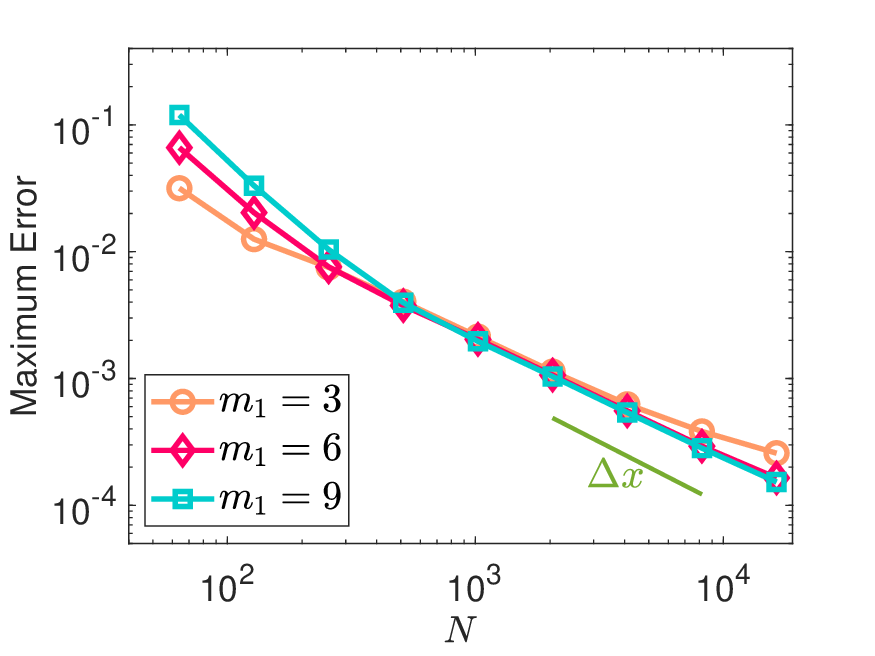}
\caption{\normalsize  Finite difference, Equation \eqref{pde again}}
\label{sin FD}
\end{subfigure}
\begin{subfigure}{0.495\textwidth}
\centering
\includegraphics[width=\textwidth]{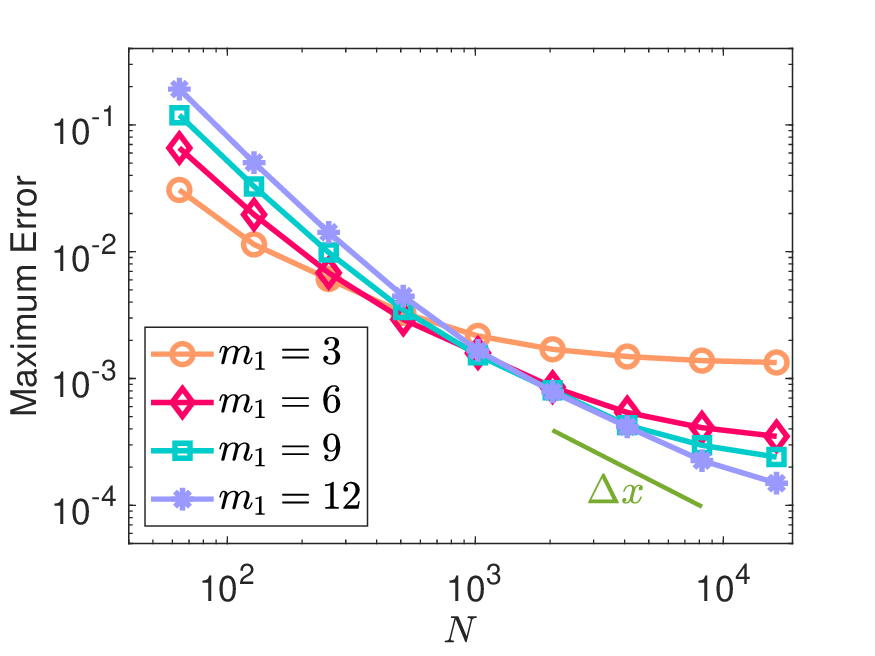}
\caption{\normalsize Fourier spectral, Equation \eqref{pde again} }
\label{sin spectral}
\end{subfigure}
\caption{The $L^{\infty}$ refinement studies for solving Equation \eqref{pde again} with the IBDL method. The computational domain is the periodic box $[-0.5, 0.5]^2$, $\Omega$ is the interior of a circle of radius 0.25, the boundary point spacing is $\Delta s \approx 0.75 \Delta x$, and $m_1$ is varied. Figure \ref{sin FD} is found using a finite difference method for discretization of the PDE, and Figure \ref{sin spectral} is found using a Fourier spectral method. }\label{FD vs spectral 1}
\end{figure}

We found that one key factor affecting the spread of error is the numerical method used to discretize the PDE operator. We therefore start by revisiting Equation \eqref{pde again} and utilizing both finite difference and Fourier spectral methods for discretization. We extend our refinement study to an even finer mesh of $N=2^{14}$ and explore different values for $m_1$, the number of meshwidths from the boundary for which interpolation is used. For all of these, we use an interior interpolation point that is $m_2=m_1+2$ meshwidths from the boundary. Figure \ref{FD vs spectral 1} shows the corresponding $L^{\infty}$ refinement studies. In Figure \ref{sin FD}, we can see that $m_1=6$ is sufficient to maintain first-order pointwise convergence when using a finite difference discretization. Using a Fourier spectral method, on the other hand, requires a larger value of $m_1$. For low, fixed values of $m_1$, we see that once the error is small enough, the Fourier spectral method gives a diminishing rate of pointwise convergence and eventually a failure to converge.

\begin{figure}
\centering
\begin{subfigure}{0.495\textwidth}
\centering
\includegraphics[width=\textwidth]{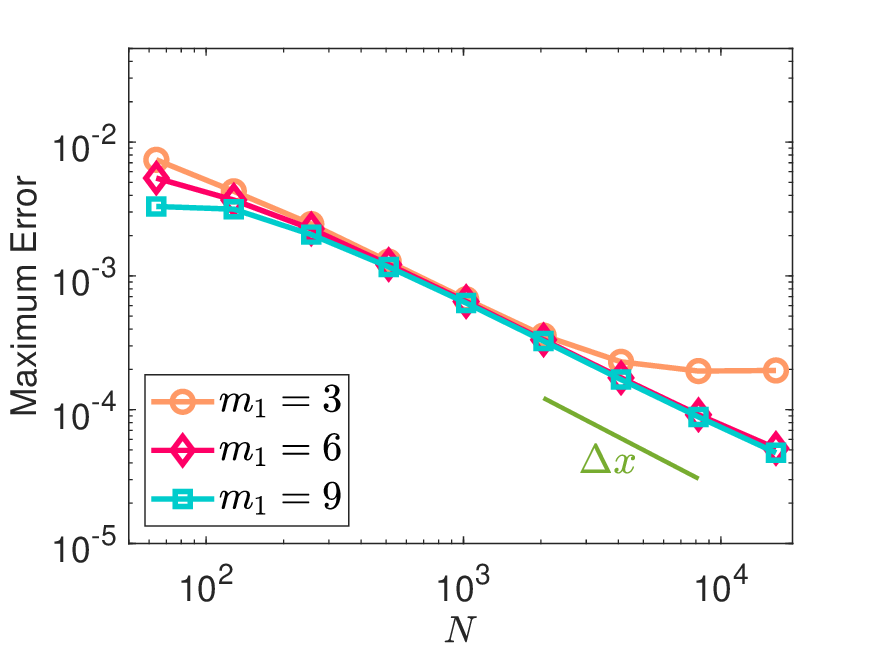}
\caption{\normalsize  Finite difference, Equation \eqref{pde again 2}}
\label{x FD}
\end{subfigure}
\begin{subfigure}{0.495\textwidth}
\centering
\includegraphics[width=\textwidth]{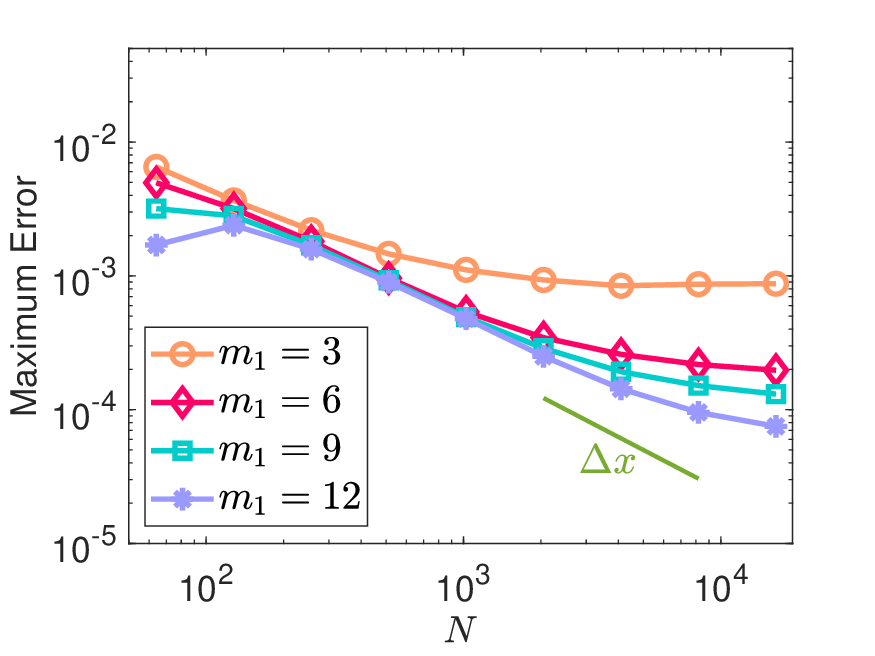}
\caption{\normalsize Fourier spectral, Equation \eqref{pde again 2} }
\label{x spectral}
\end{subfigure}
\begin{subfigure}{0.495\textwidth}
\centering
\includegraphics[width=\textwidth]{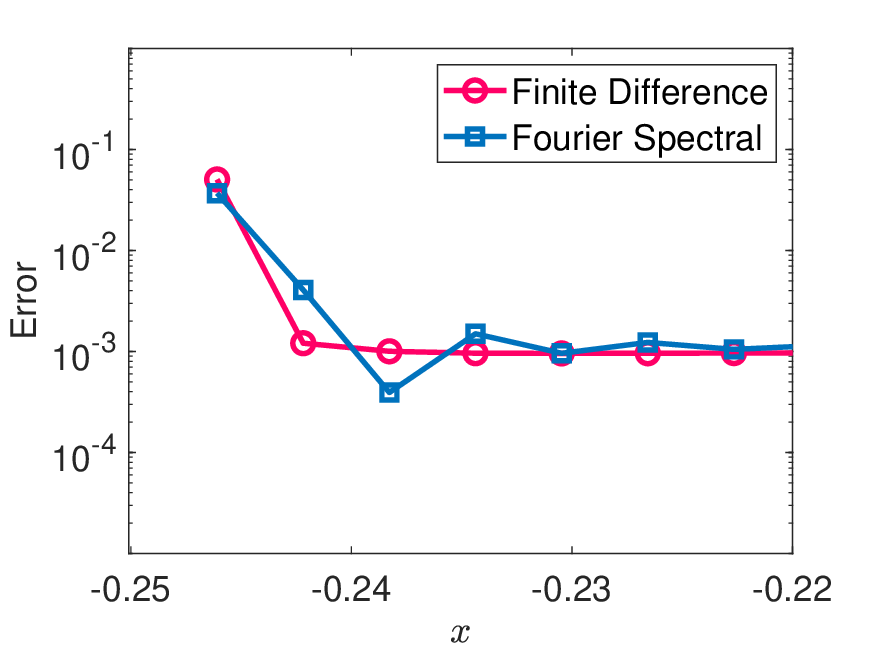}
\caption{\normalsize  Pointwise errors on coarse grid}
\label{coarse slice}
\end{subfigure}
\begin{subfigure}{0.495\textwidth}
\centering
\includegraphics[width=\textwidth]{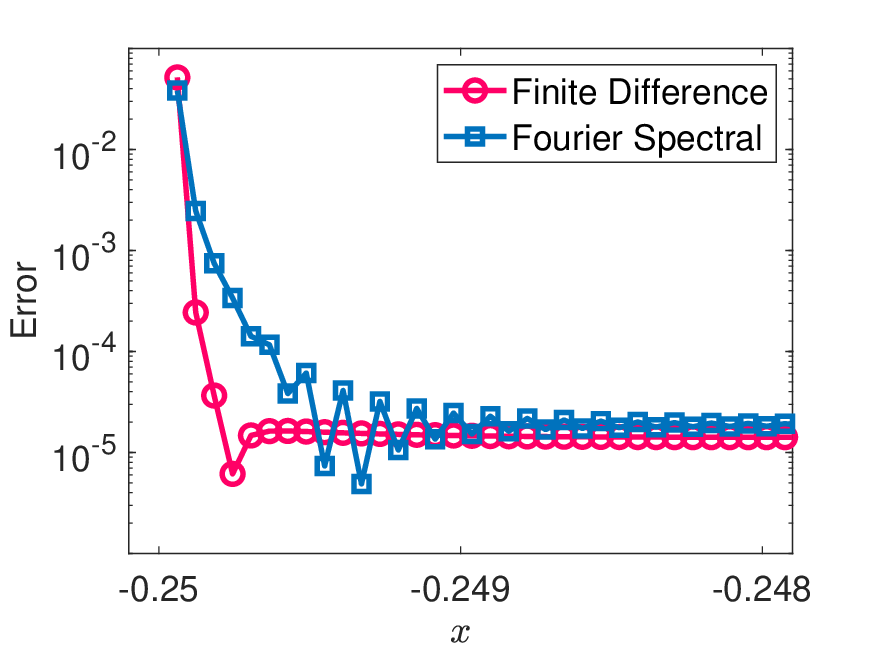}
\caption{\normalsize Pointwise errors on fine grid }
\label{fine slice}
\end{subfigure}
\caption{Figures \ref{x FD}-\ref{x spectral} give the $L^{\infty}$ refinement studies for solving Equation \eqref{pde again 2} with the IBDL method for varying values of $m_1$. The computational domain is the periodic box $[-0.5, 0.5]^2$, $\Omega$ is the interior of a circle of radius 0.25, and the boundary point spacing is $\Delta s \approx 0.75 \Delta x$. Figure \ref{x FD} is found using a finite difference method for discretization of the PDE, and Figure \ref{x spectral} is found using a Fourier spectral method. Figures \ref{coarse slice}-\ref{fine slice} show the errors along the line $y=0$ when solving Equation \eqref{pde again 2} with the IBDL method with no interpolation and using both finite difference and Fourier spectral methods for discretization of the PDE. Figure \ref{coarse slice} uses $N=2^8$, and Figure \ref{fine slice} uses $N=2^{14}$. Each plot marker represents a grid point on the corresponding mesh.}\label{FD vs spectral 2}
\end{figure}

To further explore this phenomenon, we next look at the PDE
\begin{subequations} \label{pde again 2}
\begin{alignat}{2}
& \Delta u -  u = -(x+y) \qquad && \text{in } \Omega  \label{pde1 again 2}\\
&u=x+y \qquad && \text{on } \Gamma,  \label{pde2 again 2}
\end{alignat}
\end{subequations}
where $\Omega$ is the interior of a circle of radius 0.25, centered at the origin. The analytic solution is given by $u=x+y$. Figures \ref{x FD}-\ref{x spectral} show the $L^{\infty}$ refinement studies for this problem. Due to the linearity of the solution, the primary source of error is the discontinuity, and this allows us to observe the aforementioned behavior on coarser grids. We can see that the lack of convergence with the Fourier spectral method is worse for this problem, but we see that $m_1=6$ is still sufficient to maintain convergence using a finite difference method. 

To explore this \ch{phenomenon} a little more closely, we also omit the interpolation step of the method and examine cross-sections of the error, located at $y=0$, for which the boundary point is located at $x=-0.25$. Figures \ref{coarse slice}-\ref{fine slice} show these error cross-sections on a relatively coarse grid of $N=2^8$ and a fine grid of $N=2^{14}$. Each marker on the plot represents a grid point on the mesh. For the coarse grid, we see in Figure \ref{coarse slice} that the pointwise errors for both methods decrease to the size of the interior error within a couple meshwidths. For the fine grid, however, the Fourier spectral method takes more meshwidths to decrease to the level of the interior error than the finite difference method. This illustrates the further spreading out of the error that takes place when using the Fourier spectral method due to the slow convergence of errors in a truncated Fourier series for a discontinuity. It is important to note, however, that while the number of meshwidths affected increases as the grid is refined, the physical distance affected still decreases. This is clear when observing the different scaling on the $x$-axes.

\begin{figure}
\centering
\begin{subfigure}{0.495\textwidth}
\centering
\includegraphics[width=\textwidth]{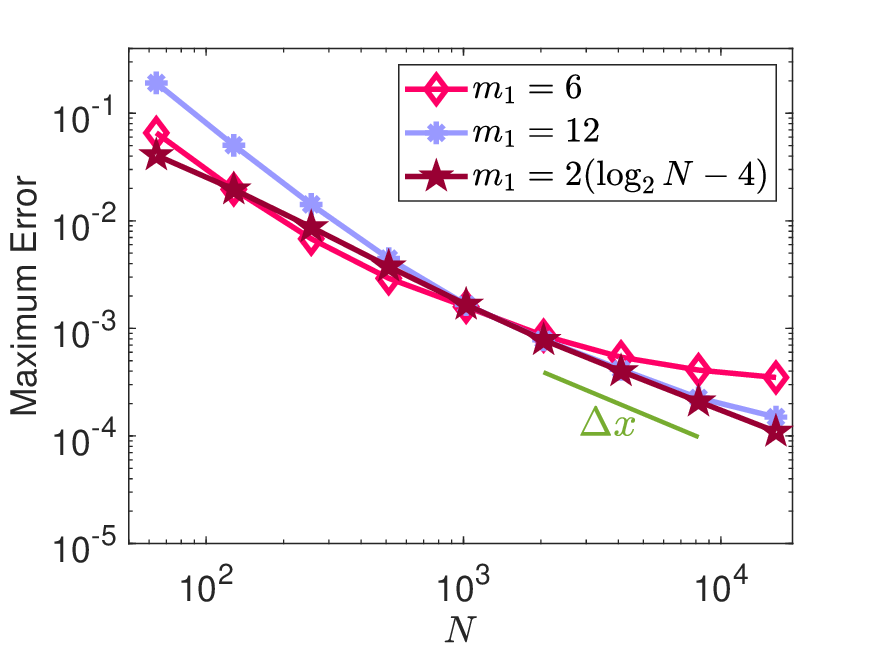}
\caption{\normalsize Fourier spectral, Equation \eqref{pde again}}
\label{sin spectral 2}
\end{subfigure}
\begin{subfigure}{0.495\textwidth}
\centering
\includegraphics[width=\textwidth]{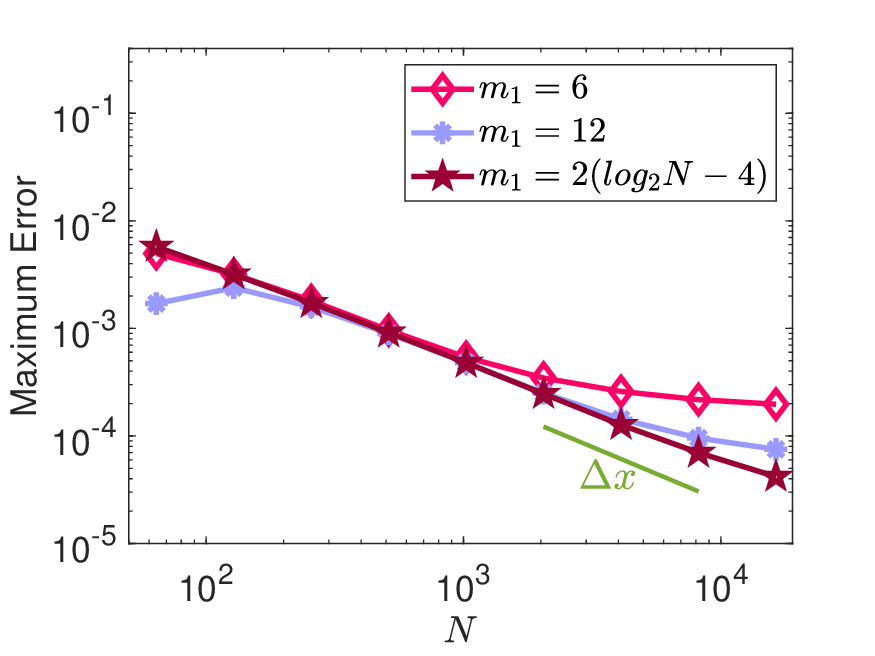}
\caption{\normalsize Fourier spectral, Equation \eqref{pde again 2}}
\label{x spectral 2}
\end{subfigure}
\caption{The $L^{\infty}$ refinement studies for solving Equation \eqref{pde again 2} using a Fourier spectral method for discretization of the PDE. The value of $m_1$ is again varied, but also included is an increasing function for $m_1$ given by $m_1=2(\log_2{(N)}-4)$. }\label{FD vs spectral 3}
\end{figure}

In order to recover pointwise convergence when utilizing a Fourier spectral method, \chIII{we} can choose a function for $m_1$ that increases as the mesh is refined. We have found that increasing logarithmically seems sufficient for these problems. For example, \chIII{we} could simply increase $m_1$ by 2 each time the number of grid points, $N$, is doubled. By doing this, we use interpolation at a greater number of meshwidths as we refine the grid, but, again, the physical distance on which we use interpolation is still decreasing relatively quickly towards 0. Figures \ref{sin spectral 2} and \ref{x spectral 2} illustrate using such a function for $m_1$ for Equations \eqref{pde again} and \eqref{pde again 2}, respectively. The function we use is $m_1=2(\log_2{N}-4)$, which gives $m_1$ values ranging from $ 4$ to $20$ on these grid sizes. We see that this option recovers pointwise convergence. 

%%%%%%%%%%%%%%%%%%%%%%%%%%%%%%%%%%%%%%%%%%%%%%%%%%
\subsubsection{Boundary point spacing}
%%%%%%%%%%%%%%%%%%%%%%%%%%%%%%%%%%%%%%%%%%%%%%%%%%

\begin{figure}
\centering
\begin{subfigure}{0.495\textwidth}
\centering
\includegraphics[width=\textwidth]{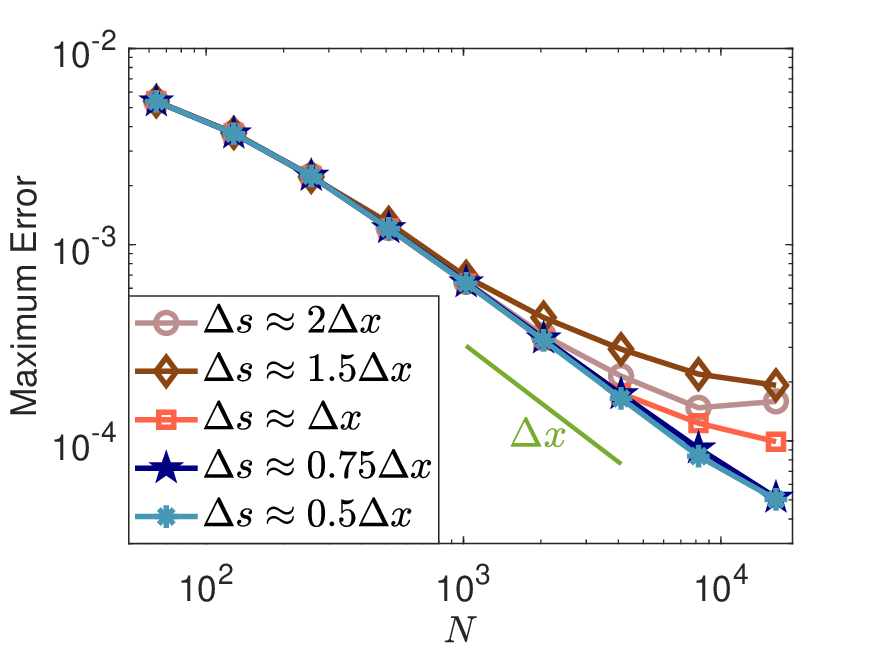}
\caption{\normalsize  Fixed $m_1=6$}
\label{x FD fixed m}
\end{subfigure}
\begin{subfigure}{0.495\textwidth}
\centering
\includegraphics[width=\textwidth]{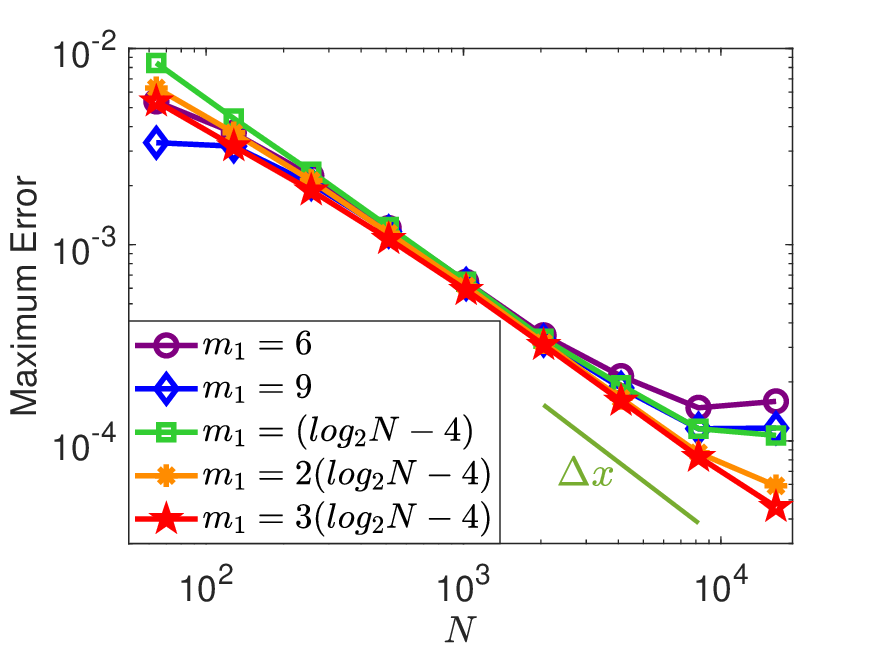}
\caption{\normalsize Fixed $\Delta s \approx 2 \Delta x$ }
\label{x FD fixed 2dx}
\end{subfigure}
\caption{The $L^{\infty}$ refinement studies for solving Equation \eqref{pde again 2} with the IBDL method using a finite difference discretization of the PDE. Figure \ref{x FD fixed m} uses a fixed interpolation width of $m_1=6$ meshwidths and a second interpolation point $m_2=8$ meshwidths from the boundary, and the boundary point spacing is varied. Figure \ref{x FD fixed 2dx} has a wide boundary point spacing of  $\Delta s \approx 2 \Delta x$, and the interpolation width $m_1$ is varied.}\label{FD point spacing}
\end{figure}

In addition to the PDE discretization and interpolation width, one other factor affecting the error spread is the boundary point spacing. We have observed that in general, using tightly spaced boundary points gives better pointwise convergence for a fixed interpolation width, especially in the case of the finite difference method. Specifically, we have seen that keeping $\Delta s \leq 0.75 \Delta x$ will help to avoid a slowing pointwise convergence rate. This is illustrated in Figure \ref{x FD fixed m}, which gives the $L^{\infty}$ refinement studies for Equation \eqref{pde again 2} using a finite difference discretization, $m_1=6$, and various the boundary point spacings. \chI{A conventional choice for IB methods is to choose $\Delta s \approx 0.5 \Delta x$ in order to avoid fluid leakage through the immersed boundary \cite{Peskin02}. This constraint is typically considered to be less important for rigid structures since the boundary points do not move relative to each other; therefore those utilizing an IBSL-like method typically choose $\Delta s \approx \Delta x$ or $\Delta s \approx 2 \Delta x$ to control conditioning \cite{Taira, GriffithDonev}. However, in the IBDL method, since we saw in Section \ref{interior circle} that the iteration count remains essentially unchanged by tightening the boundary points, there is no numerical drawback to using this tighter point spacing. }

If \chIII{we do} want to use more widely spaced points, \chIII{we} can simply increase the interpolation widths \chIII{from} the previous section. Figure \ref{x FD fixed 2dx} shows the $L^{\infty}$ refinement studies for Equation \eqref{pde again 2} when the boundary point spacing is given by $\Delta s \approx 2 \Delta x$ and the interpolation width is varied. In this case of more widely spaced boundary points, increasing $m_1$ logarithmically as the mesh is refined again recovers pointwise convergence.

In accordance with the observations seen in Section \ref{max norm stuff}, for the majority of this paper, we have elected to use a finite difference PDE discretization, a boundary point spacing of $\Delta s \approx 0.75 \Delta x$, and an interpolation width of $m_1=6$ meshwidths. Similar results can be found by using a Fourier spectral method, a boundary point spacing of $\Delta s \approx 0.75 \Delta x$, and an interpolation width given by $m_1=2(\log_2{(N)}-4)$. We reiterate that using these options is not essential if \chIII{we are} only concerned with the solution away from the boundary because the width of the region on which the error does not converge pointwise does approach $0$ as the grid is refined.

%%%%%%%%%%%%%%%%%%%%%%%%%%%%%%%%%%%%%%%%%%%%%%%%%%
%%%%%%%%%%%%%%% 6.3 POISSON%%%%%%%%%%%%%%%%%%%%%%%%
%%%%%%%%%%%%%%%%%%%%%%%%%%%%%%%%%%%%%%%%%%%%%%%%%%
\subsection{\chIII{Poisson equation and completed IBDL method}}\label{poisson1}

\begin{figure}
\centering
\begin{subfigure}{0.5\textwidth}
\centering
\includegraphics[width=\textwidth]{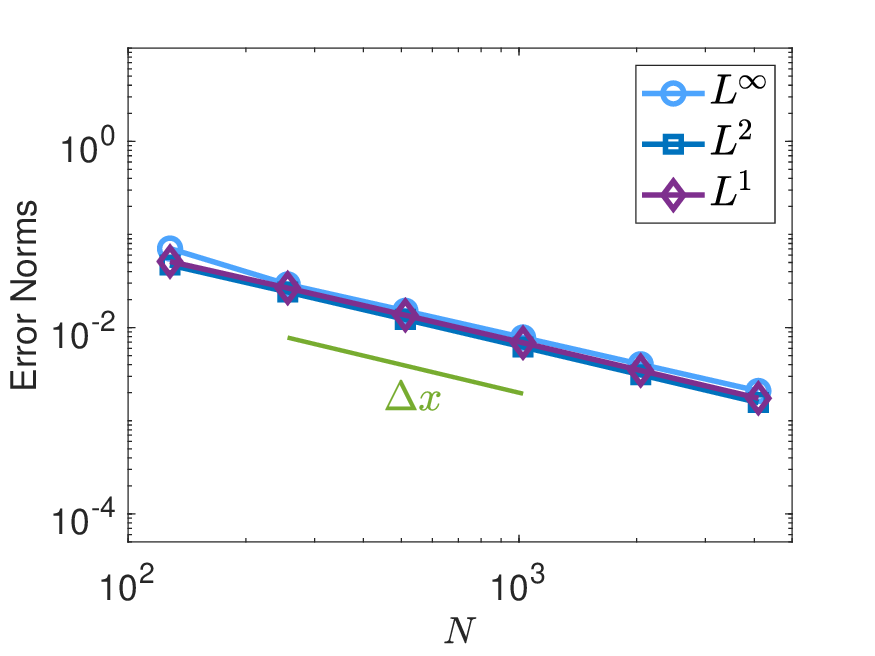}
\caption{\normalsize  IBDL, $L=1$}
\label{smallLexp}
\end{subfigure}
\begin{subfigure}{0.495\textwidth}
\centering
\includegraphics[width=\textwidth]{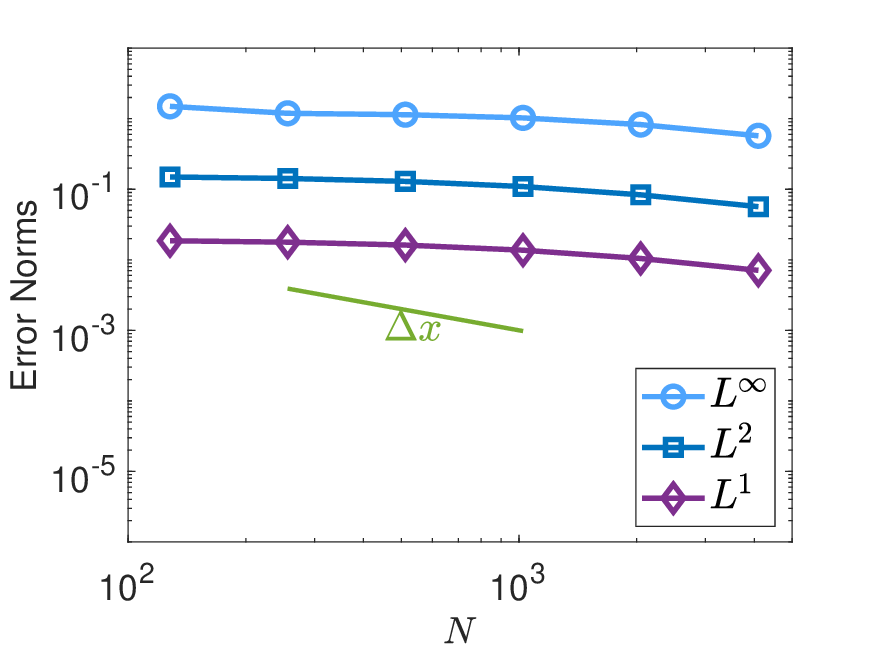}
\caption{\normalsize  IBDL, $L=8$}
\label{bigLexp}
\end{subfigure}
\begin{subfigure}{0.495\textwidth}
\includegraphics[width=\textwidth]{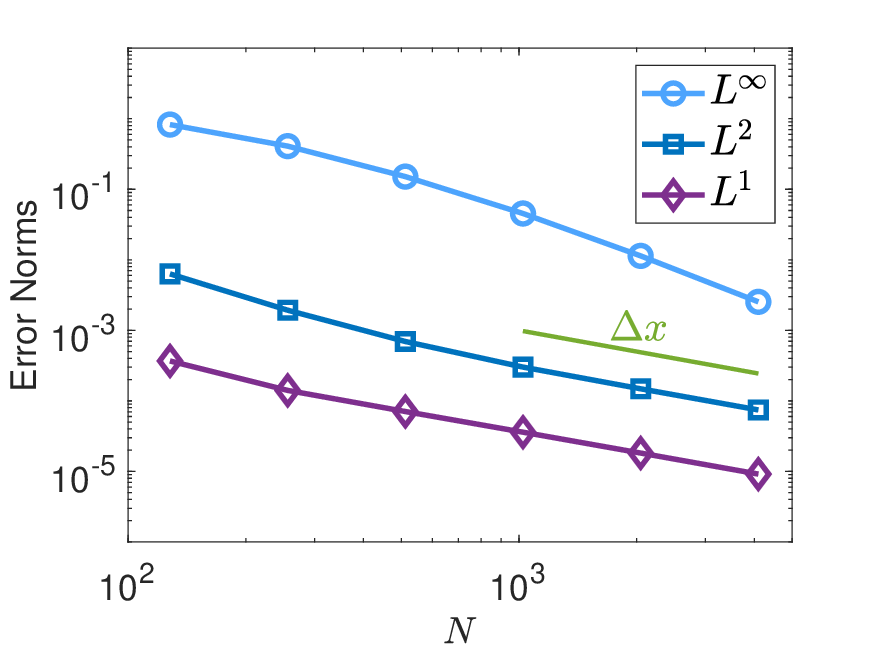}
\caption{\normalsize  Completed IBDL, $L=8$}
\label{bigLexpcomp}
\end{subfigure}
\caption{\chIII{Refinement studies for solutions to Equation \eqref{poisson shit pde exp} found using the IBDL method, with finite differences. For Figure \ref{smallLexp}, the computational domain is the periodic box $[-0.5, 0.5]^2$, and for Figures \ref{bigLexp}-\ref{bigLexpcomp}, the computational domain is the periodic box $[-4, 4]^2$. Figures \ref{smallLexp}-\ref{bigLexp} use the original IBDL method, and Figure \ref{bigLexpcomp} uses the completed IBDL method. $\Omega$ is the region in the computational domain exterior to a circle of radius $0.25$, centered at the origin. A boundary point spacing of $ \Delta s \approx 0.75  \Delta x$ is used, and solution values within $m_1= 6$ meshwidths from the boundary are replaced using a second interpolation point $m_2=8$ meshwidths away from the boundary. All errors reported are absolute errors.}}\label{poisson shit plots}
\end{figure}

\chIII{In this section, we demonstrate the use of the completed IBDL method for large exterior domains, formulated in Section \ref{poisson explanation}. To do this, we examine the PDE
\begin{subequations} \label{poisson shit pde exp}
\begin{alignat}{2}
& \Delta u  =  4\pi^2 e^{\sin{(2 \pi x )}}\Big(\cos^2{(2\pi x)}-\sin{(2\pi x)}\Big) \qquad && \text{in } \Omega \\
&u= e^{\sin{(2\pi x)}}&& \text{on } \Gamma, 
\end{alignat}
\end{subequations}
for which the analytical solution is given by $u = e^{\sin{(2\pi x)}}$. We take the PDE domain $\Omega$ to be the region in the periodic box $[-L/2, L/2]^2$ that is exterior to a circle of radius $0.25$, centered at the origin. We look at both $L=1$ and $L=8$. We use finite differences and interpolate values within $m_1=6$ meshwidths of the boundary. We first solve these problems using the IBDL method described in Equation \eqref{ibdl}, without the addition of the single layer potential. Figure \ref{smallLexp} illustrates that on the small periodic domain of length $1$, the computed solution has first-order convergence. However, Figure \ref{bigLexp} demonstrates that for a larger periodic domain of length $8$, the IBDL method is unable to achieve convergence to the analytical solution to Equation \eqref{poisson shit pde exp}. We next apply the completed IBDL method, described in Equation \eqref{ibdl completed}, with $\eta=10$ for this large computational domain. Figure \ref{bigLexpcomp} illustrates we get first order convergence with the completed method. }

%%%%%%%%%%%%%%%%%%%%%%%%%%%%%%%%%%%%%%%%%%%%%%%%%%%%%%%%%%%%%%%%%%%%%%%%%%%%%%%%%%%%%%%%%%%%%%%%%%%%%%%%%%%

\subsection{\chIII{Non-convex boundary}}\label{poisson2}

\begin{figure}
\centering
\begin{subfigure}{0.495\textwidth}
\centering
\includegraphics[width=\textwidth]{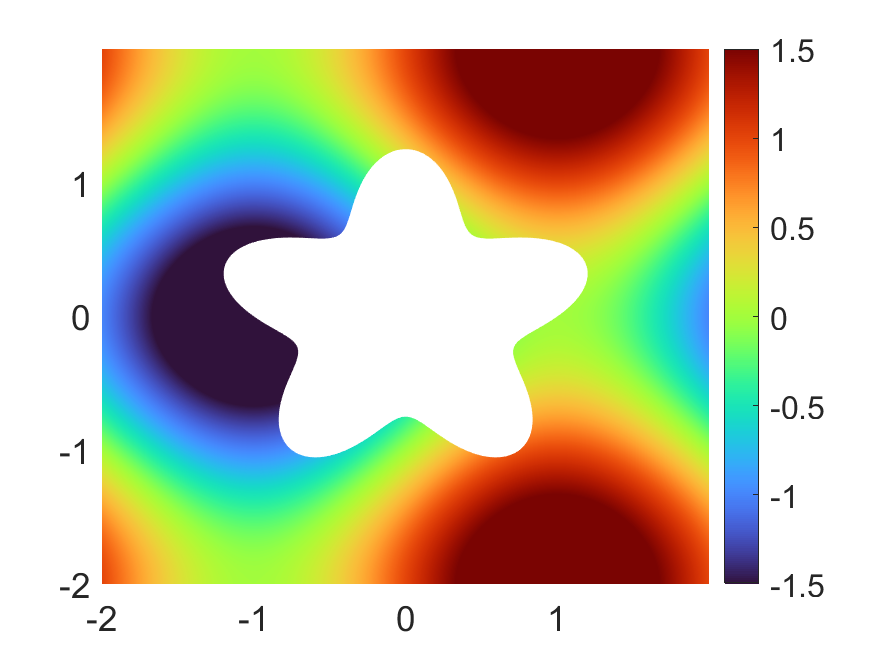}
\caption{\normalsize  IBDL Solution}
\label{ibdl poisson solution}
\end{subfigure}
\begin{subfigure}{0.495\textwidth}
\centering
\includegraphics[width=\textwidth]{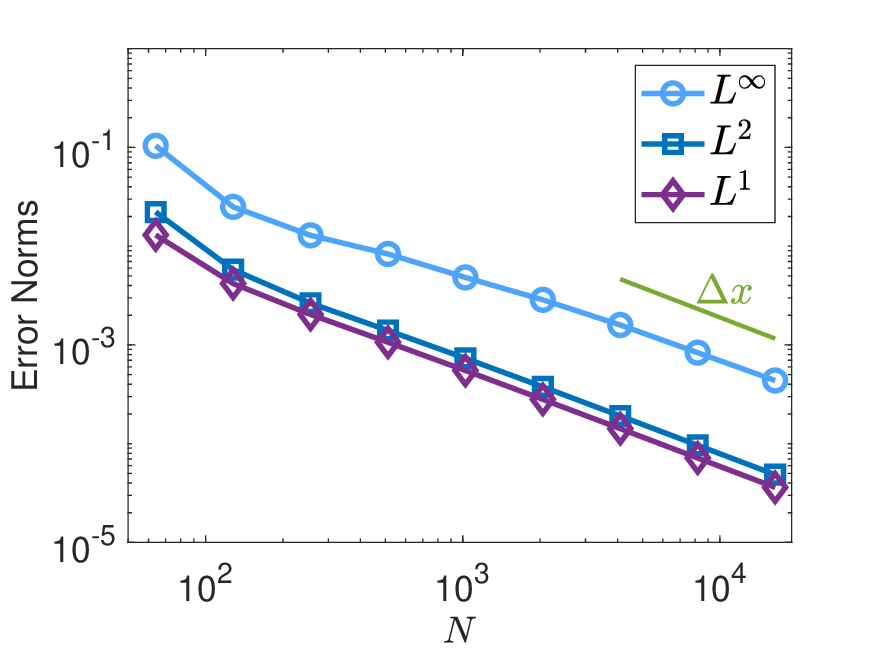}
\caption{\normalsize  IBDL Refinement}
\label{ibdl poisson refinement}
\end{subfigure}
\caption{A solution plot and refinement study for solving Equation \eqref{pde again 3} with the IBDL method. The computational domain is the periodic box $[-2, 2]^2$, $\Omega$ is the region exterior to a ``starfish'' shape, and the prescribed boundary values are given by $U_b=\sin{(\pi x/2)} - \cos{(\pi x/2)}$. A boundary point spacing of $ \Delta s \approx 0.75  \Delta x$ is used, and solution values within $m_1= 6$ meshwidths from the boundary are replaced using a second interpolation point $m_2=8$ meshwidths away from the boundary.}\label{poisson plots}
\end{figure}

\chIII{In this section, we demonstrate the effectiveness of the IBDL method for a non-convex boundary. To do this, we use the original IBDL method, without the completion, to solve the PDE }
\begin{subequations} \label{pde again 3}
\begin{alignat}{2}
& \Delta u  =\frac{\pi^2}{4}\Bigg( \cos{\Big(\frac{\pi y}{2}\Big)}-\sin{\Big(\frac{\pi x}{2}\Big)}\Bigg)\qquad && \text{in } \Omega  \label{pde1 again 3}\\
&u=\sin{\Big(\frac{\pi x}{2}\Big)} - \cos{\Big(\frac{\pi y}{2}\Big)} && \text{on } \Gamma,  \label{pde2 again 3}
\end{alignat}
\end{subequations}
where $\Omega$ is the region inside the periodic box $[-2,2]^2$ that is exterior to the `starfish'' shape used in \cite{QBX}, which is given by
\begin{equation} \label{starfish eqn}
\begin{pmatrix}
x(\theta)\\
y(\theta)
\end{pmatrix}
= 
\Bigg( 1+\frac{\sin{(10\pi \theta)}}{4}\Bigg)\begin{pmatrix}
\cos{(2\pi \theta)}\\
\sin{(2\pi \theta)}
\end{pmatrix}, 
\end{equation}
for $0 \leq \theta \leq 1$. We use equally spaced boundary points, with $\Delta s \approx 0.75 \Delta x$, and the solutions are computed for grid sizes ranging from $N=2^6$ to $2^{14}$.  The analytic solution is given by $ u=\sin{(\pi x/2)} - \cos{(\pi y/2)}$. For this problem, the IBDL method only requires  $13-14$ iterations of \texttt{gmres}. Figure \ref{ibdl poisson solution} shows a plot of the solution obtained from the IBDL method, and Figure \ref{ibdl poisson refinement} demonstrates the first-order convergence in the refinement study.

%%%%%%%%%%%%%%%%%%%%%%%%%%%%%%%%%%%%%%%%%%%%%%%%%%
%%%%%%%%%%%%%%%%%% F vs Q   %%%%%%%%%%%%%%%%%%%%%%
%%%%%%%%%%%%%%%%%%%%%%%%%%%%%%%%%%%%%%%%%%%%%%%%%%

\subsection{Convergence of potential strength }\label{Q convergence}

In addition to the increased speed with which we solve the system for the potential strength, $Q$, another benefit of the better conditioning of the IBDL method is that we are more easily able get convergence in this potential strength. It is a well-known problem that the poor conditioning of the Schur complement results in a noisy force distribution in the IBSL method \cite{Goza}, and this noise is larger for tighter boundary point spacing. Here, we illustrate that in the IBDL method, we obtain much smoother distributions with little need for filtering. 

To demonstrate this, we compute the IBSL and IBDL force distributions, $F$ and $Q$, obtained by solving the interior problem from Section \ref{interior circle}, where the computational domain is instead the periodic box $[-10,10]^2$. For estimating convergence, we approximate the exact distributions using a boundary element method, and this larger computational domain allows us to approximate the periodic Green's function with the free-space Green's function. 

\begin{figure}
\centering
\begin{subfigure}{0.41\textwidth}
\centering
\includegraphics[width=\textwidth]{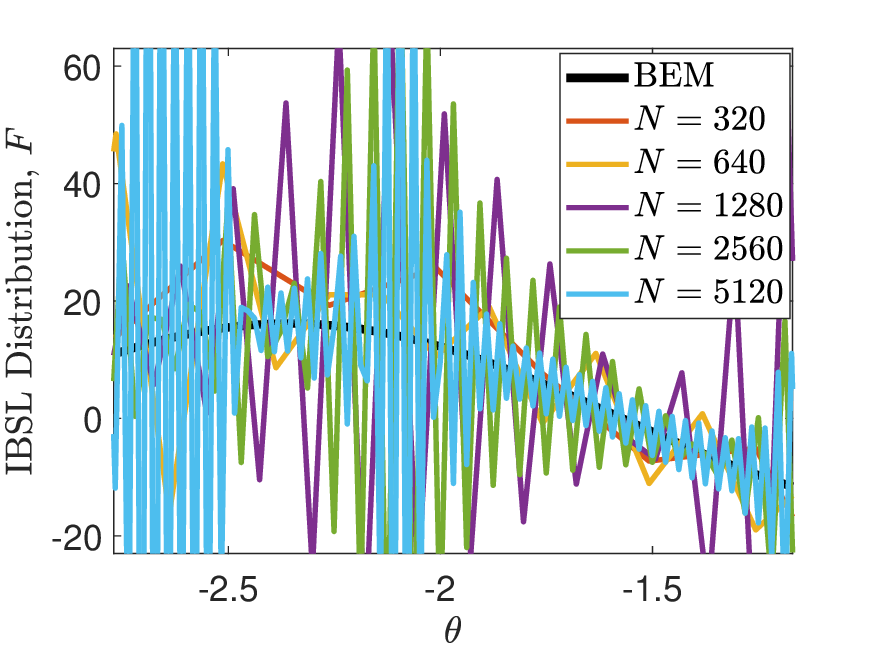}
\caption{\normalsize IBSL, $\Delta s \approx \Delta x$}
\label{ibsl f plot ds1}
\end{subfigure}
\begin{subfigure}{0.41\textwidth}
\centering
\includegraphics[width=\textwidth]{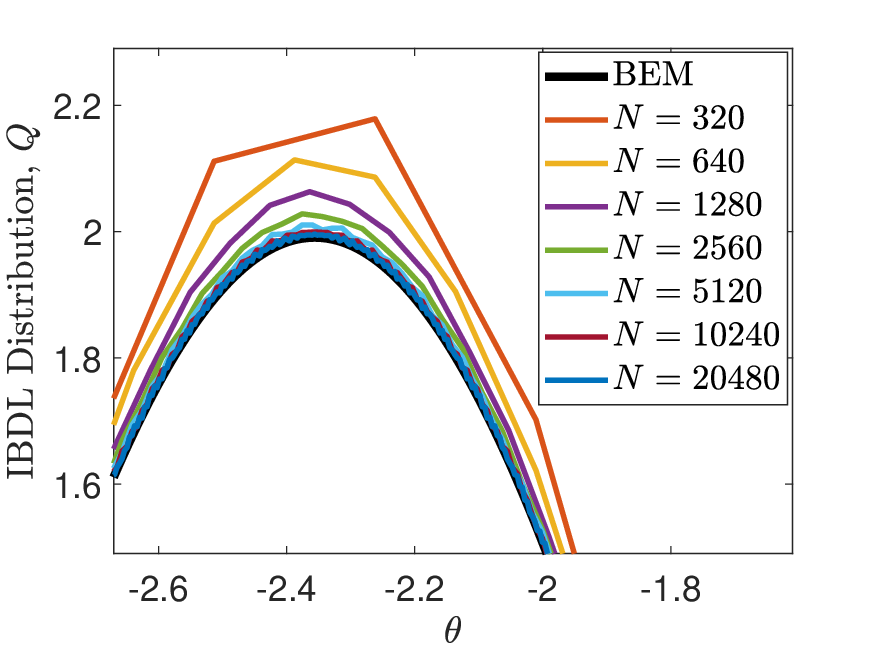}
\caption{\normalsize IBDL, $\Delta s \approx \Delta x$}
\label{ibdl q plot ds1}
\end{subfigure}
\begin{subfigure}{0.41\textwidth}
\centering
\includegraphics[width=\textwidth]{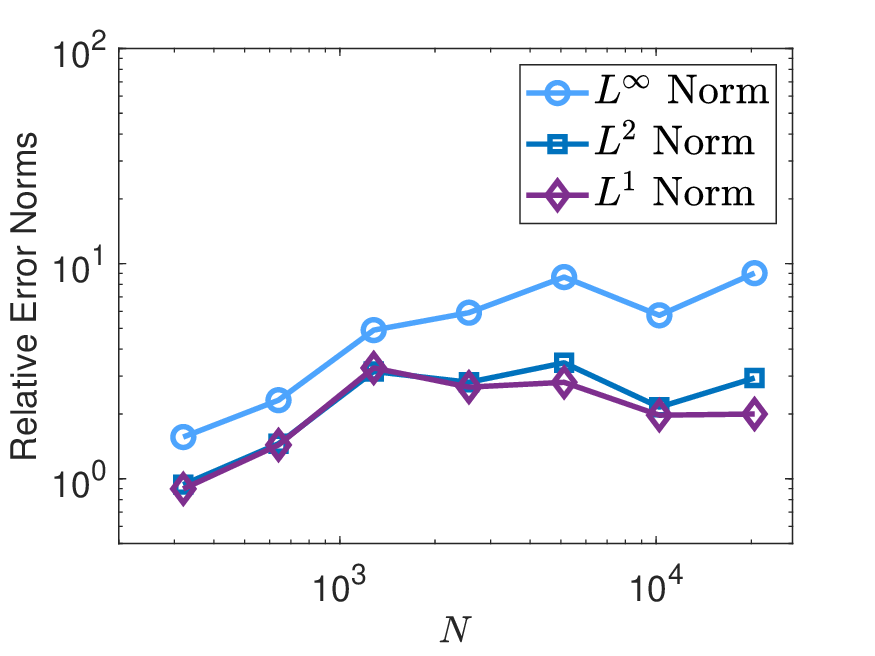}
\caption{\normalsize IBSL, $\Delta s \approx \Delta x$}
\label{ibsl f refinement}
\end{subfigure}
\begin{subfigure}{0.41\textwidth}
\centering
\includegraphics[width=\textwidth]{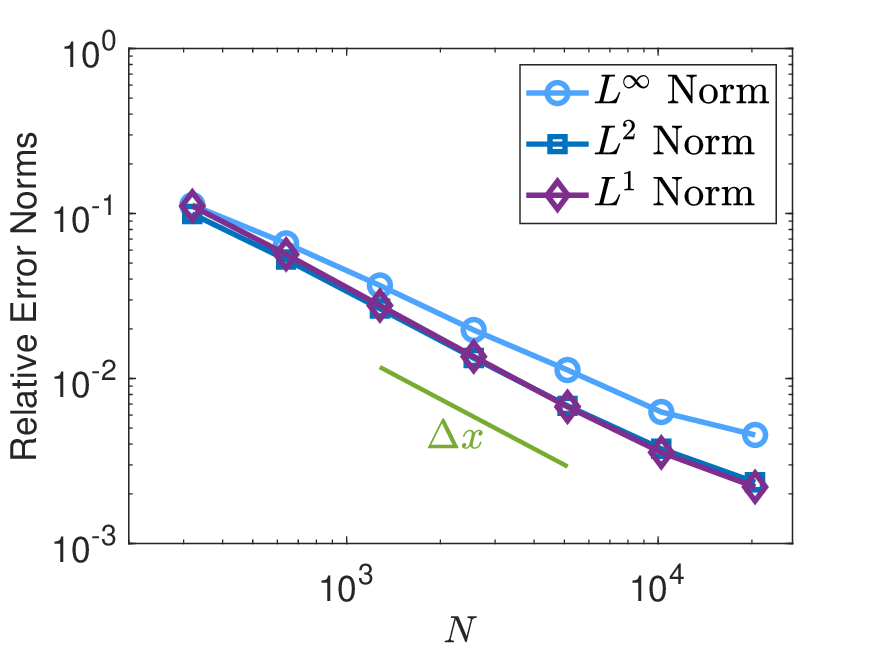}
\caption{\normalsize IBDL, $\Delta s \approx \Delta x$}
\label{ibdl q refinement}
\end{subfigure}
\begin{subfigure}{0.41\textwidth}
\centering
\includegraphics[width=\textwidth]{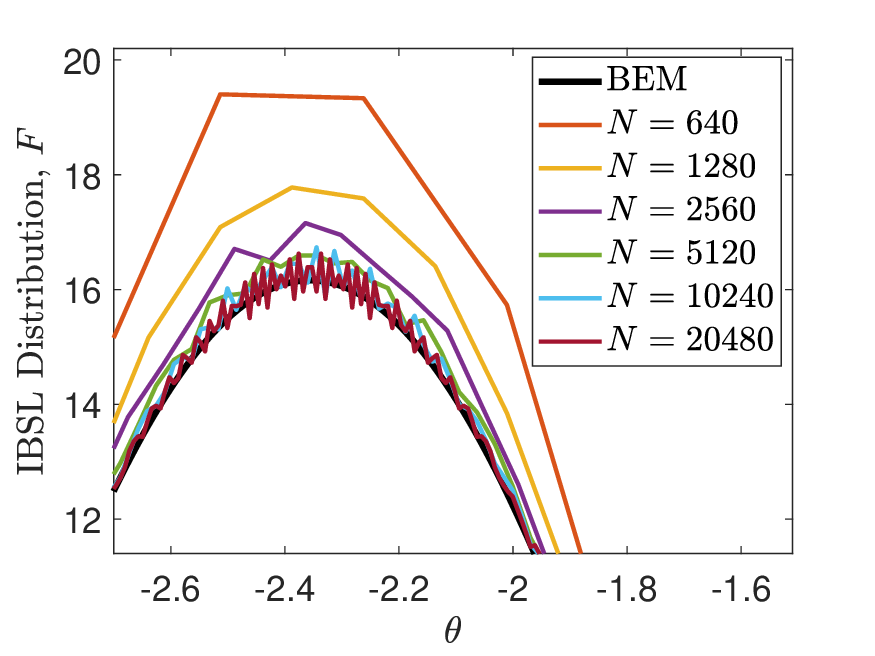}
\caption{\normalsize IBSL, $\Delta s \approx 2\Delta x$}
\label{ibsl f plot ds2}
\end{subfigure}
\begin{subfigure}{0.41\textwidth}
\centering
\includegraphics[width=\textwidth]{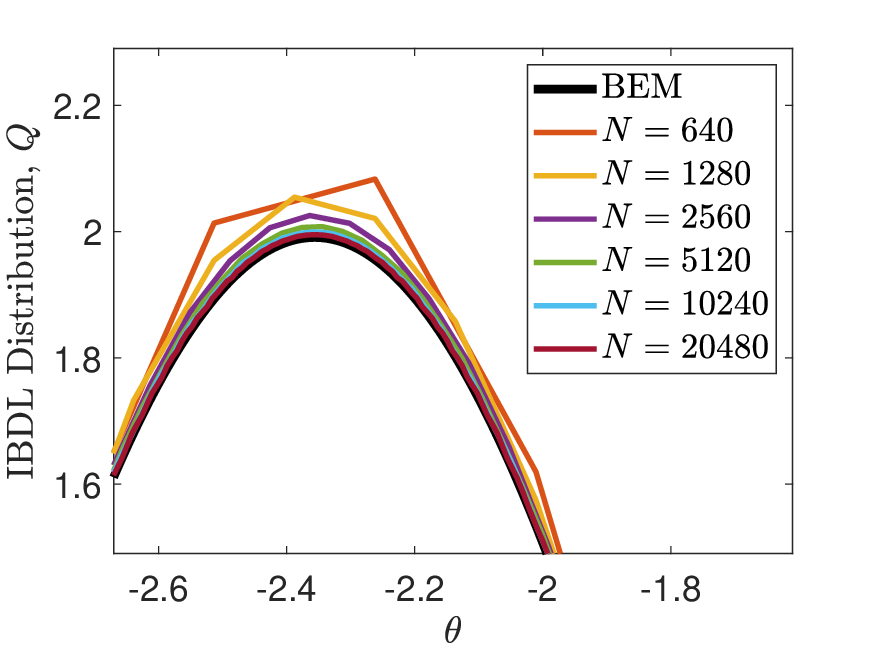}
\caption{\normalsize IBDL, $\Delta s \approx2 \Delta x$}
\label{ibdl q plot ds2}
\end{subfigure}
\begin{subfigure}{0.41\textwidth}
\centering
\includegraphics[width=\textwidth]{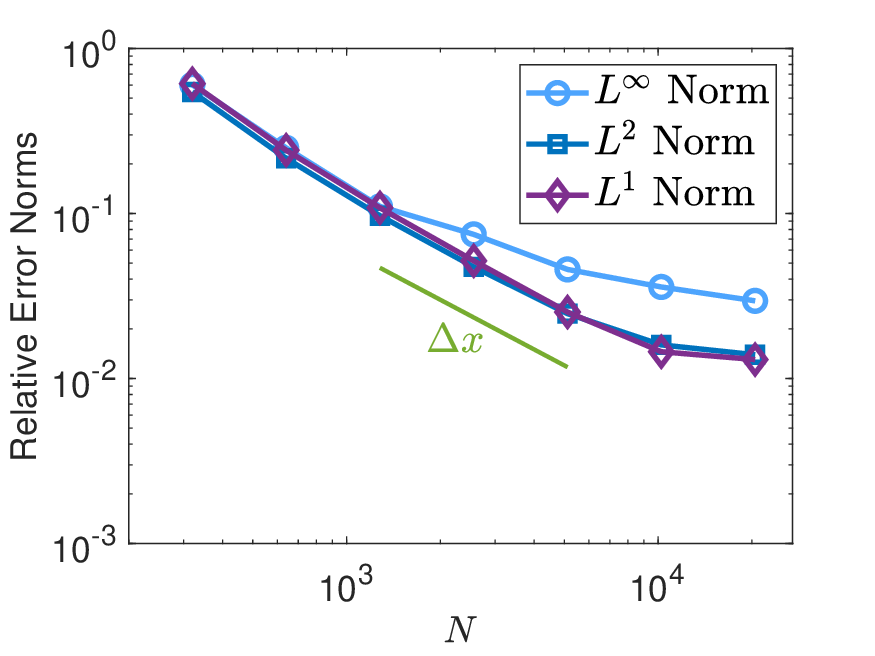}
\caption{\normalsize IBSL, $\Delta s \approx 2\Delta x$}
\label{ibsl f refinement ds2}
\end{subfigure}
\begin{subfigure}{0.41\textwidth}
\centering
\includegraphics[width=\textwidth]{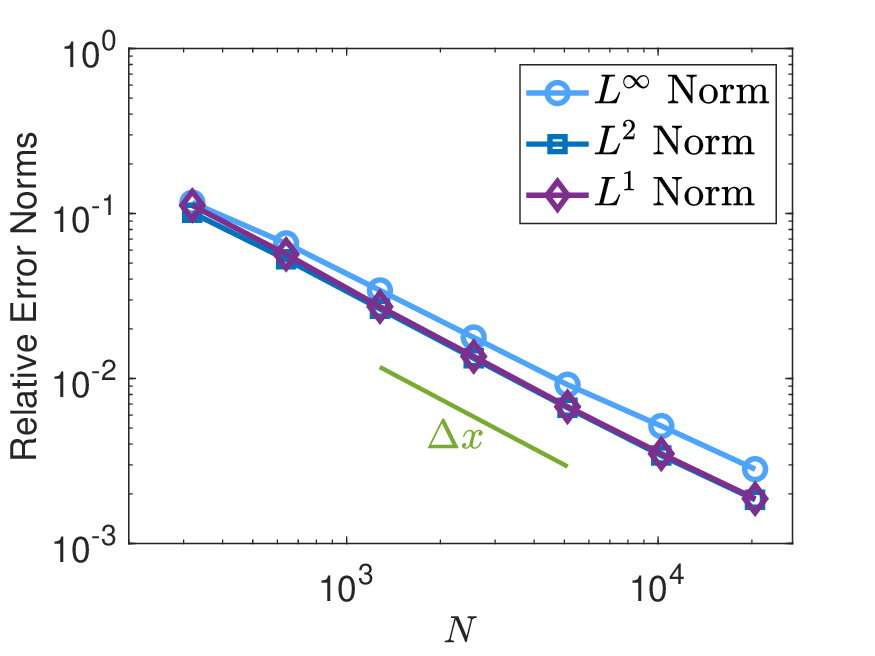}
\caption{\normalsize IBDL, $\Delta s \approx2 \Delta x$}
\label{ibdl q refinement ds2}
\end{subfigure}
\caption{Figures \ref{ibsl f plot ds1}-\ref{ibdl q plot ds1} and  \ref{ibsl f plot ds2}-\ref{ibdl q plot ds2} show portions of the the IBSL and IBDL potential strengths, plotted against the angle $\theta$, for Equation \eqref{pde again}.  Plots are shown for different mesh sizes, and the boundary point spacing is given by $\Delta s \approx \Delta x$ for Figures \ref{ibsl f plot ds1}-\ref{ibdl q plot ds1} and $\Delta s \approx 2\Delta x$ for Figures \ref{ibsl f plot ds2}-\ref{ibdl q plot ds2} . The black curves correspond to the distributions found using a boundary element method with uniform, straight line elements and Gauss-Legendre quadrature. These BEM distributions are used to approximate the error, normalized by the maximum value, and Figures \ref{ibsl f refinement}-\ref{ibdl q refinement} and \ref{ibsl f refinement ds2}-\ref{ibdl q refinement ds2} show the corresponding refinement studies. }\label{Q and F stuff }
\end{figure}

In Figures \ref{ibsl f plot ds1}-\ref{ibdl q refinement}, for which we use $\Delta s \approx \Delta x$, we see approximately first-order convergence for the IBDL distribution $Q$ and no convergence for the IBSL distribution $F$. Figure \ref{ibsl f plot ds1} shows that $F$ quickly becomes very noisy as we refine the grid. As seen in Figure \ref{ibdl q plot ds1}, there is much less noise in $Q$.  Note that the plots only show a portion of the distributions in order to more easily view this noise. 

For Figures \ref{ibsl f plot ds2}-\ref{ibdl q refinement ds2}, we increase the spacing for the boundary points to $\Delta s \approx 2 \Delta x$. Then we see some convergence in the IBSL distribution $F$. However, the IBDL potential strength $Q$ achieves greater convergence and a much lower relative error than $F$. Therefore, we see that the IBDL method in general allows us to achieve a much more accurate distribution, and the next section demonstrates an application for which this convergence is vital.

%%%%%%%%%%%%%%%%%%%%%%%%%%%%%%%%%%%%%%%%%%%%%%%%%%
%%%%%%%%%%%%%%%%%% Neumann   %%%%%%%%%%%%%%%%%%%%%%
%%%%%%%%%%%%%%%%%%%%%%%%%%%%%%%%%%%%%%%%%%%%%%%%%%

\subsection{Neumann boundary conditions}\label{neumann}

As discussed in Section \ref{Neumann formulation}, we are able to use the connection to integral equations to solve a PDE with Neumann boundary conditions by solving 
\begin{equation}
(-S^* \L^{-1}\widetilde S)U_b-\frac12 U_b=S^*\L^{-1}SV_b-S^*\L^{-1}\chII{\tilde g} \label{Neumann saddle point again}
\end{equation} 
for the unknown boundary values $U_b$. By doing this, we are again inverting an operator corresponding to a second-kind integral operator with a small condition number. Therefore, we are able to get convergence in this boundary distribution, as discussed in the previous section. When we have Neumann boundary conditions, the convergence of the distribution, $U_b$ is especially important because $U_b$ is used when performing our interpolation step for grid points near the boundary. In the Dirichlet case, $U_b$ is known, providing us with more accurate interpolation data, but in the case of Neumann, the values both in the interior and on the boundary are approximate. 

\begin{figure}
\centering
\begin{subfigure}{0.495\textwidth}
\centering
\captionsetup{justification=centering}
\includegraphics[width=\textwidth]{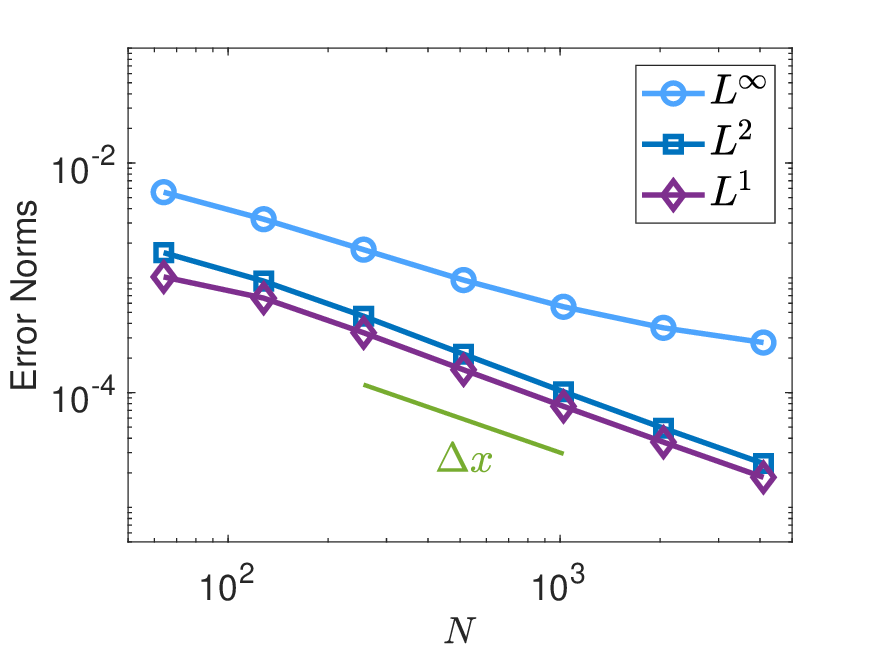}
\caption{\normalsize  Solution refinement study}
\label{neumann no filter 075}
\end{subfigure}
\begin{subfigure}{0.495\textwidth}
\centering
\captionsetup{justification=centering}
\includegraphics[width=\textwidth]{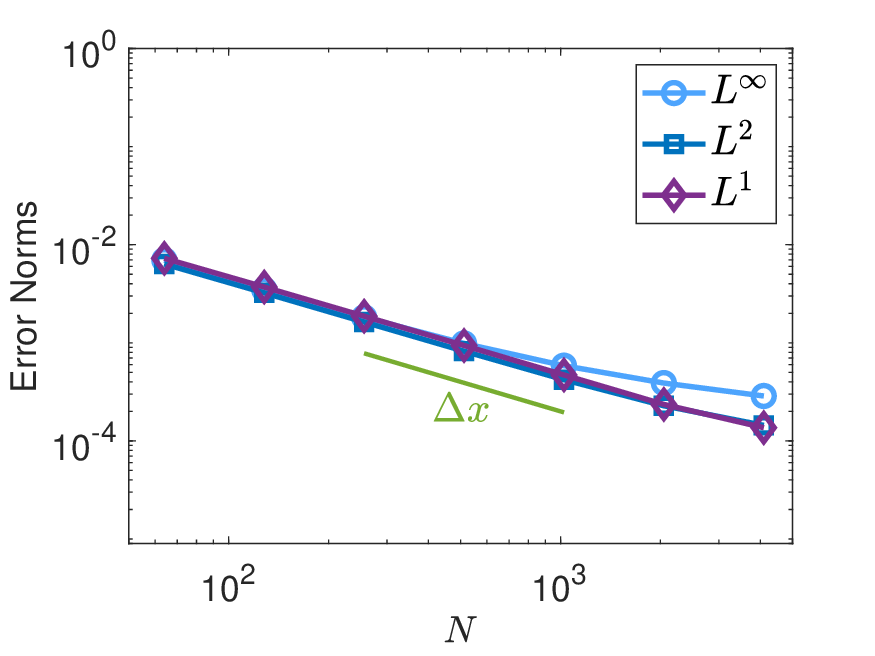}
\caption{\normalsize  $U_b$ refinement study}
\label{neumann ub no filter 075}
\end{subfigure}
\caption{Refinement studies for solving Equation \eqref{pde again 4} with the method discussed in Section \ref{Neumann formulation}. The computational domain is the periodic box $[-0.5, 0.5]$, $\Omega$ is the interior of a circle of radius 0.25, and the prescribed normal derivatives on the boundary are given by $\partial u/\partial n|_{\Gamma}= 8(x^2-y^2)$. We us a boundary point spacing of $\Delta s \approx 0.75 \Delta x$ and replace solution values within $m_1 = 6$ meshwidths from the boundary using interior interpolation points $m_2=8$ meshwidths away from the boundary. Figures \ref{neumann no filter 075} and \ref{neumann ub no filter 075} show the refinement studies for the solution errors and the boundary values $U_b$, respectively.}\label{Neumann plots}
\end{figure}

We consider the PDE given by 
\begin{subequations} \label{pde again 4}
\begin{alignat}{2}
& \Delta u -  u =-(x^2-y^2)\hspace{0.2cm} && \text{in } \Omega  \label{pde1 again 4}\\
&\frac{\partial u}{\partial n}= 8(x^2-y^2)&& \text{on } \Gamma,  \label{pde2 again 4}
\end{alignat}
\end{subequations}
where $\Omega$ is the interior of a circle of radius 0.25, centered at the origin. The analytic solution is given by $u=x^2-y^2$. Our computational domain here is the periodic box $[-0.5,0.5]^2$, our boundary point spacing is given by $\Delta s \approx0.75 \Delta x$, and we replace solution values within $m_1=6$ meshwidths of the boundary using interior interpolation points located $m_2=8$ meshwidths from the boundary. For this problem, our method only requires $4-6$ iterations of \texttt{gmres}. Figure \ref{Neumann plots} illustrates the refinement studies for the solution values and $U_b$, and we see approximately first-order convergence in both of these.

%%%%%%%%%%%%%%%%%%%%%%%%%%%%%%%%%%%%%%%%%%%%%%%%%%%%%%%%%%%%%%% 3 D %%%%%%%%%%%%%%%%%%%%%%%%%%%%%%%%%%%%%%%%%%%%%%%%%%%%%%%%%%%%%%%%%%%%%%%%%%%%%%%%%%%%%%%%%%%%%%

\subsection{\cha{Dirichlet Helmholtz in 3-D}}\label{3D}
\cha{In order to show that this method can be applied to three dimensions, we now apply the method to the PDE 
\begin{subequations} \label{3d pde}
\begin{alignat}{2}
& \Delta u -  u = (2-r^2)\cos{\theta} \sin{2\phi} - 2\cos{\theta}\cos{2\phi} \qquad && \text{in } \Omega  \\
&u=r^2\cos{\theta}\sin{2\phi} \qquad && \text{on } \partial\Omega, 
\end{alignat}
\end{subequations}
where $\Omega$ is the interior of a sphere of radius 0.25, centered at the origin. The analytic solution is given by $u=r^2\cos{\theta}\sin{2\phi}$. Our computational domain here is the periodic box $[-0.5, 0.5]^3$. }

\cha{We use the maximum determinant pionts \cite{spherepointspaper, spherepoints} for the discretization of the sphere and the associated quadrature weights for $dA_i$ in the generalized form of the spread operation in Equation \eqref{spread operator again2}. We use the analytical unit normal vectors to the sphere, and for the interpolation step, we do a linear interpolation using a known boundary value and an approximate solution value $m_2$ meshwidths away from the boundary in a direct generalization of the interpolation for the 2-D problem presented in \ref{interpolation appendix}.}

\begin{figure}
\centering
\includegraphics[width=0.5\textwidth]{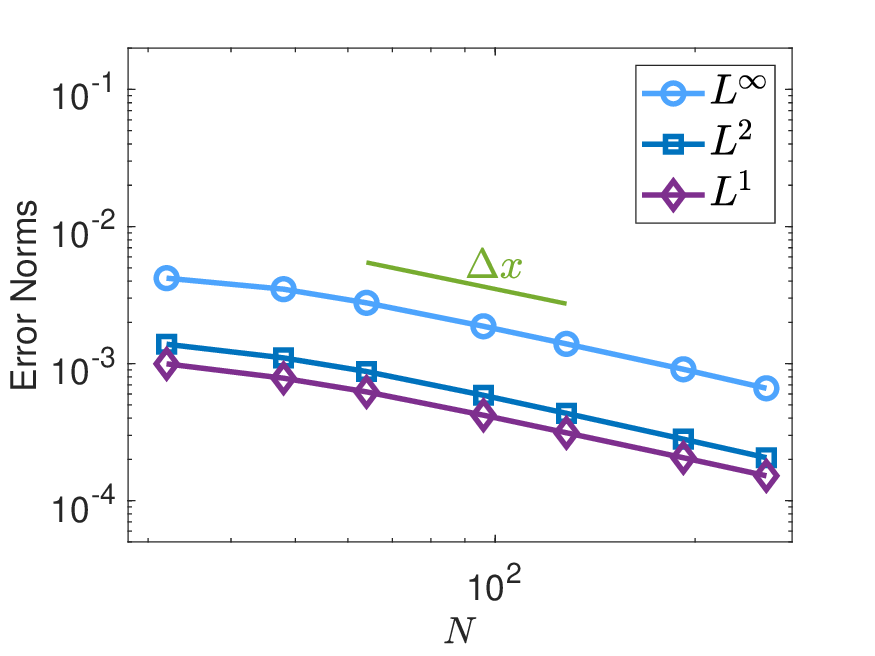}
\caption{\chIII{Refinement study to Equation \eqref{3d pde} found using the IBDL method. The computational domain is the periodic box $[-0.5, 0.5]^3$, $\Omega$ is the interior of a sphere of radius 0.25, and a Fourier spectral method is used to discretize the PDE. The boundary point spacing is $\Delta s  \approx 2 \Delta x$, and for the interpolation step, we replace solution values within $m_1= 2(\log_2{Nx}-4)$ meshwidths from the boundary using an interior interpolation point $m_2=m_1+2$ meshwidths away from the boundary. }}\label{3drefine}
\end{figure}

\cha{Table \ref{3d table} gives the iteration counts for the Krylov methods used to solve Equation \eqref{3d pde} for the IBSL and IBDL methods. For the coarsest mesh, we choose the number of points on the sphere so that in a triangulation of these points, the mean spacing between points is $\Delta s \approx c \Delta x$, for $c=2, 1.5, \text{ and } 1$. To maintain $\Delta s \approx c \Delta x$ as we refine the mesh, we quadruple the number of points on the sphere each time the Eulerian mesh spacing is halved. We can compare these iteration counts to the those in the similar 2-D problem in Equation \eqref{pde again}. Looking at $\Delta x = 2^{-7}$ and $\Delta s \approx \Delta x$ in Table \ref{iteration table 2}, we see that the IBSL method required $691$ iterations and the IBDL method required $5$ iterations. For the 3-D problem, the IBSL method required $9853$ iterations and the IBDL method still required only $5$ iterations. Since each iteration is more expensive in 3-D, this higher number of iterations required by the IBSL method can make a big difference. Maintaining low iteration counts in 3-D is a key advantage of the IBDL method. Figure \ref{3drefine} shows the refinement study for the IBDL method, again showing first-order convergence.}

\begin{figure}
\begin{center}
 \begin{tabular}{|| c || c | c | c | c | c | c | c ||}
 \hline
 \multicolumn{7}{||c||}{Iteration Counts - 3-D} \\
 \hline
   &\multicolumn{2}{c|}{$\Delta s\approx 2 \Delta x$ }&\multicolumn{2}{c|}{$\Delta s\approx 1.5 \Delta x$ } &\multicolumn{2}{c|}{$\Delta s\approx \Delta x$ } \\
 \hline
 $\Delta x$ &\textcolor{blue}{IBSL} & \textcolor{cyan}{IBDL} &\textcolor{blue}{IBSL} & \textcolor{cyan}{IBDL} &\textcolor{blue}{IBSL} & \textcolor{cyan}{IBDL} \\
 \hline
$2^{-5} $&  \textcolor{blue}{51}& \textcolor{cyan}{6}   &   \textcolor{blue}{144}& \textcolor{cyan}{6}   &  \textcolor{blue}{2320}& \textcolor{cyan}{6} \\
$2^{-6 }$&  \textcolor{blue}{72}& \textcolor{cyan}{6}   &   \textcolor{blue}{230}& \textcolor{cyan}{6}   &   \textcolor{blue}{6793}& \textcolor{cyan}{6} \\
$2^{-7} $&  \textcolor{blue}{101}& \textcolor{cyan}{6}   &   \textcolor{blue}{333}& \textcolor{cyan}{6}   &  \textcolor{blue}{9853}& \textcolor{cyan}{5}  \\
$2^{-8 }$&  \textcolor{blue}{135}& \textcolor{cyan}{5}   &   \textcolor{blue}{454}& \textcolor{cyan}{6}   &   \textcolor{blue}{--}& \textcolor{cyan}{--}   \\
%$2^{-7}$&  \textcolor{blue}{100}& \textcolor{cyan}{4}   & \textcolor{blue}{ 1669}&\textcolor{cyan}{4}      &   \textcolor{blue}{--}& \textcolor{cyan}{--}&  \textcolor{blue}{--}& \textcolor{cyan}{--}   \\
%$ 2^{-8}$&  \textcolor{blue}{133}& \textcolor{cyan}{4}   & \textcolor{blue}{1740}&\textcolor{cyan}{4}     &   \textcolor{blue}{--}& \textcolor{cyan}{--}&  \textcolor{blue}{--}& \textcolor{cyan}{--}   \\
%$ 2^{-9}$ &  \textcolor{blue}{214}& \textcolor{cyan}{4}   & \textcolor{blue}{4028}&\textcolor{cyan}{4}    &   \textcolor{blue}{-- }& \textcolor{cyan}{--}&  \textcolor{blue}{--}& \textcolor{cyan}{--}   \\
%$ 2^{-10}$ &  \textcolor{blue}{273}& \textcolor{cyan}{4}   & \textcolor{blue}{4088}&\textcolor{cyan}{4}    &   \textcolor{blue}{-- }& \textcolor{cyan}{--}&  \textcolor{blue}{--}& \textcolor{cyan}{--}   \\
%$ 2^{-12}$  & \textcolor{blue}{142}&\textcolor{cyan}{4}  &  \textcolor{blue}{ 4535}& \textcolor{cyan}{4}   \\
  \hline
  \end{tabular}
\captionof{table}[Number of iterations of \texttt{minres} and \texttt{gmres} needed to solve the Helmholtz PDE using the IBSL and IBDL methods, respectively]{\footnotesize{Number of iterations of \texttt{minres} and \texttt{gmres}, with tolerance $10^{-8}$, needed to solve Equation \eqref{3d pde} using the IBSL and IBDL methods, respectively. The computational domain is the periodic box $[-0.5, 0.5]^3$, and $\Omega$ is the interior of a sphere of radius 0.25. For $\Delta s\approx 2 \Delta x$, the  chosen maximum determinant boundary meshes from \cite{spherepoints} have 225, 900, 3600, and 14400 points. For $\Delta s\approx 1.5 \Delta x$, the number of boundary points are 400, 1600, 6400, and 25600. For $\Delta s\approx \Delta x$, the number of boundary points are 900, 3600, and 14400.}} 
\label{3d table}
%\end{center}
\end{center}
\end{figure}

%%%%%%%%%%%%%%%%%%%%%%%%%%%%%%%%%%%%%%%%%%%%%%%%%%%%%%%%%%%%%%%%   DISCUSSION  %%%%%%%%%%%%%%%%%%%
%%%%%%%%%%%%%%%%%%%%%%%%%%%%%%%%%%%%%%%%%%%%%%
\section{Discussion}\label{discussion}

We have developed the Immersed Boundary Double Layer method, a numerical method for linear PDE on complex domains with Dirichlet or Neumann boundary conditions that is more efficient than the constraint formulation of the Immersed Boundary method. We achieved this greater efficiency by reformulating the IB constraint method to correspond to a regularized double layer integral equation, which has better conditioning than the single layer integral equation to which the constraint formulation corresponds. With this better conditioning, we can solve for the Lagrangian dipole force distribution $Q$, in a small number of iterations of a Krylov method, and the iteration count does not increase as we refine the mesh. Furthermore, the computed $Q$ is relatively smooth and converges reasonably to the exact potential strength. Both of these are in stark contrast to the IB constraint method, in which the Lagrangian delta force distribution, $F$, requires many iterations to compute, the number of which increases with finer meshes and tighter boundary point spacing. Additionally, this distribution often fails to converge due to high-frequency noise. 

The IBDL method retains much of the flexibility and robustness of the IB method. The communication between the Lagrangian coordinate system and the underlying Cartesian grid is achieved with convolutions with regularized delta functions, and no analytical Green's functions are needed in the computation. Furthermore, minimal geometric information is needed for the immersed boundary. We only need boundary points, unit normals, an indicator function flagging points as inside the PDE domain, $\Omega$, and an indicator function flagging points as near-boundary points. Furthermore, we can reduce this input list to just the set of \chIII{boundary} points. The unit normals can be estimated as described in Section \ref{discret}, and the interior indicator function can be calculated as described in \ref{in out appendix}. Since such minimal geometric information is needed, this method can be easily used for complex PDE domains. \chII{The IBDL method also maintains flexibility in the PDE discretization, which is a typical advantage of an IB method. In this paper, we have compared two discretizations, and while the width of the interpolation region is affected by the choice, the IBDL method works with either option.}

Unlike the IB constraint method, the IBDL method produces a discontinuous solution across the boundary, \chII{resulting in a loss of pointwise convergence in a region near the boundary, the width of which approaches $0$ as the mesh is refined. If we wish to recover pointwise convergence up to the boundary, we can replace solution values for near-boundary points. In this paper, we have utilized a simple linear interpolation to replace these values. However, there is much work in the realm of boundary integral methods to harness analytical information of the PDE to achieve accurate near-boundary evaluations \cite{BealeLai, QBX, ShilpaBI}. The convergence of the potential strength $Q$ opens the door to the possibility of utilizing similar methods in the IBDL context. }

The IBDL method is related to the method of immersed layers, presented in the recent work of Eldredge \cite{eldredge}. In this work, Eldredge uses the indicator function discussed in \ref{in out appendix} to develop extended forms of PDEs that govern variables that define different functions on either side of the boundary. Jump quantities corresponding to the strengths of single and double layer potentials naturally emerge in the PDEs. The goal of this method is to be able to enforce different constraints on either side of the immersed boundary in order to obtain solutions on both sides and to accurately predict the surface traction on one side of the boundary. With the inclusion of jumps in both the solution and derivative, this method is able to achieve these goals. We should also note here that our application of the IBDL method to a Neumann problem includes a known delta force distribution and was derived from an integral representation containing both single and double layer potentials. Therefore, it can be seen as a specific case of the method of immersed layers \cite{eldredge}. However, if we focus on Dirichlet problems in which the solution is only desired on one side of the immersed boundary, the method of immersed layers again requires inverting an operator corresponding to a first-kind integral equation that is poorly conditioned. 

The IBDL method for Dirichlet boundary conditions uses a double layer potential alone, allowing a jump in only the solution values. This allows us to develop a method that corresponds to a second-kind integral equation that can be solved much more efficiently. One can see the IB constraint method of Taira and Colonius \cite{Taira} as corresponding to a single layer integral representation, the IBDL method presented in this paper as corresponding to a double layer integral representation, and the method of immersed layers of Eldredge \cite{eldredge} as corresponding to an integral representation that includes both single and double layers. And in the case of Dirichlet boundary conditions, the better conditioning of the IBDL operator gives a more efficient method. 

We also make note of another difference between the IBDL method and the method of immersed layers. While \cite{eldredge} uses a lattice Green's function for unbounded external flows, the IBDL method does not involve an explicit Green's function. Instead the method obtains the convolution of the boundary potential with a regularized Green's function by solving the PDE after spreading the boundary potential with the regularized delta function. Thus the IBDL method can be used on more general domains. In this paper we use a periodic computational domain, but \chIII{we} could also easily use other boundary conditions and domains.

The increased efficiency of this method makes it practical for use on time-dependent PDEs. For example, with an implicit time-discretization of the diffusion equation, \chIII{we} could use the IBDL method to solve a Helmholtz equation at each time step. \chIII{We} can also use the IBDL method for a PDE with a nonlinearity by using an implicit-explicit time-stepping scheme, as is done for the IBSE method in Stein et al. \cite{Stein}. 

\chII{In this paper, we have introduced the IBDL methodology through its application to the Helmholtz and Poisson equations, as a first step towards applying the method to fluid problems involving rigid bodies. Formulating the IBDL method for Stokes and Navier-Stokes equations introduces nontrivial challenges. First, an arbitrary external flow cannot be represented by a double layer potential alone \cite{Pozred}. However, the completed IBDL method required for large exterior domains for the periodic Poisson equation in this work will be utilized to overcome this obstacle. Second, the incompressibility condition introduces challenges for spatial discretization. We have obtained preliminary results applying the IBDL method to Stokes and Navier-Stokes equations, and this work will be presented in a forthcoming paper.}

\section*{\ch{Acknowledgment}}

This work was supported in part by NSF grant DMS-1664679 to R.D.G. 

\appendix

%%%%%%%%%%%%%%%%%%%%%%%%%%%%%%%%%%%%%%%%%%%%%%%%%%%%%%%%%%%%%%%%   interpolation appendix  %%%%%%%%%%%%%%%%%%%
%%%%%%%%%%%%%%%%%%%%%%%%%%%%%%%%%%%%%%%%%%%%%%
\section{Interpolation for solution values near the boundary}\label{interpolation appendix}

\begin{figure}
\centering
\includegraphics[width=0.69\textwidth]{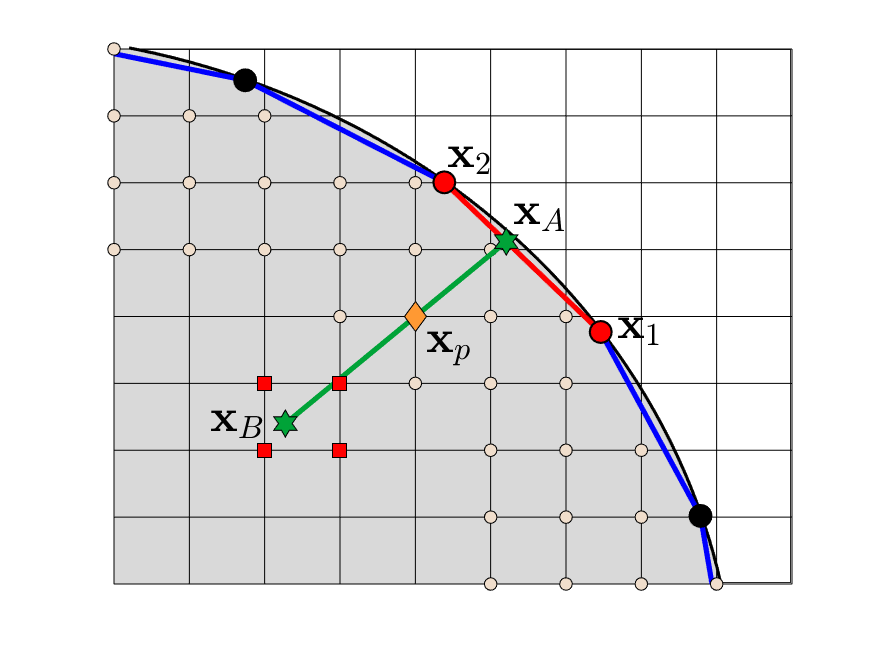}
\caption{Illustration of the interpolation step of the IBDL method. The PDE domain $\Omega$ is shown in grey. The exact boundary curve is shown in black with black circles indicating immersed boundary points. Lines connecting the IB points are shown in blue. Grid points within $m_1=3$ meshwidths of the boundary are indicated by beige circles. The diagram illustrates the process for interpolating the value for one point, denoted with an orange diamond and labeled $\x_p$. The two closest boundary points are marked with red circles and labeled $\x_1$ and $\x_2$. The closest point on the line segment connecting $\x_1$ and $\x_2$ is called $\x_A$ and is marked with a green pentagram. Also marked with a green pentagram is the point $m_2=4$ meshwidths away from $\x_A$, and it is labeled $\x_B$. These are the two interpolation points for $\x_p$.   }\label{interp pic}
\end{figure}

In this section, we discuss the details of the linear interpolation used to replace solution values within $m_1$ meshwidths of the boundary. Figure \ref{interp pic} illustrates the interpolation for one point, $\x_p$. Given such a point, we first locate the two nearest boundary points, labeled in the picture as $\x_1$ and $\x_2$. We use a projection to find $\x_A$, the point on the segment between $\x_1$ and $\x_2$ that is closest to $\x_p$. It is given by 
\begin{equation}
\x_A=\x_2+\frac{(\x_2-\x_p)\dotp (\x_1-\x_2)}{||\x_1-\x_2||^2}(\x_1-\x_2). 
\end{equation}
If the curvature and point placements are such that this calculation results in a point not located on the line segment between $\x_1$ and $\x_2$, we simply find the third closest boundary point and use the two outermost points to find $\x_A$. 
Next, we use a simple interpolation of the values at $\x_1$ and $\x_2$ to approximate the value at $\x_A$, and this is given by
\begin{equation}
u(\x_A)\approx U_b(\x_1)\frac{||\x_A-\x_2||}{||\x_2-\x_1||} + U_b(\x_2)\frac{||\x_A-\x_1||}{||\x_2-\x_1||} ,
\end{equation}
where the values at $\x_1$ and $\x_2$ are denoted with $U_b$ because these are known boundary values.

We next find the second point of interpolation, located $m_2$ meshwidths away from point $\x_A$. It is labeled in the figure as $\x_B$, and it is given by 
\begin{equation}
\x_B=\x_A+\frac{\x_p-\x_A}{||\x_p-\x_A||}(m_2\Delta x). \label{find B}
\end{equation}
Note that this step assumes that points are labeled as interior or exterior according to the \textit{discretized} shape created by the IB points. One way of doing this is presented in \ref{in out appendix}. Otherwise, \chIII{we} would need to verify that $\x_B$ is located on the correct side of the boundary and if not, alter Equation \eqref{find B} to use subtraction instead of addition.

We then estimate the value for $\x_B$ with a bilinear interpolation of the computed solution values at the four grid points located around $\x_B$, marked in the figure with red squares. Let $\x_{00}$ be the lower left point, $\x_{01}$ be the upper left, and so on, and let $d_x$ and $d_y$ be the $x$ and $y$ distances from $\x_B$ to $\x_{00}$. Then our bilinear interpolation is 
\begin{multline}
u(\x_B)\approx \Big(\frac{d_x}{\Delta x}\Big) \Big(\frac{d_y}{\Delta y}\Big) u(\x_{11}) + \Big(\frac{d_y}{\Delta y}\Big)\Big(1-\frac{d_x}{\Delta x}\Big)u(\x_{01})\\+\Big(\frac{d_x}{\Delta x}\Big)\Big(1-\frac{d_y}{\Delta y}\Big)u(\x_{10}) + \Big(1-\frac{d_x}{\Delta x}-\frac{d_y}{\Delta y}+ \Big(\frac{d_x}{\Delta x}\Big) \Big(\frac{d_y}{\Delta y}\Big) \Big)u(\x_{00})
\end{multline}
Note that if $\x_B$ is located directly on a gridline, this can be \ch{simplified}. 

Once we have the approximated values at $\x_B$ and $\x_A$, we again use a simple linear interpolation to approximate $u(\x_p)$, and this is given by
\begin{equation}
u(\x_p)\approx u(\x_A)\frac{||\x_B-\x_p||}{||\x_B-\x_A||} + u(\x_B)\frac{||\x_A-\x_p||}{||\x_B-\x_A||} .
\end{equation}

%%%%%%%%%%%%%%%%%%%%%%%%%%%%%%%%%%%%%%%%%%%%%%%%%%%%%%%%%%%%%%%%   neumann appendix  %%%%%%%%%%%%%%%%%%%
%%%%%%%%%%%%%%%%%%%%%%%%%%%%%%%%%%%%%%%%%%%%%%
\section{Immersed Boundary formulation of the Neumann Helmholtz problem}\label{neumann appendix}

As stated in Section \ref{Neumann formulation}, the IB formulation of the Neumann problem is given by 
\begin{subequations} \label{ib neumann again}
\begin{alignat}{2}
& \L u +\widetilde S U_b+SV_b = g \qquad && \text{in } \mathcal{C}  \label{ib neumann again 1}\\
&S^* u = \frac12 U_b \qquad && \text{on } \Gamma,   \label{ib neumann again 2}
\end{alignat}
\end{subequations}
where $V_n\equiv \partial u/\partial n |_{\Gamma}$ gives the known normal derivatives on $\Gamma$, and $U_b$ gives the \chI{unknown} distribution of solution values on $\Gamma$. $V_b$ corresponds to the strength of a single layer potential and $U_b$ to the strength of a double layer potential. 

Starting with Equation \eqref{ib neumann again 1} , we have 
\begin{equation}
\L u = - \widetilde S U_b - S V_b= - \grad \dotp \int_{\Gamma} U_b(s)\n(s) \delta_h(\x-\X(s)) ds  - \int_{\Gamma} V_b(s) \delta_h(\x-\X(s)) ds. \label{neum connect 1}
\end{equation}
Inverting the operator and using that $\L G_h(\x,\x_0)=-\delta_h(\x-\x_0)$, we get
\begin{equation}
u(\x) = \grad \dotp  \int_{\Gamma} U_b(s)\n(s)G_h(\x,\X(s)) ds  + \int_{\Gamma} V_b(s) G_h(\x,\X(s)) ds \label{neum connect 2}
\end{equation}
for $\x\in\Omega\setminus \Gamma$. Bringing in the divergence and manipulating the expression, we get
\begin{equation}
u(\x) =   \int_{\Gamma} U_b(s) \grad G_h(\x,\X(s)) \dotp\n(s) ds+ \int_{\Gamma} V_b(s) G_h(\x,\X(s)) ds \label{neum connect 3}
\end{equation}
The second equation of the method, Equation \eqref{ib neumann again 2} gives us 
\begin{equation}
\frac12 U_b(s') =\int_{\C} u(\x)\delta_h(\x-\X(s')) dx. \label{neum connect 4}
\end{equation}
Combining Equation \eqref{neum connect 3} with Equation \eqref{neum connect 4}, changing the order of integration, and recognizing the presence of $G_{hh} = G_h*\delta_h$, we get 
\begin{equation}
\frac12 U_b(s') =\int_{\Gamma} U_b(s) \grad G_{hh}(\X(s'),\X(s)) \dotp\n(s) ds +\int_{\Gamma} V_b(s) G_{hh}(\X(s'),\X(s)) ds. \label{neum connect 5}
\end{equation}
Note that due to the presence of a double layer potential, as discussed in Section \ref{formulation}, we assume that $\X(s)$ is the arclength parametrization of $\Gamma$ so that $  \big|\partial \X(s)/\partial s\big| =1$.

The symmetry of the Green's function, which is preserved through convolutions with the regularized delta function, gives us that $G_{hh}(\X(s'),\X(s)) = G_{hh}(\X(s),\X(s'))$. Additionally, the odd symmetry of the gradient of the Green's function gives us  $ \grad G_{hh}(\X(s'),\X(s))=- \grad G_{hh}(\X(s),\X(s'))$ \cite{Pozblue}. Using these to switch the arguments of these functions, we get
\begin{equation}
\frac12 U_b(s') =-\int_{\Gamma} U_b(s) \grad G_{hh}(\X(s),\X(s')) \dotp\n(s) ds +\int_{\Gamma} V_b(s) G_{hh}(\X(s),\X(s')) ds. \label{neum connect 6}
\end{equation}
Appropriately redefining $U_b$ and $V_b$ as functions of $\x$, we get 
\begin{multline}
\frac12 U_b(\X(s')) =-\int_{\Gamma} U_b(\X(s)) \grad G_{hh}(\X(s'),\X(s)) \dotp\n(s) ds \\
+ \int_{\Gamma} V_b(\X(s)) G_{hh}(\X(s'),\X(s)) ds. \label{neum connect 7}
\end{multline}

Recall that the integral equation for $\X_0$ on the boundary is given by 
\begin{multline}
u(\X_0) = \int_{\Gamma}^{PV} \grad u(\x) \dotp \n (\x)G(\x, \X_0)  dl(\x) \\ -\int_{\Gamma} u(\x) \grad G(\x,\X_0) \dotp \n(\x) dl(\x) +
\frac{1}{2} u(\X_0). \label{integral rep again appendix}
\end{multline} 
Using our arclength parametrization and combining the $u(\X_0)$ terms, we can rewrite this as
\begin{multline}
\frac12 u(\X(s')) =- \int_{\Gamma}^{PV} u(\X(s)) \grad G(\X(s),\X(s')) \dotp \n(\X(s)) ds\\
+ \int_{\Gamma} G(\X(s), \X(s')) \grad u(\X(s)) \dotp \n (\X(s)) ds . \label{neum connect 8}
\end{multline}
By comparing Equations \eqref{neum connect 7} and \eqref{neum connect 8}, we can see that the IB formulation presented in Equation \eqref{ib neumann again} is equivalent to a \emph{regularized} integral equation.

%%%%%%%%%%%%%%%%%%%%%%%%%%%%%%%%%%%%%%%%%%%%%%%%%%%%%%%%%%%%%%%%%%%%%%%%%%%%%%%%%%%%%%%%%%%%%%%%% NULLSPACE APP

\section{\chIII{Nullspace of periodic Laplacian}}\label{nullspace appendix}

\chIII{In the case that the differential operator is the periodic Laplacian, the original boundary value problem on $\Omega$ has a unique solution, but, by using the IB framework, we embed the PDE into a periodic computational domain on which $\Delta$ is not invertible. To recover the unique solution to the PDE, we follow a method similar to that presented by Stein et al.\ \cite{Stein}.  }

\chIII{We begin by decomposing  the solution $u$ into 
\begin{equation}
u=u_0+\bar u,
\end{equation}
where $u_0$ has mean $0$ on $\C$ and $\bar u$ gives the mean value of $u$ on $\C$. The PDE on $\C$ that results from the completed IBDL formulation in Equation \eqref{ibdl 1 completed} is given by 
\begin{equation}
\Delta u_0 + \eta SQ +\grad\dotp (SQ\n)  =\tilde g. \label{laplaces with ibdl}
\end{equation}
To derive the solvability condition, we integrate Equation \eqref{laplaces with ibdl} over a general computational domain $\C$. Noting that our spread operator is defined using a regularized delta function, the functions are smooth, and we can use the divergence theorem. Let us use $\partial \C$ as the boundary of $\C$ and $\n_{\C}$ as the unit normal on $\partial \C$, to distinguish it from $\n$, which we continue to use as the unit normal on the immersed boundary $\Gamma$. We then get
\begin{equation}
\int_{\partial \C} \grad u_0 \dotp \n_{\C} dl(\x)  +\eta \int_{\C} SQ d\x+ \int_{\partial \C} S Q \n\dotp \n_{\C} dl(\x) = \int_{\C} \tilde g d\x.  \label{first step after divergence theorem again}
\end{equation}
If $\Gamma$ is away from $\partial \C$, the third term vanishes due to the compact support of the integrand. Then, if we take $\C$ to be the periodic box used in this work, the first term also disappears. Using that $\int_{\C} SQ d\x=\int_{\Gamma}Qds$, the solvability condition on the periodic computational domain $\C$ then has the form
\begin{equation}
\eta \int_{\Gamma} Qds=\int_{\C}\tilde g(\x) d\x .\label{poisson completed solv}
\end{equation}
When we are not using the completed formulation, since $\eta=0$, we can easily satisfy this condition with our choice of $ g_e$. For example, we can choose 
\begin{equation}
g_e(\x) =-\frac{1}{|\C\setminus \Omega|} \int_{\Omega} g(\x)d\x, 
\end{equation} 
where we use $|\dotp | $ to denote the area of the domain. With this constraint satisfied, the solution to 
\begin{subequations}
\begin{alignat}{2}
&\Delta u +\grad\dotp (SQ\n)  =\tilde g\\
& S^* u +\frac12Q=U_b
\end{alignat}
\end{subequations}
 can be found, and it is unique up to an additive constant. Therefore, since we are not interested in the solution on the nonphysical domain, $\C\setminus \Omega$, we can use the solution with mean $0$ on $\C$. As discussed in Section \ref{nullspace1}, any other solution would give the same unique solution on the PDE domain $\Omega$.  }

\chIII{If, on the other hand, we use the completed IBDL method, we can use the constraint in Equation \eqref{poisson completed solv} to find the value of $\bar u$. Discretizing the constraint, we can then summarize the discrete IBDL system as 
\begin{subequations} \label{ibdl completed again}
\begin{alignat}{2}
& \Delta u_0 +\eta S Q+ \widetilde S Q = \tilde g \label{ibdlcomleted1a}\\
&S^* u_0 + \bar u \mathds{1}_{N_{IB}}+ \frac12 Q = U_b\label{ibdlcomleted1b}\\
& \eta(\Delta s) \mathds{1}^\intercal_{N_{IB}}Q = (\Delta x\Delta y)\mathds{1}^\intercal_{N^2} \tilde g,
\end{alignat}
\end{subequations}
where $\mathds{1}_{N^2}$ denotes a vector of length $N^2$ consisting of all ones. With the solvability constraint satisfied, we can then find the unique solution to Equation \eqref{ibdl completed again}. Letting $\Delta_0^{-1}$ denote the operation that inverts the Laplacian by returning a solution with mean $0$ on $\C$, we have $\Delta_0^{-1}\Delta u_0=u_0$. Using this to invert the Laplacian in  Equation \eqref{ibdlcomleted1a} and then applying $S^*$, we get
\begin{equation}
S^* u_0 +\eta S^* \Delta_0^{-1}SQ+S^* \Delta_0^{-1}\widetilde S Q = S^* \Delta_0^{-1} \tilde g.
\end{equation}
Then, we can replace $S^*u_0$ using Equation \eqref{ibdlcomleted1b}, and we get the the following system of discrete equations, for which we use \texttt{gmres} to solve for $Q$ and $\bar u$:
\begin{subequations}
\begin{alignat}{2}
&-(S^* \Delta_0^{-1}\widetilde S)Q -\eta(S^*\Delta_0^{-1}S)Q+\bar u \mathds{1}_{N_{IB}} +\frac12 Q=U_b - S^*\Delta_0^{-1} \tilde g\\
&\eta(\Delta s) \mathds{1}^\intercal_{N_{IB}}Q = (\Delta x\Delta y)\mathds{1}^\intercal_{N^2} \tilde g
\end{alignat}
\end{subequations}
We can then obtain $u$ by computing
\begin{equation}
u=-\Delta_0^{-1}(\widetilde S +\eta S)Q +\Delta_0^{-1}\tilde g+\bar u \mathds{1}_{N^2}.
\end{equation}}

%%%%%%%%%%%%%%%%%%%%%%%%%%%%%%%%%%%%%%%%%%%%%%%%%%%%%%%%%%%%%%%%   flag appendix  %%%%%%%%%%%%%%%%%%%
%%%%%%%%%%%%%%%%%%%%%%%%%%%%%%%%%%%%%%%%%%%%%%
\section{Using IBDL method to flag interior grid points}\label{in out appendix}

In this section, we show that by using analytical features of the indicator function for the interior of the immersed boundary $f\Gamma$, we arrive at a familiar IBDL expression that enables us to use the framework already set forth in order to flag grid points as either exterior or interior to $\Gamma$. It is useful to have the points flagged as such in order to identify the physical domain of the solution, as well as for completing the interpolation step discussed in \ref{interpolation appendix}. Eldredge also discusses this in Appendix A.2 of \cite{eldredge}, and his work uses such an indicator function as a foundation on which to build PDEs that govern different variables on different sides of the immersed boundary. Here, we present the derivation and application of this mechanism in the context of the IBDL method and also discuss how to apply it to a periodic computational domain. 

Let the indicator function be defined as
\begin{equation}
\chi(\x)=\begin{cases}
      1 & \x \text{ in } \Omega  \\
      0 & \x \text{ in } \C\setminus\Omega, \\
\end{cases}\label{indic}
\end{equation}
where we assume here that $\Omega$ is the interior of $\Gamma$. Then, using a test function, $\psi$, we have
\begin{equation}
\int_{\C} \chi \Delta \psi d\x = \int_{\Omega} \Delta \psi d\x= \int_{\Gamma} \grad \psi \dotp \n dl(\x) , \label{B2}
 \end{equation}
where we first use the definition of $\chi$ and then the divergence theorem. Then, since the weak derivative of $\chi$ is $\grad \chi$ such that 
\begin{equation}
\int_{\C} \grad \chi \dotp \grad \psi d\x = - \int_{\C} \chi \Delta \psi d\x , 
\end{equation}
we see from the last term in Equation \eqref{B2} that the weak derivative is given by 
\begin{equation}
\grad \chi = - \delta (\x-\X) \n(\x), 
\end{equation}
where the $\X$ is a point on $\Gamma$, indicating that the support of $\grad \chi$ is $\Gamma$. 

We therefore see that 
\begin{equation}
- \Delta \chi = -\grad \dotp \grad \chi = \grad \dotp (\delta \n) . 
\end{equation} 
We can recognize $\grad \dotp (\delta \n)$ as the continuous version of the term found in the IBDL method, $\widetilde S Q = \grad \dotp( S Q\n)$, where $Q$ is a constant distribution of $1$. This then gives us a way to construct the indicator function $\chi$ by solving the PDE on $\C$ given by
\begin{equation}
\Delta u + \widetilde S \n=0,
\end{equation} 
using the method already in place for the IBDL method. 

Then, if \chIII{we are} using a periodic computational domain, as we do in this paper, $u$ is only unique up to an additive constant. Therefore, we can first find the solution with mean $0$ on $\C$, as discussed in \ref{nullspace appendix}. Then, to shift the exterior values to $0$, we simply need to subtract the appropriate constant, which can be approximated using a value of $u$ far from $\Gamma$. A likely choice is the corner of the periodic domain. Subtracting this value from $u$ then results in a regularized indicator function, $\chi_h$, that is approximately equal to $1/2$ on the boundary $\Gamma$. To retrieve the final indicator function, we simply set all values of $\chi_h$ that are less than $1/2$ equal to $0$ and all those above $1/2$ to 1. Then $\chi$ is the indicator function described in Equation \eqref{indic} that is equal to 1 for all points on the interior of the discretized boundary $\Gamma$.

\newpage
%REFERENCES

\bibliography{ibdl_bibfile.bib}

\end{document}